\documentclass[a4paper]{article}
\usepackage[latin1]{inputenc}
\usepackage[T1]{fontenc}
\usepackage{amssymb}
\usepackage[dvips]{graphics}
\usepackage{hyperref}

\def\vell{{\vec\ell}}
\def\dx{{\dot x}}
\def\dy{{\dot y}}
\def\EK{{E\!/\!K}}

\def\RT{{\Im(T)}}
\def\RTp{{\Im(T')}}
\def\be{\begin{equation}}
\def\calL{{\cal L}}
\def\calT{{\cal T}}
\def\dd{(\cdot,\cdot)}
\def\demi{{1\over 2}}
\def\diag{{\rm diag}}
\def\DO{{\cal D}(\Omega)}
\def\dvg{{\rm div}}
\def\ee{\end{equation}}
\def\eqdef{{\;\mathop{=}^{\hbox{\footnotesize def}}\;}}
\def\eqnote{{\;\mathop{=}^{\hbox{\footnotesize denoted}}\;}}
\def\eref#1{(\ref{#1})}
\def\Es{{E'}}
\def\Ess{{E''}}

\def\grad{{\rm grad}}
\def\Hdemig{{H^{{1\over2}}(\Gamma)}}
\def\Hmuo{{H^{-1}(\Omega)}}
\def\Hu{{H^1}}

\def\Huo{{H^1(\Omega)}}
\def\Huz{{H^1_0}}
\def\Huzo{{H^1_0(\Omega)}}
\def\Im{{\rm Im}}
\def\Ker{{\rm Ker}}
\def\Ld{{L^2}}
\def\Ldg{{L^2(\Gamma)}}
\def\Ldo{{L^2(\Omega)}}

\def\Ldzo{{L^2_0(\Omega)}}
\def\matrarrow{\mathop{\longrightarrow}}
\def\mrar{\mathop{\longrightarrow}}
\def\NN{\mathbb{N}}
\def\NNs{\mathbb{N}^*}
\def\PVh{{\Pi_{V_h}\!}}
\def\PQh{{\Pi_{Q_h}\!}}
\def\RR{\mathbb{R}}
\def\RRn{{\mathbb{R}^n}}

\def\ve{{\vec{e}}}
\def\Vect{{\rm Vect}}
\def\vu{{\vec{u}}}
\def\vv{{\vec{v}}}
\def\vw{{\vec{w}}}
\let\ds\displaystyle
\let\eps\varepsilon
\let\la\langle
\let\pa\partial
\let\ra\rangle
\let\rar\rightarrow
\long\def\comment#1{}
\def\Qh{{Q_h}}
\def\Vh{{V_h}}
\def\Xh{{X_h}}
\catcode`@=11
\def\@afterheading{\global\@nobreaktrue
      \everypar{\if@nobreak
                   \global\@nobreakfalse
                   \clubpenalty \@M
                 \else \clubpenalty \@clubpenalty
                    \everypar{}\fi}}

\@addtoreset{equation}{section} 
\@addtoreset{figure}{section} 
\def\eqalign#1{\null\,\vcenter{\openup1\jot \m@th
   \ialign{\strut \hfil$\displaystyle{##}$ & $\displaystyle{{}##}$\hfil
      \crcr#1\crcr}}\,}
\def\eqalignrll#1{\null\,\vcenter{\openup1\jot \m@th
   \ialign{\strut \hfil$\displaystyle{##}$ & $\displaystyle{{}##}$\hfil
    & $\displaystyle{{}##}$\hfil
      \crcr#1\crcr}}\,}
\def\eqalignrllq#1{\null\,\vcenter{\openup1\jot \m@th
   \ialign{\strut \hfil$\displaystyle{##}$ & $\displaystyle{{}##}$\hfil
    & \qquad$\displaystyle{{}##}$\hfil
      \crcr#1\crcr}}\,}

\renewcommand\tableofcontents{\section*{\contentsname}\@starttoc{toc}}

\catcode`@=12
\def\urltilde{\raise.3ex\hbox{\footnotesize$\!\sim$}}

\newtheorem{lemme}{Lemma}[section]  
\newtheorem{corollary}[lemme]{Corollary}
\newtheorem{definition}[lemme]{Definition}
\newtheorem{example}[lemme]{Example}
\newtheorem{exercise}[lemme]{Exercise}
\newtheorem{prop}[lemme]{Proposition}  
\newtheorem{remark}[lemme]{Remark}
\newtheorem{theorem}[lemme]{Theorem}

   \def\boxending{\leavevmode\nolinebreak\hfill{$\tblackbox$}}

\def\debcor{\begin{corollary}\sl }
   \def\fincor{\end{corollary}}
\def\debdef{\begin{definition}\rm }
   \def\findef{\end{definition}}
\def\debexa{\begin{example}\rm}
   \def\finexa{\boxending\end{example}}
\def\debexe{\begin{exercise}\rm}
   \def\finexe{\boxending\end{exercise}}
\def\deblem{\begin{lemme}\sl }
   \def\finlem{\end{lemme}}
\def\debprop{\begin{prop}\sl }
   \def\finprop{\end{prop}}
\def\debrem{\begin{remark}\rm}
   \def\finrem{\boxending\end{remark}}
\def\debthm{\begin{theorem}\sl }
   \def\finthm{\end{theorem}}

\def\debproof{\medskip\par\noindent{\bf Proof.\ }}
   \def\finproof{\boxending\par\medskip}
\let\debdem\debproof
   \let\findem\finproof
\def\tbox{\leavevmode\vrule height 3pt width 3pt depth 0pt\relax}
\def\tblackbox{\tbox\kern-.75pt\raise4pt\hbox{\tbox}\kern-.75pt\tbox}

\def\cf{cf{.}}
\def\Eg{E{.}g{.}}
\def\eg{e{.}g{.}}

\def\ie{i{.}e{.}}
\def\st{s{.}t{.}}

\def\qwith{\quad\hbox{with}\quad}
\def\qand{\quad\hbox{and}\quad}
\def\qandthen{\quad\hbox{and then}\quad}

\def\qe{\quad\hbox{et}\quad}

\def\qie{\quad\hbox{\ie}\quad}

\def\qwhere{\quad\hbox{where}\quad}

\def\qtq{\quad\hbox{t.q.}\quad}
\def\qst{\quad\hbox{\st}\quad}
\def\qwhen{\quad\hbox{when}\quad}

\def\sumin{{\sum_{i=1}^n}}

\def\vcurl{{\mathop{\vec{\rm curl}}}}
\def\curl{{\mathop{\rm curl}}}

\def\calF{{\cal F}}
\def\calI{{\cal I}}

\def\Hd{{H^2}}
\def\Hdo{{H^2(\Omega)}}
\def\Hdzo{{H^2_0(\Omega)}}
\def\Hdz{{H^2_0}}
\def\Hdemi{{H^{\demi}}}

\def\Hmdemig{{H^{-\demi}(\Gamma)}}
\def\Hmu{{H^{-1}}}
\def\Hmd{{H^{-2}}}
\def\Hmuo{{H^{-1}(\Omega)}}
\def\Hmuon{{\Hmuo^n}}
\def\Huon{{\Huo^n}}
\def\Huzon{{\Huzo^n}}
\def\Ldon{{\Ldo^n}}
\def\ld{{\ell^2}}
\def\NNs{{\NN^*}}
\def\rang{{\rm Rang}}
\def\RRm{{\RR^m}}

\def\va{{\vec a}}

\def\vf{{\vec f}}
\def\vg{{\vec g}}
\def\vgrad{{\vec \grad}}
\def\vlambda{{\vec \lambda}}
\def\vn{{\vec n}}
\def\vphi{{\vec\phi}}

\def\vu{{\vec u}}
\def\vw{{\vec w}}
\def\vx{{\vec x}}
\def\vy{{\vec y}}
\def\vz{{\vec z}}

\def\vgamma{{\vec\gamma}}
\def\vmu{{\vec\mu}}
\def\Ldz{{L^2_0}}
\def\Hmdo{{H^{-2}(\Omega)}}

\def\Hdvgo{{H^\dvg(\Omega)}}   \def\Hdvg{{H^\dvg}}
\def\Hdvgzo{{H^\dvg_0(\Omega)}}   \def\Hdvgz{{H^\dvg_0}}

\def\Hcurlo{{H^\curl(\Omega)}}   \def\Hcurl{{H^\curl}}
\def\Hcurlzo{{H^\curl_0(\Omega)}}   \def\Hcurlz{{H^\curl_0}}

\def\anneecours#1{\ifcase #1\or première\or deuxième\or troisième\fi}

\title{\vspace{-2cm}
Inf-sup condition and locking: Understanding and circumventing.
\\ \smallskip
\Large %Examples : 
Stokes, Laplacian, bi-Laplacian, Kirchhoff--Love and Mindlin--Reissner locking type, boundary conditions}
\date{April 22, 2020}
\author{Gilles Leborgne\\ Isima, University of Clermont-Ferrand, France}

\begin{document}
%\montitre{3}

\maketitle

\begin{abstract}
The inf-sup condition, also called the Ladyzhenskaya--Babu\v ska--Brezzi (LBB) condition, ensures the well-posedness of a saddle point problem, relative to a partial differential equation.
Discretization by the finite element method gives the discrete problem which must satisfy the discrete inf-sup condition. But, depending on the choice of finite elements, the discrete condition may fail. This paper attempts to explain why it fails from an engineer's perspective, and reviews current methods to work around this failure. The last part recalls the mathematical bases.
\end{abstract}

\tableofcontents

%%%%%%%%%%%%%%%%%%%%%%%%%%%%%%%%%%%%%%%%%%%%%%%%%%%%%%%%%%%%%%%
%%%%%%%%%%%%%%%%%%%%%%%%%%%%%%%%%%%%%%%%%%%%%%%%%%%%%%%%%%%%%%%

%%%%%%%%%%%%%%%%%%%%%%%%%%%%%%%%%%%%%%%%%%%%%%%%%%%%%%%%%%%%%%%%%%%%%%%%%%%%%%%%%%%

\newpage

\part{Introduction}
%%%%%%%%%%%%%%%%%%%%%%%%%%%%%%%%%%%%%%%%%%%%%%%%%%%%%%%%%%%%%%%%%%%%%%%%%%%%%%%%%%%

\section{The inf-sup condition, what is that ?}

%%%%%%%%%%%%%%%%%%%%%%%%%%%%%%%%%%%%%%%%%%%%%%%%%%%%%%%%%%%%%%%%%%%%%%%%%%%%%%%%%%%

\subsection{The constrained problem under concern}

Let $(V,\dd_V)$ be a Hilbert space and $(Q,||.||_Q)$ be a Banach space.
Let $V'=\calL(V;\RR)$ and $Q'=\calL(Q;\RR)$ be the associated dual spaces
(spaces of the linear continuous forms).
Let $f\in V'$ and $g\in Q'$ be linear continuous forms,
and let $a\dd : V \times V\rar\RR$ and $b\dd : V\times Q \rar \RR$ be bilinear continuous forms.

The problem under concern is: Find $(u,p)\in V \times Q$ \st\
\be
\label{eqgenw0}
\left\{\eqalignrll{
& a(u,v) + b(v,p) &= \la f,v \ra_{V',V}, \quad \forall v \in V , \cr
& b(u,q) &= \la g,q\ra_{Q',Q}, \quad \forall q \in Q.
}\right.
\ee
\Eg, with $\Omega$ an open regular bounded set in~$\RRn$
and $\Ldzo\simeq \Ldo/\RR$ (that is the space of $\Ld$ functions defined up to a constant), find $(\vu,p)\in \Huzo^n \times \Ldzo$ \st\
\be
\label{eqgenw1}
\left\{\eqalign{
& (\grad \vu , \grad\vv)_\Ld - (p,\dvg \vv)_\Ld = \la\vf,\vv\ra_{\Hmu,\Huz},\quad \forall \vv \in \Huzo^n \cr
& -(\dvg\vu,q)_\Ld = 0, \quad \forall q \in \Ldzo.
}\right.
\ee
%(Here $p$ is unique up to a constant, and $\Ldzo = \{ f \in \Ldo,\; \int_\Omega p\;d\Omega\} \simeq \Ldo/\RR$.)

Let $A\in \calL(V;V')$, $B\in\calL(V;Q')$ and $B^t\in\calL(Q,V')$
be the associated linear continuous mapping (bounded operators)
given by
\be
\label{eqDefB}
\la Au,v \ra_{V',V} = a(u,v), \quad \la Bv,p \ra_{Q',Q} = b(v,p) = \la B^t p,v\ra_{V',V}.
\ee
Then~\eref{eqgenw0} also reads
\be
\label{eqgens0}
\left\{\eqalignrll{
& Au + B^t p &= f \in V', \cr
& Bu &= g  \in Q',
}\right.
\ee
the equation $Bu=g$ being the constraint. \Eg, find $(\vu,p)\in \Huzo^n \times \Ldzo$ \st\
\be
\left\{\eqalignrll{
& {-}\Delta \vu  + \vgrad p &= \vf , \cr
& \dvg\vu &= 0.
}\right.
\ee

%%%%%%%%%%%%%%%%%%%%%%%%%%%%%%%%%%%%%%%%%%%%%%%%%%%%%%%%%%%%%%%%%%%%%%%%%%%%%%%%%%%

\subsection{The control on~$p$ to get a well-posed problem, and inf-sup}

The simplest numerical finite element simulations show non admissible results for~$p$ (the pressure) in~\eref{eqgenw1}. And $p$ being only present in~\eref{eqgenw1}$_1$,
we need to study~$B^t$, \cf~\eref{eqgens0}.
The needed result will be the closure of $\Im(B^t)$ in~$V'$: In that case the open mapping theorem gives
the control on~$p$ thanks to:
\be
\label{eq0is}
\exists \beta > 0 ,\;
\forall p\in Q,\quad  ||B^t p||_{V'} \ge \beta ||p||_{Q/\Ker(B^t)},
%\forall p\in Q,\quad   ||p||_{Q/\Ker(B^t)} \le {1\over \beta}||B^t p||_{V'},
\ee
also written as the ``inf-sup condition'': 
$\ds \exists \beta>0, \;
\inf_{p\in Q}\sup_{v\in V} {|b(v,p)| \over ||v||_{V} ||p||_{Q/\Ker(B^t)}}  \ge \beta$,
since $||B^t p||_{V'} = \sup_{v\in V} { |\la B^tp,v\ra_{V',V}| \over ||v||_V}
= \sup_{v\in V} { |b(v,p)| \over ||v||_V}$.

\Eg\ for~\eref{eqgenw1}, with $B^t = \vgrad : \Ldzo \rar \Hmuon$, we have $\Im(B^t)$ closed in~$\Hmuon$
(this is ``the'' difficult theorem to establish, see next~\S), thus
\be
\exists \beta > 0 ,\;
\forall p\in \Ldo,\; ||\vgrad p||_\Hmu \ge \beta ||p||_\Ldz.
\ee
Which can be written as the ``inf-sup condition'':
$\ds \exists \beta>0, \;
\inf_{p\in \Ldo}\sup_{v\in \Huzon} {|(\dvg\vv,p)_\Ld| \over ||v||_\Huz ||p||_\Ldz}  \ge \beta$.

\medskip
And we get the theorem: The problem~\eref{eqgenw1} is well-posed,
that is, the solution $(u,p)$ exists, is unique, and depends continuously on $f$ and~$g$. See~\eref{eqcis2}.
%(The hypothesis $a\dd$ coercive on~$\Ker B$ can be relaxed into two inf-sup condition on~$\Ker B$ stating that $A $ is invertible on~$\Ker B$.)

\medskip
Remark: Thanks to the closed range theorem~\ref{thmcr},
the closure of~$\Im(B^t)$ is equivalent to the closure of~$\Im(B)$
(under usual hypotheses). This result is needed to get the existence of~$u$.

%%%%%%%%%%%%%%%%%%%%%%%%%%%%%%%%%%%%%%%%%%%%%%%%%%%%%%%%%%%%%%%%%%%%%%%%%%%%%%%%%%%

\subsection{The loss of control on~$p$ for the discrete problem}

\def\Eh{{E_h}}
\def\Fh{{F_h}}
\def\PEh{{\Pi_{E_h}\!}}
\def\PFh{{\Pi_{F_h}\!}}
\def\Vh{{V_h}}
\def\Qh{{Q_h}}
\def\PVh{{\Pi_{V_h}\!}}
\def\PQh{{\Pi_{Q_h}\!}}

Let $\Vh\subset V$ and $\Qh \subset Q$ be finite dimensional subspaces (conform finite elements to simplify).
The discretization of~\eref{eqgenw0} (for computation purposes) reads: Find $(u_h,p_h)\in \Vh \times \Qh$ \st\
\be
\label{eqgenh0}
\left\{\eqalignrll{
& a(u_h,v_h) + b(v_h,p_h) &= \la f,v_h \ra_{V',V}, \quad \forall v_h \in \Vh , \cr
& b(u_h,q_h) &= \la g,q_h\ra_{Q',Q}, \quad \forall q_h \in \Qh.
}\right.
\ee
\Eg, with $\Vh \subset \Huzo^n$, $\Qh \subset \Ldzo$, and $\vf\in\Ldon$, 
\be
\label{eqgenh0eg}
\left\{\eqalignrll{
& (\grad \vu_h,\grad \vv_h)_\Ld - (p_h,\dvg\vv_h)_\Ld &= (\vf,\vv_h)_\Ld, \quad \forall \vv_h \in \Vh , \cr
& -(\dvg\vu_h,q_h)_\Ld &= 0, \qquad \forall q_h \in \Qh.
}\right.
\ee
The $h$-discrete inf-sup condition is
\be
\label{eqdish}
\exists \beta_h>0,\; \forall p_h\in Q_h,\quad ||B_h^t p_h||_{V'} \ge \beta_h ||p_h||_{Q_h/\Ker(B^t)}.
\ee
And $\beta_h$ should satisfy $\beta_h > \gamma$ is satisfied for some $\gamma>0$:
we get the so-called discrete inf-sup condition:
\be
\label{eqdis}
\exists \gamma>0,\; \forall h>0,
\forall p_h\in Q_h,\quad ||B_h^t p_h||_{V'} \ge \gamma ||p_h||_{Q_h/\Ker(B^t)}.
\ee
Fortin~\cite{fortin2} gives a general useful method to check if the discrete inf-sup condition is satisfied.

Unfortunately, in many situations the stability condition~\eref{eqdis} is not satisfied.
\Eg\ $P_1$-continuous finite elements for both the velocity and pressure.

\medskip
The associated matrix problem relative to~\eref{eqgenh0} reads: Find $(U_h,P_h)\in\RR^{n_V} \times \RR^{n_Q}$ \st\ :
\be
\label{eqgenh0m}
\pmatrix{[A_h] & [B_h]^T \cr [B_h] & 0}\pmatrix{[U_h] \cr [P_h]} = \pmatrix{ [F_h] \cr [G_h]},
\ee
$[B_h]^T$ being $[B_h]$ transposed.
And if~\eref{eqdis} is not satisfied then the matrix $\pmatrix{[A_h] & [B_h]^T \cr [B_h] & [0]}$
is non invertible for some~$h$.
%(its inverse is not bounded as $h\rar0$).

%%%%%%%%%%%%%%%%%%%%%%%%%%%%%%%%%%%%%%%%%%%%%%%%%%%%%%%%%%%%%%%%%%%%%%%%%%%%%%%%%%%

\subsection{Where is the problem?}

\Eg, with~\eref{eqgenh0eg} and continuous $P_1$-continuous finite elements for $\vv_h$ and $p_h$ we have
\be
\label{eqgenh0eg2}
b(\vv_h,p_h) = (\vgrad p_h,\vv_h)_\Ld = (\PVh \vgrad p_h,\vv_h)_\Ld,
%(\dvg\vu_h,q_h)_\Ld = -(\vu_h,\vgrad q_h)_L_d
\ee
with $\PVh :\Ldo^n \rar \Vh$ the $\dd_\Ld$-orthogonal projection on~$\Vh$;
Here $\vgrad p_h$ is constant by element, and $\PVh$ is the projection on continuous $P_1$ functions.

This projection $\PVh$, as any projection, looses information:
Here we would like to consider $\vgrad p_h$ (to control~$p_h$),
but~\eref{eqgenh0eg2} tell us that only $\PVh \vgrad p_h$ is taken into account (is computed):
Since $\vgrad p_h = \PVh \vgrad p_h + (\vgrad p_h-\PVh \vgrad p_h)$,
we have lost $\vgrad p_h-\PVh \vgrad p_h$.
And, \eg\ with $P_1$-continuous finite elements for $\vv_h$ and $p_h$,
if nothing is done then the computation fails to give a good result (and it get worse as $h\rar0$).

To recover a well-posed problem, the missing term $\vgrad p_h-\PVh \vgrad p_h$ can be reintroduced, see~\eref{eqpbstokes2h}, and we then get an optimal result (\eg, order~$h$ for convergence for $P_1$-continuous finite elements).

%%%%%%%%%%%%%%%%%%%%%%%%%%%%%%%%%%%%%%%%%%%%%%%%%%%%%%%%%%%%%%%%%%%%%%%%%%%%%%%%%%%

%%%%%%%%%%%%%%%%%%%%%%%%%%%%%%%%%%%%%%%%%%%%%%%%%%%%%%%%%%%%%%%%%%%%%%%%%%%%%%%%%%%

%\newpage 

\part{Examples and preventions}

\section{Stokes model}
\label{secs}

%%%%%%%%%%%%%%%%%%%%%%%%%%%%%%%%%%%%%%%%%%%%%%%%%%%%%%%%%%%%%%%%%%%%%%%%%%%%%%%%%%%

\subsection{A first model}

Let $\Omega$ be a bounded open set.
Let $\dvg : \Huzon \rar\Ldzo$ with $\dd_\Huz = (\grad(.),\grad(.))_\Ld$,
and let
\be
V = \{ \vv \in \Huzon : \dvg\vv=0\}.
\ee

%\medskip
Problem (homogeneous Dirichlet type):
for $\vf\in \Hmuon$, find $\vu\in \Huzo$ \st
\be
\label{eqrsf00}
 {-}\Delta \vu = \vf .
\ee
Associated weak problem : find $\vu\in V$ \st
\be
\label{eqrsf}
(\grad\vu , \grad\vv)_\Ld = (\vf,\vv)_\Ld, \quad \forall \vv\in \Huzo .
\ee
The Lax--Milgram gives the well-posedness in~$(\Huzo,\dd_\Huz)$.

The optimized associated problem is: Find $\vu \in \Huzo$ realizing the minimum of % $J(\vu) = \min_{\vv\in \Huzo} J(\vv)$ of
\be
J(\vv) := \demi ||\vgrad\vv||_\Ld^2 - (\vf,\vv)_\Ld.
\ee

%%%%%%%%%%%%%%%%%%%%%%%%%%%%%%%%%%%%%%%%%%%%%%%%%%%%%%%%%%%%%%%%%%%%%%%%%%%%%%%%%%%

\subsection{Constrained associated problem}

The constraint $\dvg\vu=0$ is imposed with a Lagrangian multiplier~$p$:
The problem~\eref{eqrsf00} is transformed into: Find $(\vu,p)\in\Huzo \times \Ldzo$ \st
\be
\label{eqsto2}
\left\{\eqalign{
& (\grad \vu , \grad\vv)_\Ld - (p,\dvg \vv)_\Ld = (\vf,\vv)_\Ld, \quad \forall \vv \in \Huzon , \cr
& -(\dvg\vu,q)_\Ld = 0, \quad \forall q \in \Ldzo.
}\right.
\ee
We have obtained~\eref{eqgenw0} with
$V=\Huzon$, $Q=\Ldzo$, $g=0$, $B=\dvg : \Huzon \rar \Ldzo$ and
\be
\label{eqdefabsto}
\left\{\eqalign{
& a(\vu,\vv) = (\grad\vu , \grad\vv)_\Ld \quad\hbox{sur } \Huzon\times\Huzon, \cr
& b(\vv,q) = - (\dvg \vv,q)_\Ld \quad\hbox{sur } \Huzon \times \Ldzo.
}\right.
\ee

Since $B = \dvg : \Huzo \rar \Ldzo$, $\Ker B=V$,
$a\dd$ is coercive on~$\Ker(B)$ (it is on~$\Huzon$),
and $B^t = \vgrad : \Ldo \rar \Hmuon$ is surjective, \cf~theorem~\ref{thmrgo},
the theorem~\ref{thmis} applies, and the problem~\eref{eqsto2} is well-posed.

The associated weak problem reads: Find $(\vu,p)\in\Huzo \times \Ldzo$ \st
\be
\label{eqsto2F}
\left\{\eqalign{
& {-}\Delta \vu  + \vgrad p = \vf \in \Hmuo , \cr
& \dvg\vu = 0.
}\right.
\ee

The associated Lagrangian reads, find the saddle point in $\Huzon\times\Ldzo$ of
\be
\label{eqlagr}
\calL(\vv,q)
= \demi ||\grad\vv||_\Ld^2  - (q,\dvg \vv)_\Ld - (\vf,\vv)_\Ld.
\ee

%%%%%%%%%%%%%%%%%%%%%%%%%%%%%%%%%%%%%%%%%%%%%%%%%%%%%%%%%%%%%%%

\section{Numerical approximation of the Stokes model}

%%%%%%%%%%%%%%%%%%%%%%%%%%%%%%%%%%%%%%%%%%%%%%%%%%%%%%%%%%%%%%%

\subsection{Approximation}

Let $\Vh\subset \Huzo$ and $Q_h\subset\Ldzo$ (conform approximation to simplify) be finite dimension subspaces.
The discretization of~\eref{eqsto2} is: Find $\vu_h\in (\Vh)^n$ and $p_h\in Q_h$ \st
\be
\label{eqpbstokes1h}
\left\{\eqalign{
& (\grad\vu_h,\grad\vv_h)_\Ld - (p_h,\dvg\vv_h)_\Ld = (\vf,\vv_h)_\Ld,
           \quad\forall\vv_h\in (\Vh)^n,\cr
& (\dvg\vu_h,q_h)_\Ld = 0,\quad\forall q_h\in Q_h.\cr
}\right.
\ee

%%%%%%%%%%%%%%%%%%%%%%%%%%%%%%%%%%%%%%%%%%%%%%%%%%%%%%%%%%%%%%%

\subsection{Projections (finite element method)}

If $X_h$ is a subspace in~$\Ldo$, let
$\Pi_\Xh :
\left\{\eqalign{
\Ldo & \rar \Xh \cr
f & \rar \Pi_\Xh f
}\right\}$
be the $\dd_\Ld$-orthogonal projection on~$X_h$, that is,
\be
\forall f\in\Ldo,\quad (\Pi_\Xh f , x_h)_\Ld = (f,x_h)_\Ld,\quad \forall x_h\in X_h.
\ee
\Eg, if $X_h=P_1$ then $\Pi_{P_1} f \in P_1$ is the best approximation $P_1$ of~$f$
for the $\dd_\Ld$ inner product.
Similar notation for $X_h$ a subspace in~$\Ldon$.

Let
\be
\vgrad_h \eqdef \Pi_\Vh \circ \vgrad : 
\left\{\eqalign{
\Ldo & \rar \Vh \cr
p & \rar  \vgrad_h p = \Pi_\Vh(\vgrad p).
}\right.
\ee
So, $\vgrad_h p$ is characterized by  $(\vgrad_h p , \vv_h)_\Ld = ( \vgrad p , \vv_h)_\Ld$ pour tout $\vv_h\in \Vh$.
And~\eref{eqpbstokes1h} reads
\be
\label{eqpbstokes1h3}
\left\{\eqalign{
 & (\grad\vu_h,\grad\vv_h)_\Ld + (\vgrad_h p_h,\vv_h)_\Ld = (\vf,\vv_h)_\Ld,
           \quad\forall\vv_h\in \Vh,\cr
 & (\vu_h,\vgrad q_h)_\Ld = 0,\quad\forall q_h\in Q_h.,\cr
}\right.
\ee

\comment{
Let
\be
\label{eqdvgh}
\dvg_h = \Pi_\Qh \circ \dvg :
\left\{\eqalign{
\Huzon & \rar \Qh \cr
\vv & \rar  \dvg_h\vv = \Pi_\Vh(\dvg \vv).
}\right.
\ee
So, $\dvg_h\vv$ is characterized by $(\dvg_h\vv_h , q_h)_\Ld = (\dvg\vv , q_h)_\Ld$ for all $q_h\in Q_h$.
\be
\label{eqpbstokes1h2}
\left\{\eqalign{
 & (\grad\vu_h,\grad\vv_h)_\Ld - (p_h,\dvg_h\vv_h)_\Ld = (\vf,\vv_h)_\Ld,
           \quad\forall\vv_h\in \Vh,\cr
 & (\vu_h,\vgrad q_h)_\Ld = 0,\quad\forall q_h\in Q_h,\cr
}\right.
\ee
Donc~\eref{eqpbstokes1h} se lit aussi :
}

%%%%%%%%%%%%%%%%%%%%%%%%%%%%%%%%%%%%%%%%%%%%%%%%%%%%%%%%%%%%%%%

\subsection{Matrix representation}

With given bases in~$V_h$ and~$Q_h$, \eref{eqpbstokes1h} become
\be
\label{eqpbmat}
\pmatrix{A & B^T\cr B&0}.\pmatrix{\vx\cr \vy}=\pmatrix{F\cr0}.
\ee

%%%%%%%%%%%%%%%%%%%%%%%%%%%%%%%%%%%%%%%%%%%%%%%%%%%%%%%%%%%%%%%

\subsection{The problematic pressure}
\label{secpbp}

In many cases there is no problem with the computation of~$u_h$
(the $A$ matrix i~\eref{eqpbmat} is well conditioned since $a\dd$ is continuous and coercive).

But the results obtained for~$p_h$ can be absurd. To see why,
suppose that $u_h$ is known,
let $(g,\vv_h)_\Ld := (\vgrad\vu_h,\vgrad\vv_h)_\Ld-(\vf,\vv_h)_\Ld$,
and try to find $p_h\in Q_h$ \st
\be
(\vgrad p_h,\vv_h)_{\Hmu,\Huz} = -(g,\vv_h)_\Ld, \quad\forall\vv_h\in \Vh,
\ee
that is, \eg\ with continuous finite elements where
$\la \vgrad p_h,\vv_h\ra_{\Hmu,\Huz} = (\vgrad p_h,\vv_h)_\Ld$,
\be
\label{eqpgsjar}
(\vgrad_h p_h,\vv_h)_\Ld = -(g,\vv_h)_\Ld, \quad\forall\vv_h\in \Vh.
\ee

1- Nice case: $\vgrad_h = \Pi_\Vh \circ \vgrad : \Vh \rar \Qh$ is surjective (onto)
with a constant independent of~$h$, cf.~\eref{eqthmrgo}, that is,
\be
\label{eqisdish0}
\exists k>0,\; \forall h>0,\; \forall p_h\in Q_h,\quad
||\vgrad_h p_h||_{\Hmu} \ge k\, ||p_h||_\Ldz.
%    \inf_{\vv_h\in \Vh}\sup_{p_h\in Q_h} {(\dvg\vv_h,p_h)_\Ldo\over||\vv_h||_\Huz||p_h||_\Ld}\ge k.
\ee
And~\eref{eqisdish0} is called the ``discrete inf-sup condition''.

Then the problem~\eref{eqpbstokes1h3} is well-posed, \ie\ the matrix i~\eref{eqpbmat}
is well-conditioned, \cf\ theorem~\ref{thmis}.
See Fortin~\cite{fortin2} for $V_h$ and $Q_h$ finite element spaces that can satisfy~\eref{eqisdish0}.

(Remark: the problem~\eref{eqpgsjar} cannot be solved on its own in general,
since it is surjective but not bijective.
But~\eref{eqpbstokes1h3} can be solved, the matrix $\pmatrix{A & B^t\cr B&0}$ being invertible and well-conditioned if~\eref{eqisdish0} is satisfied.)

\medskip
2- Bad case: In~\eref{eqisdish0}, $k>0$ does not exists, \eg
\be
 \exists k_h>0,\quad
   \inf_{p_h\in Q_h}\sup_{\vv_h\in \Vh} {(\dvg\vv_h,p_h)_\Ldo\over||\vv_h||_\Huz||p_h||_\Ld}\ge k_h,
	\quad \hbox{but}\quad k_h\matrarrow_{h\rar0}0.
\ee
Then $\pmatrix{A & B^t\cr B&0}$ is not invertible (at least not numerically invertible as $h\rar0$:
bad conditioning).

\debexa
\label{eqadmloss}
A useful criteria to check the discrete inf-sup condition~\eref{eqisdish0}
is given by Fortin~\cite{fortin2}. \Eg, the discrete inf-sup condition is satisfied with the classical:

$P_2$,$P_1$ (velocity-pressure) Taylor--Hood finite elements
(see \eg\ Bercovier--Pironneau~\cite{bercovier-pironneau}).

$P_1$-bubble,$P_1$ (velocity-pressure) finite elements, named the mini-elements,
see Arnold--Brezzi--Fortin~\cite{arnold-brezzi-fortin}.

$P_2$,$P_0$ (velocity-pressure) finite elements, see Crouzeix--Raviart~\cite{crouzeix-raviart}.

(And for non conformity, the $P_1$-discontinuous velocity, $P_0$-pression,
see Crouzeix--Raviart~\cite{crouzeix-raviart}.)
\finexa

\debexa
\label{eqnadmloss}
No convergence \eg\ for the $P_1$,$P_1$ continuous finite elements,
or the $P_1$-continuous,$P_0$ elements (checkerboard instability).
\finexa
%\medskip
%Eléments finis $P_2$-continus  en vitesse, $P_0$ en pression, appelés éléments finis de Crouzeix et Raviart.

%%%%%%%%%%%%%%%%%%%%%%%%%%%%%%%%%%%%%%%%%%%%%%%%%%%%%%%%%%%%%%%

\subsection{What has been lost...}

\eref{eqpbstokes1h3} reads
\be
(\PVh\vgrad p_h,\vv_h)_\Ld = -(g,\vv_b),\quad\forall\vv_h\in \Vh.
\ee
So we want $\vgrad p_h$, but we can only compute $\vgrad_h p_h = \PVh\vgrad p_h$,
which in many cases is different from $\vgrad p_h$.
Since
\be
\vgrad p_h = \PVh\vgrad p_h + \bigl(\vgrad p_h- \PVh\vgrad p_h\bigr),
\ee
we have lost
\be
\label{eqperdu}
\hbox{loss} = (\vgrad p_h- \PVh\vgrad p_h) = (\vgrad p_h- \vgrad_h p_h).
\ee
This can be an admissible loss, see \eg\ example~\ref{eqadmloss},
or not, see \eg\ example~\ref{eqnadmloss}.

%%%%%%%%%%%%%%%%%%%%%%%%%%%%%%%%%%%%%%%%%%%%%%%%%%%%%%%%%%%%%%%

\subsection{... and a reintroduction of the loss}
\label{seccp}

\def\LdK{{L^2(K)}}

To recover the loss~\eref{eqperdu},
we modify~\eref{eqpbstokes1h} to get the new problem:
Find $\vu_h\in \Vh$ et $p_h\in Q_h$ \st
\be
\label{eqpbstokes2h}
\left\{\eqalign{
 & (\vgrad\vu_h,\vgrad\vv_h)_\Ld - (p_h,\dvg\vv_h)_\Ld = (\vf,\vv_h)_\Ld,
           \quad\forall\vv_h\in \Vh,\cr
 & -(\dvg\vu_h,q_h)_\Ld - \sum_{K=1}^{n_K} h_K^2(\vgrad p_h{-}\vgrad_h p_h,\vgrad q_h{-}\vgrad_h q_h)_\LdK = 0,
  \quad\forall q_h\in Q_h.\cr
}\right.
\ee
where $n_K$ is the number of elements constituting the mesh,
$h$ is ``the size of an element'', and
the $h_K^2$ coefficient to get optimal results, see Leborgne~\cite{leborgne}
(we are interested in~$p_h$ and, for quasi-uniform meshes, $p_h$ is of the same order than $h\,\vgrad p_h$).
\Eg\ for $P_1$,$P_1$ continuous finite elements for both the velocity and the pressure,
we get order~1 convergence results (classic for $P_1$ finite elements).

\eref{eqpbstokes2h} can also be written
\be
\label{eqpbstokes21h}
\left\{\eqalign{
 & (\vgrad\vu_h,\vgrad\vv_h)_\Ld - (p_h,\dvg\vv_h)_\Ld = (\vf,\vv_h)_\Ld,
           \quad\forall\vv_h\in \Vh,\cr
 & -(\dvg\vu_h,q_h)_\Ld - \sum_{K=1}^{n_K} h_K^2(\vgrad p_h{-}\PVh\vgrad p_h,\vgrad q_h)_\LdK = 0,
  \quad\forall q_h\in Q_h,\cr
}\right.
\ee
since $(\vgrad p_h{-}\PVh\vgrad p_h,\vw_h)=0$ for $\vw_h\in \Vh$ (definition of~$\Pi_\Vh$).

Computation:
we have to compute a new unknown $\vz_h=\PVh\vgrad p_h\in \Vh$ (luckily very cheap for $P_1$ finite elements):
Find $\vu_h,\vz_h\in \Vh$ et $p_h\in Q_h$ \st
\be
\label{eqpbstokes3h}
\left\{\eqalign{
 & (\vgrad\vu_h,\vgrad\vv_h)_\Ld - (p_h,\dvg\vv_h)_\Ld = (\vf,\vv_h)_\Ld,
           \quad\forall\vv_h\in \Vh,\cr
 & -(\dvg\vu_h,q_h)_\Ld - \sum_{K=1}^{n_K} h_K^2 (\vgrad p_h,\vgrad q_h)_\Ld + h^2 ( \vz_h,\vgrad q_h)_\Ld = 0,
  \quad\forall q_h\in Q_h,\cr
 & (\vgrad p_h,\vz_h')_\LdK -(\vz_h,\vz_h')_\LdK =0,\quad\forall\vz_h'\in \Vh,\; \forall K.\cr
}\right.
\ee
\Eg with $P_1$ finite elements, the $(\vz_h,\vz_h')_\Ld$ associated matrix
can be made diagonal thanks to the ``mass lumping'' technique:
Thus the last equation (in $\vz_h'$) gives $\vz_h$ explicitly as a function of~$\vgrad_h p_h$
(order~1 precision).

\debrem
The associated Lagrangian, \cf~\eref{eqlagr}, is now:
\be
\label{eqlso}
   L_h(\vv_h,p_h)=\demi||\grad\vv_h||^2_\Ld - (p_h,\dvg\vv_h)_\Ld - (f,v_h)_\Ld
  -\demi \sum_{K=1}^{n_K} h_K^2 ||\vgrad p_h{-}\PVh\vgrad p_h||^2_\LdK.
\ee
\finrem

%%%%%%%%%%%%%%%%%%%%%%%%%%%%%%%%%%%%%%%%%%%%%%%%%%%%%%%%%%%%%%%

\subsection{Brezzi and Pitkäranta's method}

A previous method proposed by Brezzi and Pitkäranta~\cite{brezzi-pitkaranta}
consists in penalizing the initial problem with the Laplacian of the pressure
(to ``control the oscillations'' of~$p_h$):
Find $\vu_h\in V_h$ and $p_h\in Q_h$ \st
\be
\label{eqpbstokesbph}
\left\{\eqalign{
 & (\vgrad\vu_h,\vgrad\vv_h)_\Ld - (p_h,\dvg\vv_h)_\Ld = (\vf,\vv_h)_\Ld,
           \quad\forall\vv_h\in V_h,\cr
 & -(\dvg\vu_h,q_h)_\Ld -\eps \sum_{K=1}^{n_K} h_K^2(\vgrad p_h,\vgrad q_h)_\LdK = 0,
  \quad\forall q_h\in Q_h,\cr
}\right.
\ee
with some $\eps>0$.
We however get a spurious limit condition ${\pa p\over\pa n}=0$ independent of~$\eps$ (by integration by parts).
(This spurious limit condition is lessen with~\eref{eqpbstokes21h}.)

\debrem
The associated Lagrangian is now:
\be
   L_h(\vv_h,p_h)=\demi||\grad\vv_h||^2_\Ld - (p_h,\dvg\vv_h)_\Ld - (f,v_h)_\Ld
  -\demi \eps \sum_{K=1}^{n_K} h_K^2 ||\vgrad p_h||^2_\LdK,
\ee
to compare with~\eref{eqlso}.
\finrem

%%%%%%%%%%%%%%%%%%%%%%%%%%%%%%%%%%%%%%%%%%%%%%%%%%%%%%%%%%%%%%%

\subsection{Hughes, Franca and Balestra's method}

Hughes, Franca and Balestra~\cite{hughes-franca-balestra} proposed a
``Galerkin Least-squares'' method: The pressure is stabilized ``with the solution''.
The problem reads, with the associated Lagrangian,
\be
\label{eqlhfb}
   L(\vv,p)=\demi||\grad\vv||^2_\Ldo - (p,\dvg\vv)_\Ldo - (f,v)_\Ldo
  -{\eps\over2} \sum_{K=1}^{n_K} h_K^2||-\Delta u{+}\vgrad p{-}f||^2_{L^2(K)}.
\ee
(For $P_1$ finite elements, this method is similar to Brezzi and Pitkäranta's method.)

Here $\eps$ has to be small enough not to destroy the coercivity in~$u$,
see the term $(\vgrad u,\vgrad v)_\Ldo-\eps \sum_kh^2(\Delta u,\Delta v)_{L^2(K)}$,
the control being done thanks to the inverse inequality (quasi-uniform mesh)
$$
||\Delta u_h||_\LdK \le C h ||\vgrad u_h||_\LdK,\quad\forall u_h\in V_h.
$$
(So $0< \eps<{1\over\sqrt C}$.)

%%%%%%%%%%%%%%%%%%%%%%%%%%%%%%%%%%%%%%%%%%%%%%%%%%%%%%%%%%%%%%%

\subsection{Douglas and Wang's method}

To avoid the eventual destruction of the coercivity for~$\vu$,
Douglas et Wang consider
\be
\label{eqldw}
\underbrace{
(\grad \vu,\grad \vv)_\Ld - (p,\dvg \vv) + (q,\dvg \vu) 
+\eps \sum_{K=1}^{n_K} h_K(-\Delta \vu+\vgrad p-\vf,-\Delta \vv+\vgrad q)_\LdK
}_{c((\vu,p),(\vv,q))}
= 
\underbrace{
(f,v)_\Ld
}_{\ell(\vv,q)}.
\ee
This preserves the stability since $c\dd$ is coercive, but the symmetry is lost.
So this method is adapted to the generalization of the Stokes equations to the Navier--Stokes equations.

%%%%%%%%%%%%%%%%%%%%%%%%%%%%%%%%%%%%%%%%%%%%%%%%%%%%%%%%%%%%%%%%%%%%%%%%%%%%%%%%%%%

\newpage

\section{Laplacian (harmonic problem)}

The linear spaces needed are described in~\S~\ref{secspaces}.

%%%%%%%%%%%%%%%%%%%%%%%%%%%%%%%%%%%%%%%%%%%%%%%%%%%%%%%%%%%%%%%%%%%%%%%%%%%%%%%%%%%

\subsection{A dirichlet problem}

Let $f\in\Hmuo$.
Problem: Find $p\in\Huzo$ \st
\be
\label{eqlap}
-\Delta p = f.
\ee
That is,
\be
\label{eqlap2}
(\vgrad p , \vgrad q)_\Ld = \la f, q \ra, \quad \forall  q\in\Huzo .
\ee
The associated minimum problem is:
Find the minimum of $J(p) = \min_{ q\in \Huzo} J(q)$,
where 
\be
J(q) := \demi ||\vgrad q||_\Ld^2 - \la f, q \ra.
\ee

To get $\vgrad p$ during the computation, introduce
\be
\label{eqvyu}
\vu = \vgrad p\; \in \Ldo, \quad \hbox{and then}\quad -\dvg\vu = f.
\ee
%(Donc $\vy$ dérive du potentiel~$\phi$.)
And~\eref{eqlap} becomes: Find $(\vu,p)\in\Ldo \times \Huzo$ \st
\be
\label{eqvyu2}
\left\{\eqalign{
& (\vu,\vv)_\Ld - (\vgrad p , \vv)_\Ld = 0, \quad \forall \vv \in \Ldo, \cr
& -(\vu, \vgrad q)_\Ld = - \la f, q \ra_{\Hmu,\Huz}, \quad \forall  q\in \Huzo .
}\right.
\ee
And $p$ is now the Lagrangian multiplier for the constraint $\dvg\vu = -f$, \cf\ the integration by parts.

And if $(\vu,p)\in\Ldo \times \Huzo$ is a solution, then
$\vu = \vgrad p \in \Ldon$ in~$\Omega$, and
$\dvg\vu = -f$ in~$\Hmuo$.
So $\Delta p= f$ in~$\Hmuo$ with $p \in \Huzo$: This is~\eref{eqlap}.

With
\be
\left\{\eqalign{
& a(\vu,\vv) = (\vu,\vv)_\Ld \quad\hbox{sur } \Ldon\times\Ldon, \cr
& b(\vv,q) = -( \vv , \vgrad q)_\Ld\quad\hbox{sur } \Ldon\times\Huzo,
}\right.
\ee
\eref{eqvyu2} has the appearance of~\eref{eqgenw0} with $V=\Ldon$ and $Q=\Huzo$.

Here $b(\vv,p) = \la B\vv,p \ra_{\Hmu,\Huo}  = (B^tp,\vv)_\Ldo$,
so $B = \dvg : 
\left\{\eqalign{
\Ldon & \rar \Hmuo \cr
\vv & \rar B\vv = \dvg(\vv)
}\right\}
$
and $B^t  = -\vgrad:
\left\{\eqalign{
\Huzo & \rar \Ldo \cr
p & \rar B^tp = -\vgrad p
}\right\}
$.

Thus $\Ker(B) = \Ker(\dvg) = \{ \vv\in\Ldon : \dvg\vv=0\}$,
and thanks to the Helmholtz decomposition~\eref{eqHe0}
$\Ldon = \vgrad(\Huzo) \oplus^{\perp_\Ld} \Ker(\dvg)$,
the bilinear form $a\dd$ is $\dd_\Ld$ coercive on~$\Ker(B)$.

And $B^t$ is surjective since $B$ is, \cf~\eref{eqcisdiv3b}
and the closed range theorem~\ref{thmcr}.

Thus~\eref{eqvyu2} is well-posed, see theorem~\ref{thmis}.

%%%%%%%%%%%%%%%%%%%%%%%%%%%%%%%%%%%%%%%%%%%%%%%%%%%%%%%%%%%%%%%%%%%%%%%%%%%%%%%%%%%

\subsection{A Neumann problem}

Let $f\in\Ldo$.
Problem: Find $p\in\Huo$ \st
\be
-\Delta p = f, \qand {\pa p \over \pa n}_{|\Gamma} = 0.
\ee
That is,
\be
(\vgrad p , \vgrad q)_\Ld = (f, q)_\Ld, \quad \forall  q\in\Huzo .
\ee
The associated minimum problem is:
Find the minimum of $J(p) = \min_{ q\in \Huo} J(q)$,
where 
\be
J(q) := \demi ||\vgrad q||_\Ld^2 -(f, q)_\Ld.
\ee
The mixed associated problem is: Find $(\vu,p)\in\Hdvgo \times \Ldo$ \st
\be
\label{eqvyu3}
\left\{\eqalign{
& (\vu,\vv)_\Ld + (p ,\dvg \vv)_\Ld = 0, \quad \forall \vv \in \Hdvgo, \cr
& (\dvg\vu , q)_\Ld = ( f, q)_\Ld, \quad \forall  q\in \Ldo .
}\right.
\ee
Indeed, if $(\vu,p)\in\Hdvgo \times \Ldo$ is a solution, then
$\dvg\vu=f \in \Ldo$, $\vu = -\vgrad p\in\Hmuo$, thus $-\Delta p=f\in\Ldo$,
with ${\pa p \over \pa n}_{|\Gamma} = 0$ (since $\Im(\gamma_n)=\Hmdemig$).

With
\be
\left\{\eqalign{
& a(\vu,\vv) = (\vu,\vv)_\Ld  \quad\hbox{sur } \Hdvgo\times\Hdvgo, \cr
& b(\vv,q) = (\dvg \vv , q)_\Ld  \quad\hbox{sur } \Hdvgo\times\Ldo.
}\right.
\ee
\eref{eqvyu2} has the appearance of~\eref{eqgenw0} with $V=\Hdvgo$ and $Q=\Ldo$.

Here $B : \vv\in\Hdvgo \rar B\vv=\dvg\vv\in\Ldo$ is surjective, \cf~\eref{eqdiv11},
and $a\dd$ est $\Hdvg$-coercive on $\Ker(B) = \{\vv\in \Hdvgo : \dvg\vv = 0\}$.
And $\Im(B)$ being closed (since it is surjective), so is $\Im(B^t)$ (closed range theorem\ref{thmcr}).
Thus~\eref{eqvyu3} is well-posed, see theorem~\ref{thmis}.

% $\vu$ est la variable ``principale'' et $p$ est le multiplicateur de Lagrange de la contrainte $\dvg\vu = f$.

%Et par intégration par parties, $(\phi ,\dvg \vz)_\Ld = -(\vgrad \phi,\vz)_\Ld + \int_\Gamma \phi\,\vz.\vn\,d\Gamma$.

%%%%%%%%%%%%%%%%%%%%%%%%%%%%%%%%%%%%%%%%%%%%%%%%%%%%%%%%%%%%%%%%%%%%%%%%%%%%%%%%%%%

\newpage

\section{Bilaplacian (biharmonic problem)}

The linear spaces needed are described in~\S~\ref{secspaces}.

%%%%%%%%%%%%%%%%%%%%%%%%%%%%%%%%%%%%%%%%%%%%%%%%%%%%%%%%%%%%%%%%%%%%%%%%%%%%%%%%%%%

\subsection{Problem}

We look here at the Dirichlet problem: If $f\in\Hmdo = (\Hdzo)'$, then find $p\in\Hdzo$ \st
\be
\label{eqbilap}
\Delta (\Delta p) = f,\quad\hbox{so with}\quad p_{|\Gamma}=0 \qand {\pa p\over\pa n}_{|\Gamma}=0.
\ee
Weak form: Find $p\in\Hdzo$ \st
\be
\label{eqbilapv}
(\Delta p , \Delta q)_\Ld = \la  f, q \ra_{\Hmd,\Hdz}, \quad \forall  q\in\Hdzo .
\ee
The Lax--Milgram theorem indicates that~\eref{eqbilapv} is well-posed.

The associated minimum problem is:
Find the minimum of $J(p) = \min_{ q\in \Huzo} J(q)$,
where 
\be
\label{eqJbil}
J(q) :=  \demi ||\Delta q||_\Ld^2 - \la  f, q \ra_{\Hmd,\Hdz}.
\ee

%%%%%%%%%%%%%%%%%%%%%%%%%%%%%%%%%%%%%%%%%%%%%%%%%%%%%%%%%%%%%%%%%%%%%%%%%%%%%%%%%%%

\subsection{Introduction of $\Delta p$}

(Not conclusive.)
A function in $\in\Hdo$, \st\ $\Delta p$ and~$\Delta q$ in~\eref{eqbilapv}, is cumbersome to approximate,
\cf\ the $C^1$ Argyris finite elements. Let
\be
\phi = \Delta p
\ee
Then problem~\eref{eqbilap} is rewritten as:
Find $(\phi,p)\in\Ldo\times \Hdzo$ \st
\be
\left\{\eqalign{
& \phi = \Delta p,  \cr
& \Delta \phi = f.
}\right.
\ee
And the weak form is, if $f\in\Hmuo$:
Find $(\phi,p)\in\Huo \times \Huzo$ \st
\be
\label{eqffbil1}
\left\{\eqalign{
& (\phi,\psi)_\Ld + (\vgrad p, \vgrad \psi)_\Ld = 0, \quad \forall \psi\in \Huo, \cr
& (\vgrad \phi , \vgrad q)_\Ld = -\la f,q\ra_{\Hmu,\Huz}, \quad \forall q\in\Huzo.
}\right.
\ee
Indeed, if $(\phi,p)\in\Huo \times \Huzo$ is as solution of~\eref{eqffbil1},
then $\phi - \Delta p=0$ (thus $\Delta p\in \Ldo$), and ${\pa p\over \pa n}_{|\Gamma}=0$.
And $p\in\Huzo$ with $\Delta\phi = f\in \Hmuo$, thus $\Delta^2p=f$.

With
\be
\left\{\eqalign{
& a(\phi,\psi) = (\phi,\psi)_\Ld \quad \hbox{sur } \Huo \times \Huo , \cr
& b(\phi,v) =  (\vgrad \phi, \vgrad v)_\Ld \quad \hbox{sur } \Huo \times \Huzo,
}\right.
\ee
\eref{eqvyu2} has the appearance of~\eref{eqffbil1} with $V=\Huo$ and $Q=\Huzo$.

Here $b(\phi,v)=- \la \Delta \phi , v \ra_{\Hmu,\Huz}$, thus $B = -\Delta : \Huo \rar \Hmuo$.

Thus $B : \phi \in \Huo \rar B\phi = -\Delta \phi \in \Hmuo$ is surjective:
Apply Lax--Milgram theorem for $g\in\Hmuo$ and $(\vgrad\phi,\vgrad\psi)_\Ld = \la g,\psi\ra_{\Hmu,\Huz}$.

And $\Ker B=\{\phi \in \Huo : \Delta\phi = 0\}$ (harmonic functions).
So $\phi\in\Ker B$ iff $(\vgrad\phi ,\vgrad v)_\Ld = 0$ for all $v\in\Huzo$
(this is not $(\vgrad\phi ,\vgrad v)_\Ld = 0$ for all $v\in\Huo$).
Thus $a\dd$ is not $\dd_\Hu$-coercive on~$\Ker B$, but only $\dd_\Ld$ coercive, and the usual theorem is not applicable:
A loss of precision (precision$||.||_\Ld$ instead of precision~$||.||_\Hu$ for~$\phi$)
is to be expected.
%D'où le~\S~suivant~\ref{secbilb}.

%%%%%%%%%%%%%%%%%%%%%%%%%%%%%%%%%%%%%%%%%%%%%%%%%%%%%%%%%%%%%%%%%%%%%%%%%%%%%%%%%%%

\subsection{Introduction of $\vgrad p$}
\label{secbilb}

\def\tdd{\mathop{\,\tilde:}\,}
\def\Guzdvgo{{G^{1,\dvg}_{0,0}(\Omega)}}

%%%%%%%%%%%%%%%%%%%%%%%%%%%%%%%%%%%%%%%%%%%%%%%%%%%%%%%%%%%%%%%%%%%%%%%%%%%%%%%%%%%

\subsubsection{Weak form}

%On~a $\Delta p = \dvg(\vgrad p)$, et 
In~\eref{eqbilapv} let us introduce
\be
\label{eqblapv301}
\vu = \vgrad p , \quad\hbox{thus}\quad \Delta p = \dvg\vu.
\ee
Notation:
\be
\label{eqhgradnot}
\hbox{if } \vv = \vgrad q\in  \vgrad(\Huzo), \quad\hbox{then}\quad
\vv \eqnote  \vv_q.
\ee
(That is, $\vv_q$ derives from a potential $q\in\Huzo$.)

Then~\eref{eqbilapv} becomes: 
Find $\vu_p = \vgrad p \in \vgrad(\Huzo)$ \st
\be
\label{eqblapv30}
(\dvg\vu_p,\dvg\vv_q)_\Ld = \la f, q \ra_{\Hmd,\Hdz}, \quad \forall q \in\Hdzo.
\ee

%%%%%%%%%%%%%%%%%%%%%%%%%%%%%%%%%%%%%%%%%%%%%%%%%%%%%%%%%%%%%%%%%%%%%%%%%%%%%%%%%%%

\subsubsection{A first constrained formulation}

(Not conclusive.)
To avoid working with the ``small space'' $\vgrad(\Huzo) \ni \vu$, 
we consider the whole space $\Huzo\ni\vu$ and
we add the constraint $\vu-\vgrad p=0$, \cf~\eref{eqblapv301}
(and the associated Lagrangian multiplier~$\vlambda$).

And $\vu \in\Huo$ and $\vu = \vgrad p$ for some $p\in\Hdzo$ solution of~\eref{eqbilapv},
give $\dvg\vu = \Delta p \in \Ldo$ with $0 = {\pa p \over \pa n}_{|\Gamma} = \vu.\vn_{|\Gamma}$,
so $\vu \in \Hdvgzo$.
Then let
\be
X = \Hdvgzo \times \Huzo, 
\ee
provided with the inner product
\be
((\vu,p),(\vv,q))_X = (\dvg\vu,\dvg\vv)_\Ld + (\vgrad p,\vgrad q)_\Ldo,
% = (\dvg\vu,\dvg\vv)_\Ld + (p,q)_\Huz,
\ee
so that $X$ is a Hilbert space.

Then, if $f\in\Hmuo$, \eref{eqblapv30} is turned into:
Find $((\vu,p),\vlambda) \in X \times \Ldon$ \st
\be
\label{eqblapv2l0}
\left\{\eqalign{
& (\dvg\vu,\dvg\vv)_\Ld + (\vlambda,\vv - \vgrad q)_\Ld = \la f,q\ra_{\Hmu,\Huz}, \quad \forall (\vv,q) \in X, \cr
& (\vu - \vgrad p,\vmu)_\Ld = 0 , \quad \forall \vmu \in\Ldon, \cr
}\right.
\ee
that is,
\be
\label{eqblapv2ld}
\left\{\eqalign{
& (\dvg\vu,\dvg\vv)_\Ld + (\vlambda , \vv)_\Ld =0,\quad \forall \vv \in \Hdvgzo,\cr
& -(\vlambda , \vgrad q)_\Ld = \la f,q\ra_{\Hmu,\Huz},\quad \forall q \in \Huzo,\cr
&  (\vu , \vmu)_\Ld - (\vgrad p , \vmu)_\Ld  = 0, \quad \forall \vmu \in \Ldon.
}\right.
\ee
Check:
If $((\vu,p),\vlambda) \in X \times \Ldon$ is a solution of~\eref{eqblapv2l0} or~\eref{eqblapv2ld},
then
$\vlambda = \vgrad(\dvg \vu)$,
$\dvg\vlambda = f$,
$\vu = \vgrad p$,
thus $\dvg\vu = \Delta p$
and
$\dvg(\vgrad(\Delta p))=f$, \ie\ $\Delta^2 p =f$.
And $\vu \in \Hdvgzo$ gives $\vu.\vn=0$, so $\vgrad p.\vn=0$,
and with $p\in\Huzo$ we get $p\in \Hdzo$.

With
\be
\label{eqblapv2i}
\left\{\eqalign{
& a((\vu,p),(\vv,q)) = (\dvg\vu,\dvg\vv)_\Ld \quad \hbox{sur } X \times X, \cr
& b((\vv,q),\vmu) = (\vv - \vgrad q,\vmu)_\Ld \quad \hbox{sur } X \times \Ldon,
}\right.
\ee
\eref{eqblapv2l0} has the appearance of~\eref{eqgenw0} with $V=X$ and $Q=\Ldon$.

Here $ B : 
\left\{\eqalign{
\Hdvgzo \times \Huzo & \rar  \Ldon \cr
(\vv,q) & \rar B(\vv,q) = \vv - \vgrad q
}\right\}
$.

And $\Ker(B) =\{(\vv,q) \in \Hdvgzo \times \Huzo : \vv = \vgrad q\}$.
Thus is $(\vv,q) \in \Ker(B)$ then $\Delta p\in\Ldo$ and
$a((\vu,p),(\vv,q)) = \demi (\dvg\vu,\dvg\vv)_\Ld + \demi(\Delta p, \Delta q)_\Ld$.
And when $p\in\Hdo\cap \Huzo$ we have $||\Delta p||_\Ld \ge C\,||p||_\Hd \ge C\,||p||_\Hu$,
\cf~\eref{eqpoinc2}.
%voir par exemple Raviart--Thomas~\cite{raviart-thomas}.
Thus $a\dd$ is coercive on $(\Ker(B),\dd_X)$.

But $B$ is not surjective:
If $\vell \in \Ldon$ we should find $(\vv,q) \in \Hdvgzo \times \Huzo$
\st\ $\vell = \vv - \vgrad q$, but we only have~\eref{eqHe0}.
So $\vlambda$ has a priori no $||.||_\Ldo$ control,
and for the discretization we expect a loss of precision.
Here $B^t :
\left\{\eqalign{
D\subset\Ldon & \rar \Hdvgzo' \times \Hmuo  \cr
\mu & \rar B^t\mu : B^t\mu(\vv,q) = \la \vv,\mu \ra_{\Hdvgz',\Hdvg} + \la q,\dvg\mu\ra_{\Hmu,\Huz}
}\right\}
$ where $D= \Hdvgzo$ (the domain of definition) is not closed in~$\Ldo$ (its closure is $\Ldo$).

%%%%%%%%%%%%%%%%%%%%%%%%%%%%%%%%%%%%%%%%%%%%%%%%%%%%%%%%%%%%%%%%%%%%%%%%%%%%%%%%%%%

\subsubsection{A second constrained formulation}
\label{secbilapb}

(Not conclusive.)
Let
\be
X_+ = \Huzon \times \Huzo,
\ee
provided with the inner product
\be
((\vu,p),(\vv,q))_{X_+} = (\grad\vu,\grad\vv)_\Ld + (\vgrad p,\vgrad q)_\Ldo,
\ee
so that $X_+$ is a Hilbert space.

%Pour un contrôle~$\Hu$ de~$\vphi$, $\dvg\vphi$ est insuffisant, voir~\eref{eqqpuX}.
With a Cartesian basis, we notice that $(\Delta p , \Delta q)_\Ld
= \sum_{ij} \int_\Omega{\pa^2 p \over \pa x_i^2}{\pa^2 q \over \pa x_j^2}\,d\Omega$,
and, for $p,q \in \Hdzo$,
\be
\label{eqblapv2}
\int_\Omega  {\pa^2 p \over \pa x_i^2}{\pa^2 q \over \pa x_j^2}\,d\Omega
= - \int_\Omega  {\pa^3 p \over \pa x_i^2\pa x_j}{\pa q \over \pa x_j}\,d\Omega
= \int_\Omega  {\pa^2 p \over \pa x_i \pa x_j}{\pa^2 q \over \pa x_i\pa x_j}\,d\Omega.
\ee
Thus, with~\eref{eqhgradnot} and $\vv_p,\vv_q \in \vgrad(\Huzo)$ \cf, we get
\be
(\dvg\vv_p, \dvg\vv_q)_\Ld = (\Delta p , \Delta q)_\Ld
= (\grad(\vgrad p),\grad(\vgrad q))_\Ld = (\grad\vv_p,\grad\vv_q)_\Ld.
\ee
(Nota Bene: Here $\vv_p$ and $\vv_q$ derives from a potential.)
Thus~\eref{eqblapv30} reads:
Find $\vu=\vv_p\in \vgrad(\Huzo)$ \st
\be
\label{eqblapv31}
(\grad\vu,\grad\vv_q)_\Ld = ( f, q)_\Ld, \quad \forall q \in \Huzo.
\ee
So, if $f\in\Hmuo$, then \eref{eqblapv30} is transformed into:
Find $((\vu,p),\vlambda) \in X_+ \times \Ldon$ \st
\be
\label{eqblapv2l0b}
\left\{\eqalign{
& (\grad\vu,\grad\vv)_\Ld + (\vlambda,\vv - \vgrad q)_\Ld = \la f,q\ra_{\Hmu,\Huz}, \quad \forall (\vv,q) \in X_+, \cr
& (\vu - \vgrad p,\vmu)_\Ld = 0 , \quad \forall \vmu \in\Ldon, \cr
}\right.
\ee
$\vlambda$ being the Lagrangian multiplier of the constraint~\eref{eqblapv301}.
That is,
\be
\left\{\eqalign{
& (\grad\vu,\grad\vv)_\Ld + (\vlambda , \vv)_\Ld =0,\quad \forall \vv \in \Huzon,\cr
& -(\vlambda , \vgrad q)_\Ld = \la f,q\ra_{\Hmu,\Huz},\quad \forall q \in \Huzo,\cr
&  (\vu , \vmu)_\Ld - (\vgrad p , \vmu)_\Ld  = 0, \quad \forall \vmu \in \Ldon.
}\right.
\ee
With
\be
\left\{\eqalign{
& a((\vu,p),(\vv,q)) = (\grad\vu,\grad\vv)_\Ld \quad \hbox{sur } X_+ \times X_+, \cr
& b((\vv,q),\vmu) = (\vv, \vmu)_\Ld - (\vgrad p , \vmu)_\Ld  \quad \hbox{sur } X_+ \times \Ldon,
}\right.
\ee
\eref{eqblapv2l0b} has the appearance of~\eref{eqgenw0} with $V=X_+$ and $Q=\Ldon$.

Et $ B : 
\left\{\eqalign{
\Huzon \times \Huzo & \rar  \Ldon \cr
(\vv,q) & \rar B(\vv,q) = \vv- \vgrad q
}\right\}
$.
And $(\vv,q) \in \Ker(B)$ iff $\vgrad q = \vv \in \Huzon$,
thus $a\dd$ is coercive on $(\Ker(B),\dd_{X_+})$.
(Compared to~\eref{eqblapv2i}, here we a $||.||_\Hu$ for~$\vu$,
not only a~$||.||_\Hdvg$ control.)

However $B$ is not surjective. And $\vlambda$ is not controlled the classical way.

%%%%%%%%%%%%%%%%%%%%%%%%%%%%%%%%%%%%%%%%%%%%%%%%%%%%%%%%%%%%%%%%%%%%%%%%%%%%%%%%%%%

\subsubsection{A third constrained formulation}
\label{secbilapa}

(``A good one''.)
Let
\be
Y = \Hdvgzo \times \Huo,
\ee
provided with the inner product
\be
((\vu,p),(\vv,q))_Y = (\dvg\vu,\dvg\vv)_\Ld + (p,q)_\Hu
\ee
so that $Y$ is a Hilbert space.

If $f\in\Ldo$ (or in $(\Huo)'$), \eref{eqblapv2l0} is transformed into:
Find $((\vu,p),\vlambda) \in Y \times \Hdvgo$ \st
\be
\label{eqblapv22}
\left\{\eqalign{
& (\dvg\vu,\dvg\vv)_\Ld + (\vlambda,\vv)_\Ld + (q, \dvg \vlambda)_\Ld = ( f, q)_\Ld, \quad \forall (\vv,q) \in Y, \cr
& (\vu,\vmu)_\Ld + (\vlambda,\dvg\vmu) = 0 , \quad \forall \vmu \in \Hdvgo, \cr
}\right.
\ee
that is,
\be
\label{eqblapv22c}
\left\{\eqalign{
& (\dvg\vu,\dvg\vv)_\Ld + (\vlambda , \vv)_\Ld =0,\quad \forall \vv \in \Hdvgzo,\cr
& (\dvg\vlambda , q)_\Ld = ( f, q)_\Ld,\quad \forall q \in \Huo,\cr
&  (\vu , \vmu)_\Ld + (p , \dvg\vmu)_\Ld  = 0, \quad \forall \vmu \in \Hdvgo.
}\right.
\ee
So
$\vlambda = \vgrad(\dvg \vu)$, $\dvg\vlambda = f$, $\vu = \vgrad p$,
thus $\dvg(\vgrad(\dvg \vgrad p))=f$, \ie\ $\Delta^2 p =f$.
And
$\int_\Gamma p\,\vmu.\vn\,d\Gamma=0 = \la p,\vmu.\vn\ra_{\Hdemig,\Hmdemig}$ for all $\vmu\in\Hdvgo$,
thus $p_{|\Gamma}=0$ in~$\Hdemig$
(the trace operator $\vv \in \Hdvgo \rar \vv.\vn \in \Hmdemig$ is surjective).
Thus $p\in\Huzo$. With $\vu\in\Hdvgzo$, thus $\vgrad p.\vn=\vu.\vn=0$, so $p\in\Hdzo$.

With
\be
\label{eqblapv22d}
\left\{\eqalign{
& a((\vu,p),(\vv,q)) = (\dvg\vu,\dvg\vv)_\Ld \quad \hbox{sur } Y \times Y, \cr
& b((\vv,q),\vmu) = (\vv, \vmu)_\Ld + (q, \dvg \vmu)_\Ld  \quad \hbox{sur } Y \times \Hdvgo,
}\right.
\ee
\eref{eqblapv22} has the appearance of~\eref{eqgenw0} with $V=Y$ and $Q=\Hdvgo$.

%Et $b((\vv,q),\vmu) = (\vv, \vmu)_\Ldo - (\vgrad q, \vmu)_\Ldo + (q,\vmu.\vn)_\Ldg $.
And $ B : 
\left\{\eqalign{
\Hdvgzo \times \Huo & \rar  \Hdvgo' \cr
(\vv,q) & \rar B(\vv,q),
}\right\}
$
with $\la B(\vv,q), \vmu\ra = (\vv, \vmu)_\Ld + (q, \dvg \vmu)_\Ld$.
So $B$ is surjective, \cf~\eref{eqhdvgd}.

And $b((\vv,q),\vmu) = (\vv, \vmu)_\Ldo - (\vgrad q, \vmu)_\Ldo + (q,\vmu.\vn)_\Ldg $
gives $(\vv,q) \in \Ker(B)$ iff $(\vv,q) \in \Hdvgzo \times \Huzo$
with $\vv = \vgrad q$. Thus $a\dd$ is coercive on $(\Ker(B),\dd_Y)$,
and the classical theorem~\ref{thmis} apply.

\debrem
(Neumann.)
With $Z = \Hdvgo\times \Huzo$, \eref{eqblapv22} is transformed into:
Find $((\vu,p),\vlambda) \in Z \times \Hdvgzo$ \st
\be
\left\{\eqalign{
& (\dvg\vu,\dvg\vv)_\Ld + (\vlambda , \vv)_\Ld =0,\quad \forall \vv \in \Hdvgo,\cr
& (\dvg\vlambda , q)_\Ld = ( f, q)_\Ld,\quad \forall q \in \Huzo,\cr
&  (\vu , \vmu)_\Ld + (p , \dvg\vmu)_\Ld  = 0, \quad \forall \vmu \in \Hdvgzo.
}\right.
\ee
So, in~$\Omega$,
$\vlambda = \vgrad(\dvg \vu)$, $\dvg\vlambda = f$, $\vu = \vgrad p$,
thus $\dvg(\vgrad(\dvg \vgrad p))=f$, \ie\ $\Delta^2 p =f$.
And on~$\Gamma$, $\vgrad(\dvg\vu).\vn=0$, thus $\vgrad(\Delta p).\vn = 0$ (Neumann limit condition), with
 $p\in\Huzo$.
\finrem

%%%%%%%%%%%%%%%%%%%%%%%%%%%%%%%%%%%%%%%%%%%%%%%%%%%%%%%%%%%%%%%%%%%%%%%%%%%%%%%%%%%

\subsubsection{A fourth constrained formulation}

Let
\be
Y_+ = \Huzo \times \Huo,
\ee
provided with the inner product
\be
((\vu,p),(\vv,q))_{Y_+} = (\vgrad\vu,\vgrad\vv)_\Ld + (p,q)_\Hu.
\ee
And~\eref{eqblapv22c} is replaced with
\be
\left\{\eqalign{
& (\grad\vu,\grad\vv)_\Ld + (\vlambda , \vv)_\Ld =0,\quad \forall \vv \in \Hdvgzo,\cr
& (\dvg\vlambda , q)_\Ld = ( f, q)_\Ld,\quad \forall q \in \Huo,\cr
&  (\vu , \vmu)_\Ld + (p , \dvg\vmu)_\Ld  = 0, \quad \forall \vmu \in \Hdvgo.
}\right.
\ee

And~\eref{eqblapv22d} is replaced with
\be
\label{eqblapv22e}
\left\{\eqalign{
& a((\vu,p),(\vv,q)) = (\grad\vu,\grad\vv)_\Ld \quad \hbox{sur } Y_+ \times Y_+, \cr
& b((\vv,q),\vmu) = (\vv, \vmu)_\Ld + (q, \dvg \vmu)_\Ld  \quad \hbox{sur } Y_+ \times \Hdvgo
}\right.
\ee

%%%%%%%%%%%%%%%%%%%%%%%%%%%%%%%%%%%%%%%%%%%%%%%%%%%%%%%%%%%%%%%
%%%%%%%%%%%%%%%%%%%%%%%%%%%%%%%%%%%%%%%%%%%%%%%%%%%%%%%%%%%%%%%

\newpage

\section{Locking}

The locking phenomenon appears when the coercivity of the approximated problem increases much faster than
the coercivity of the continuous problem.
Thus the numerical solution is close to~zero, which is absurd in general.

Let $\Omega$ be bounded open set in~$\RRn$.

%%%%%%%%%%%%%%%%%%%%%%%%%%%%%%%%%%%%%%%%%%%%%%%%%%%%%%%%%%%%%%%

\subsection{A typical situation}

Let $\lambda\in\RR$ so that $\lambda >>1$ (a ``large'' given real).
We look for $\vu \in \Huzon$ and $p \in \Huzo$
that minimize
\be
\label{eqloc0}
  M(\vv,q)= \demi ||\grad\vv||^2_\Ldo + {\lambda\over2}||\vv-\vgrad q||_\Ldo
   -(\vf,\vv)_\Ldo -(g,q)_\Ldo.
\ee

\debexa
For the Mindlin--Reissner problem,
$||\grad\vv||^2_\Ldo$ is replaced by $|a(\vv,\vv)|$
where $a\dd$ %:(\Huzo)^n\times(\Huzo)^n\rar\RR$ 
is a bilinear form that is continuous and coercive on~$\Huzo$.
\finexa

Let
\be
X=\Huzo^n\times\Huzo
\ee
provided with the inner product associated to the norm
\be
||(\vv,q)||_X=(||\grad\vv||_\Ld^2+||\vgrad q||_\Ld^2)^\demi
\ee
so that $X$ is a Hilbert space.

A solution $(\vu,p)\in X$ realizing the min of~$M$ satisfies
\be
\label{eqloc1}
\left\{\eqalignrllq{
 & (\grad\vu,\grad\vv)_\Ld + \lambda(\vu-\vgrad p,\vv)_\Ld = (\vf,\vv),&\forall\vv\in\Huzo^n,\cr
 & \lambda(\vu-\vgrad p,\vgrad q)_\Ld = (g,q),&\forall q\in\Huo.\cr
}\right.
\ee
Let
\be
\label{eqloc20}
  \Phi((\vu,p),(\vv,q))=(\grad\vu,\grad\vv)_\Ld + \lambda(\vu-\vgrad p,\vv-\vgrad q)_\Ld.
\ee
Thus~\eref{eqloc1} reads: Find $(\vu,p)\in X$ \st
\be
\label{eqloc2}
   \Phi((\vu,p),(\vv,q)) = (\vf,\vv)_\Ld + (g,q)_\Ld , \quad \forall (\vv,q)\in X.
\ee

\debprop
The bilinear form $\Phi:X\times X\rar\RR$ is coercive and continuous on~$X$:
with the Poincaré inequality~\eref{eqpoinc} we have
\be
\left\{\eqalign{
& \exists \alpha_\Phi > 0, \quad \forall (\vv,q)\in X, \quad
\Phi((\vv,q),(\vv,q)) \ge \alpha_\Phi ||(\vv,q)||_X^2, \qe
\alpha_\Phi \mathop{\sim}_{\lambda\rar\infty} {1\over c_\Omega} , \cr
& \exists C > 0, \quad \forall (\vu,p),(\vv,q)\in X, \quad
\Phi((\vu,p),(\vv,q)) \le C||(\vu,p)||_X||(\vv,q)||_X, \qe
  C \mathop{=}_{\lambda\rar\infty} O(\lambda). %\quad\hbox{au vois de }\lambda=\infty.
}\right.
\ee
And problem~\eref{eqloc2} is well posed.
\finprop

\debdem
Bi-linearity. Since $\Phi$ is symmetric (trivial), that is $\Phi((\vu,p),(\vv,q))=\Phi((\vv,q),(\vu,p))$,
we have to prove that 
$\Phi((\vu_1,p_1)+\alpha(\vu_2,p_2),(\vv,q))
=\Phi((\vu_1,p_1),(\vv,q))+\alpha\Phi((\vu_2,p_2),(\vv,q))$: trivial.

Coercivity.
If $\kappa>0$ then, with~\eref{eqpoinc} :
$$
\eqalign{
  \Phi((\vv,q),(\vv,q))
 =  & ||\grad\vv||_\Ld^2+\lambda||\vv-\vgrad q||_\Ld^2 \cr
 \ge&[(1{-}\kappa)||\grad\vv||_\Ld + c_\Omega\kappa||\vv||_\Ld^2]
  +\lambda[||\vv||_\Ld^2+||\vgrad q||_\Ld^2-2||\vv||_\Ld||q||_\Ld].
}
$$
Let $x=||\vv||_\Ld$ and $y=||\vgrad q||_\Ld$. We have
$$
   c_\Omega\kappa x^2 + \lambda(x-y)^2 \ge {\lambda c_\Omega\kappa\over\lambda+ c_\Omega\kappa} y^2,
$$
with $c_{\rm max}={\lambda c_\Omega\kappa\over\lambda+ c_\Omega\kappa}$ the largest constant possible
($c_{\rm max}$ is largest $c$ \st\ $c_\Omega\kappa x^2 + \lambda(x-y)^2 \ge cy^2$: easy check). 
Thus
$$
  \Phi((\vv,q),(\vv,q)) \ge 
(1{-}\kappa)||\grad\vv||_\Ld+{\lambda c_\Omega\kappa\over\lambda+ c_\Omega\kappa}||\vgrad q||_\Ld.
$$
And the ``best'' $\alpha_\Phi$ (the largest possible)
is obtained by choosing $\kappa$ \st\
$1{-}\kappa={\lambda c_\Omega\kappa\over\lambda+ c_\Omega\kappa}$, \ie\ $\kappa$
solution of $\kappa^2+b\kappa -{\lambda\over c_\Omega}=0$ where $b=\lambda{ c_\Omega+1\over c_\Omega}-1$.
The discriminant is $b^2+4{\lambda\over c_\Omega}$, and the positive root is
$\kappa={b\over 2}(-1+\sqrt{1+4{\lambda\over b^2 c_\Omega}})$.
And for $\lambda>>1$, we have $b\simeq\lambda$, thus
$4{\lambda\over b^2 c_\Omega}\simeq 4{1\over\lambda c_\Omega}$, thus
$-1+\sqrt{1+4{\lambda\over b^2 c_\Omega}}\simeq {2\over\lambda c_\Omega}$, so
$\kappa\simeq{1\over c_\Omega}$ in the vicinity of $\lambda=\infty$.

Continuity: Easy check.
\findem

\debrem
For the associated numerical approximation with finite elements,
it $\lambda$ is ``large'' then difficulties are expected, since $C=O(\lambda)$.
Indeed, the conditioning of the associated matrix is $\simeq{C\over\alpha_\Phi}=0(\lambda)$,
and this conditioning explodes with~$\lambda$.
However, a bad choice of the discrete spaces leads to a much faster explosion than expected,
see proposition~\ref{propcch}.
\finrem

%%%%%%%%%%%%%%%%%%%%%%%%%%%%%%%%%%%%%%%%%%%%%%%%%%%%%%%%%%%%%%%

\subsection{The coercivity constant for $p$}

For the analysis of the locking phenomenon (due to a ``bad choice'' of the discrete spaces),
let us look at the coercivity constant for~$p$ (where $\lambda$ appears):

\debprop
If $(\vv,q)\in X$ then, with~\eref{eqpoinc},
\be
\label{eqcfp0}
  \Phi((\vv,q),(\vv,q)) \ge  c_\Omega{\lambda\over\lambda+ c_\Omega} ||\vgrad q||_\Ld^2,
\ee
and
\be
\label{eqcfp}
c_\Omega{\lambda\over\lambda+ c_\Omega}\simeq c_\Omega \quad \hbox{as} \quad \lambda \rar \infty.
\ee
\finprop

\debdem
Modification or the previous proof:
$$
  \Phi((\vv,q),(\vv,q)) \ge  c_\Omega||\vv||_\Ld^2
  +\lambda[||\vv||_\Ld^2+||\vgrad q||_\Ld^2-2||\vv||_\Ld||q||_\Ld]
  \ge \alpha_p||\vgrad q||_\Ld^2,
$$
and the largest $\alpha_p$ possible is $\alpha_p= c_\Omega{\lambda\over\lambda+ c_\Omega}$
(has to satisfy ``$ c_\Omega x^2+\lambda(x-y)^2\ge \alpha_py^2$).
\findem

%%%%%%%%%%%%%%%%%%%%%%%%%%%%%%%%%%%%%%%%%%%%%%%%%%%%%%%%%%%%%%%

\subsection{The discrete problem}

\def\Pivh{{\Pi_{V_h}}}
\def\Piqh{{\Pi_{Q_h}}}

Let $V_h\subset\Huzon$ be a finite dimension subspace.
Let $\Pivh : \Ldon \rar V_h$ be the $\dd_\Ld$ projection onto~$V_h$, that is,
$$
\forall \vv\in\Huzon,\quad \forall \vw_h\in V_h, \quad
(\Pivh\vv,\vw_h)_\Ldo=(\vv,\vw_h)_\Ldo.
$$
%En particulier, pour tout $\vw\in V$, on~a $(\vv-\Pivh\vv,\Pivh\vw)_\Ldo=0$.

Let $Q_h\subset\Huzo$ be a finite dimension subspace.
\comment{
Let $\Piqh : \Ldo \rar Q_h$ be the $\dd_\Ld$ projection onto~$Q_h$, that is,
$$
\forall p\in\Huzo,\quad \forall q_h\in V_h, \quad
(\Piqh p,q_h)_\Ldo=(p,q_h)_\Ldo.
$$
}

Let $X_h=V_h\times Q_h$.
The discrete problem associated to~\eref{eqloc2} is: Find $(\vu_h,p_h)\in X_h$ \st
\be
  \Phi((\vu_h,p_h),(\vv_h,q_h))=(\vf,\vv_h)+(g,q_h),\qquad\forall (\vv_h,q_h)\in X_h.
\ee

\debprop
\label{propcch}
If $(\vv_h,q_h)\in X_h$ then
\be
\label{eqpropcch}
  \Phi((\vv_h,q_h),(\vv_h,q_h)) \ge  c_\Omega{\lambda\over\lambda+ c_\Omega} ||\vgrad q_h||_\Ld^2
  + \lambda{\lambda\over\lambda+ c_\Omega}||\vgrad q_h-\Pivh\vgrad q_h||_\Ldo^2,
\ee
to be compared with~\eref{eqcfp0}.

Illustration:
If $V_h$ is ``small relatively to~$Q_h$'' so that for some $q_h$ the real
$||\vgrad q_h-\Pivh\vgrad q_h||_\Ldo$ does not vanish (fast enough with~$h$)
then the right hand side of~\eref{eqpropcch} increases with~$\lambda$,
to compare with~\eref{eqcfp}.
And the solution $(u_n,p_h)\in X_h$ is bounded by the inverse constant that decreases with~$\lambda$,
thus $(u_h,p_h)$ decreases to zero as $\lambda$ increases:
We get the ``locking'' phenomenon.
\finprop

\debdem
Let $(\vu_h,p_h)$ and $(\vv_h,q_h)\in X_h$. Then
$$
\eqalign{
\Phi((\vu_h,p_h),(\vv_h,q_h))
= & (\grad\vu_h,\grad\vv_h)_\Ld+\lambda(\vgrad p_h - \vu_h , \vgrad q_h - \vv_h )_\Ld \cr
= & (\grad\vu_h,\grad\vv_h)_\Ld
+\lambda(\vgrad p_h-\Pivh\vgrad p_h,\vgrad q_h-\Pivh\vgrad q_h)_\Ld
	\cr
&\quad  
  +\lambda(\Pivh\vgrad p_h - \vu_h , \Pivh\vgrad q_h - \vv_h )_\Ld.
}
$$
Thus
$$
  \Phi((\vv_h,q_h),(\vv_h,q_h))\ge  c_\Omega{\lambda\over\lambda+ c_\Omega} ||\Pivh\vgrad q_h||_\Ld^2
 + \lambda||\vgrad q_h-\Pivh\vgrad q_h||_\Ld^2,
$$
see previous~\S computation.
And Pythagoras give~\eref{propcch}.
\findem

\debrem
The term $\vgrad q_h-\Pivh\vgrad q_h$ is also a problem for the Stokes equations, see~\S~\ref{seccp}.
\finrem

%%%%%%%%%%%%%%%%%%%%%%%%%%%%%%%%%%%%%%%%%%%%%%%%%%%%%%%%%%%%%%%

\subsection{An optimal correction}

Due to the choice of $V_h$ and~$Q_h$, we eventually have too much of coercivity, \cf~\eref{propcch},
so we decide to get rid of it.
That is, we modify~\eref{eqloc0} to get:
Find $(\vu,p)\in V_h\times Q_h$ realizing the minimum of
\be
\eqalign{
   M_h(\vv,q)
 =&\demi||\grad\vv_h||_\Ld^2
  + {\lambda\over2}( ||\vv_h-\vgrad q_h||_\Ld^2
 - \lambda{\lambda\over\lambda+ c_\Omega}||\vgrad q_h-\Pivh\vgrad q_h||_\Ld^2) \cr
   &-(\vf,\vv_h)_\Ld -(g,q_h)_\Ld.
}
\ee
Thus $\Phi_h$ has been transformed into
$$
\eqalign{
  \Phi_h((\vu_h,p_h),(\vv_h,q_h))
=&(\grad\vu_h,\grad\vv_h)_\Ld + \lambda(\vu_h-\vgrad p_h,\vv-\vgrad q_h)_\Ld\cr
 &- {\lambda^2\over\lambda+ c_\Omega}(\vgrad p_h-\Pivh\vgrad p_h,\vgrad q_h-\Pivh\vgrad q_h)_\Ld,
}
$$
and the problem reads:
Find $(\vu_h,p_h) \in V_h\times Q_h$ \st, for all $(\vv_h,q_h) \in V_h\times Q_h$,
$$
   \Phi_h((\vu_h,p_h),(\vv_h,q_h)) = (\vf,\vv_h)_\Ld + (g,q_h)_\Ld.
$$
To solve this problem, we need to compute $\Pivh\vgrad p_h$:
If the $V_h=P_1$-continuous finite elements is made, the computation
amounts to inverse a diagonal matrix, thanks to the mass-lumping technique,
thus is costless.

Computation: we have to compute $(\vu_h,p_h)\in V_h\times Q_h$ \st,
for all $(\vv_h,q_h)\in V_h\times Q_h$,
$$
\left\{\eqalign{
 & (\grad\vu_h,\grad\vv_h)_\Ld + \lambda(\vu_h-\vgrad p_h,\vv_h)_\Ld = (\vf,\vv_h)_\Ld,\cr
 & {-}\lambda(\vu_h-\vgrad p_h,\vgrad q_h)_\Ld
 - \lambda{\lambda\over\lambda+ c_\Omega}(\vgrad p_h{-}\Pivh\vgrad p_h,\vgrad q_h)_\Ld
=(g,q_h)_\Ld.
}\right.
$$
Introducing $\vw_h=\Pivh\vgrad p_h$, we have to find
$(\vu_h,p_h,\vw_h)\in V_h\times Q_h\times V_h$ \st,
for all $(\vv_h,q_h)\in V_h\times Q_h$,
\be
\left\{\eqalign{
 & (\grad\vu_h,\grad\vv_h)_\Ld + \lambda(\vu_h,\vv_h)_\Ld
  -\lambda(\vgrad p_h,\vv_h)_\Ld = (\vf,\vv_h)_\Ld,\cr
 & {-}\lambda(\vu_h,\vgrad q_h)_\Ld
+{ c_\Omega\lambda\over\lambda+ c_\Omega}(\vgrad p_h,\vgrad q_h)_\Ld
 + \lambda{\lambda\over\lambda+ c_\Omega}(\vw_h,\vgrad q_h)_\Ld=(g,q_h)_\Ld,\cr
 & (\vgrad p_h,\vw_h')_\Ld-(\vw_h,\vw_h')_\Ld = 0.
}\right.
\ee
This method gives optimal convergence results.
\Eg, for $P_1$-continuous finite elements of $\vu_h$ and $p_h$,
we get an $O(h)$ convergence.

%%%%%%%%%%%%%%%%%%%%%%%%%%%%%%%%%%%%%%%%%%%%%%%%%%%%%%%%%%%%%%%

\subsection{Classical treatment of the locking}

%%%%%%%%%%%%%%%%%%%%%%%%%%%%%%%%%%%%%%%%%%%%%%%%%%%%%%%%%%%%%%%

\subsubsection{Initial problem}

\def\vgamma{{\vec\gamma}}
\def\vdelta{{\vec\delta}}

See \eg\ Chapelle~\cite{chapelle}, Brezzi and Fortin~\cite{brezzi-fortin}.
The variable
\be
\label{eqvg}
\vgamma = \lambda(\vu-\vgrad p)
\ee
is introduced. Then problem~\eref{eqloc1} becomes:
Find $(\vu,p,\vgamma)\in\Huzo^n\times\Huzo\times Y$ s.t.
\be
\label{eqlocc}
\left\{\eqalign{
 & a((\vu,p),(\vv,q)) + b((\vv,q),\vgamma) = (\vf,\vv)_\Ld + (g,q)_\Ld,\qquad\forall (\vv,q)\in X,\cr
 & b((\vu,p),\vdelta) - {1\over\lambda}\la \vdelta,\vgamma \ra_{\Hmu,\Huz}=0,\qquad\forall \vdelta\in Y,\cr
}\right.
\ee
with
\be
\label{eqY}
\left\{\eqalignrllq{
 & Y = \{\vdelta\in(\Hmuo)^n : \dvg\vdelta\in\Hmuo\},\span\cr
 & a\dd: X\times X\rar\RR : & a((\vu,p),(\vv,q))=(\grad\vu,\grad\vv)_\Ld,\cr
 & b\dd: X\times Y\rar\RR : & b((\vu,p),\vdelta) %=\la\vdelta , \vu-\vgrad p\ra
        =\la \vdelta, \vu\ra_{\Hmu,\Huz} - \la \dvg\vdelta , p\ra_{\Hmu,\Huz}.\cr
}\right.
\ee
(And~\eref{eqlocc}$_2$ gives $\vu-\vgrad p - {1\over\lambda}\vgamma=0$.)
And it is shown that $Y$ is a Banach space for the norm
$$
  ||\vdelta||_Y\eqdef ||\vdelta||_{\Hmuo^n} + ||\dvg\vdelta||_\Hmuo.
$$

\debrem
$Y$ is the space corresponding to ${1\over\lambda}=0$
(\ie\ $\lambda$ ``infinitely large''), \cf\ Kirchhoff--Love shell model:
\be
\left\{\eqalign{
 & a((\vu,p),(\vv,q)) + b((\vv,q),\vgamma) = (\vf,\vv)_\Ld + (g,q)_\Ld,\qquad\forall (\vv,q)\in X,\cr
 & b((\vu,p),\vdelta)=0,\qquad\forall \vdelta\in Y,\cr
}\right.
\ee
to be compared with~\eref{eqlocc}.

For a discretization with finite elements, finite dimensional spaces are often chosen \st\
$V_h\subset\Ldon$, $Y_h\subset\Ldon$ and $Q_h\subset C^0(\Omega;\RR)$;
And then~\eref{eqlocc} is meaningful in $V_h\times Q_h\times Y_h$
with $b((\vu_h,p_h),\vdelta_h) = (\vdelta_h , \vu_h-\vgrad p_h)_\Ldo$.
\finrem

Let $B:X\rar Y'$ be the operator associated to~$b\dd$, that is,
$
\la B(\vv,q),\vdelta\ra_{\Huz,\Hmu}  := b((\vv,q),\vdelta),
$
\ie,
$$
   B(\vv,q)=\vv-\vgrad q.
$$
Then, with Poincaré inequality, it is easy to check that $a\dd$ is coercive on
$\Ker(B)=\{ (\vv,q)\in X : \vv=\vgrad q\}$.

Then is shown that $B$ is surjective (inf-sup condition), that is,
$$
  \exists k>0,\quad \forall \vdelta\in Y,\quad \sup_{(\vv,q)\in X}
   {b((\vv,q),\vdelta) \over||(\vv,q)||_X||\vdelta||_Y}\ge k.
$$
See \eg\ Chapelle~\cite{chapelle}, Brezzi et Fortin~\cite{brezzi-fortin}.

%%%%%%%%%%%%%%%%%%%%%%%%%%%%%%%%%%%%%%%%%%%%%%%%%%%%%%%%%%%%%%%

\subsubsection{discrete problem}

The discrete inf-sup condition has to be satisfied:
this lead to numerous articles.
There are two difficulties:

1- 
An adequat choice of finite element spaces to satisfy the inf-sup condition
(it the stabilization is not used),
la stabilisation de~$\vgamma_h$, ou le choix adequat d'éléments finis compatibles
 pour satisfaire la condition inf-sup,

2-
An adequat choice of finite element spaces to satisfy
the coercivity of~$a\dd$ on the kernel $\Ker(B_h)$ (with $B_h$ the discrete operator).
But this problem can be easily fixed  by modiying~\eref{eqloc0} into
\be
\tilde M(\vv,q)= \demi ||\grad\vv||^2_\Ldo
+ \demi ||\vv-\vgrad q||_\Ldo + {\lambda-1\over2}||\vv-\vgrad q||_\Ldo
   -(\vf,\vv)_\Ldo -(g,q)_\Ldo,
\ee
that is, by replacing $a\dd$, \cf~\eref{eqY}, with
$$
\eqalign{
& \tilde a((\vu,p),(\vv,q)) = (\grad\vu,\grad\vv)_\Ld + (\vu-\vgrad p,\vv-\vgrad q)_\Ld. \cr
}
$$
And we consider~\eref{eqlocc} with $\tilde a\dd$ instead of~$a\dd$.

Now $\tilde a\dd$ is coercive sur~$X$ (thanks to~\eref{eqpoinc}),
thus on $\Ker(B)=\{(\vv,q) : \vv=\vgrad q\}$,
and we get a similar problem to the Stokes problem
(choice of adequat finite element spaces, or choice of a stabilization).

%%%%%%%%%%%%%%%%%%%%%%%%%%%%%%%%%%%%%%%%%%%%%%%%%%%%%%%%%%%%%%%
%%%%%%%%%%%%%%%%%%%%%%%%%%%%%%%%%%%%%%%%%%%%%%%%%%%%%%%%%%%%%%%

\newpage

\section{Weak Dirichlet condition}

(See Babu\v ska~\cite{babuska} for the initial manuscript.)

\def\dird{{d}}

%L'article de référence est celui de Babu\v ska~\cite{babuska}.

%Voir également Pitkäranta~\cite{pitkaranta}, Stenberg~\cite{stenberg1}, Barbosa et Hughes~\cite{barbosa-hughes}.

%%%%%%%%%%%%%%%%%%%%%%%%%%%%%%%%%%%%%%%%%%%%%%%%%%%%%%%%%%%%%%%

\subsection{Initial problem}

Let $f\in\Ldo$ and $\dird\in\Hdemig$.
Let $u_\dird\in\Huo$ \st\ $u_\dird{}_{|\Gamma}=\dird$;
Such a function $u_d$ exists since the trace operator
$\Gamma_0 : \Huo \rar \Ldg$ is surjective onto~$\Gamma$,
cf.~\eref{eqngammaz}.
Let $\dird+\Huzo := u_\dird+\Huzo = \{u_\dird+v,\; v\in\Huzo\}$, affine space in~$\Huo$
independent of the choice of $u_\dird$ the reverse image of~$\dird$ by~$\gamma_0$ (trivial).

Consider the problem: Find $u\in \dird+\Huzo$ \st
\be
\label{eqcld1}
\left\{\eqalign{
 & -\Delta u + u = f\quad\hbox{dans }\Omega,\cr
 & u_{|\Gamma}=\dird\quad\hbox{sur }\Gamma.\cr
}\right.
\ee
Thanks to the Lax--Milgram theorem, this problem is well-posed.

%%%%%%%%%%%%%%%%%%%%%%%%%%%%%%%%%%%%%%%%%%%%%%%%%%%%%%%%%%%%%%%

\subsection{Mixed problem}

The aim is to impose the Dirichlet condition with a Lagrangian multiplier.
So the problem becomes:
Find $(u,\lambda)\in \Huo\times\Hmdemig$ \st
\be
\label{eqluld2}
\left\{\eqalign{
 & (u,v)_\Huo + \la \lambda,v\ra_{\Hmdemig,\Hdemig} = (f,v)_\Ldo,
           \quad\forall v\in\Huo,\cr
 & \la u,\mu\ra_{\Hdemig,\Hmdemig} = \la \dird,\mu\ra_{\Hdemig,\Hmdemig},
           \quad\forall\mu\in\Hmdemig.\cr
}\right.
\ee
If $(u,\lambda)$ exists in $\Huo\times\Hmdemig$, then we get:
\be
\label{eqluld2i}
\left\{\eqalign{
 & -\Delta u + u = f \;\in \Ldo,\cr
 & u = \dird \; \in \Hdemig,\cr
 & \lambda = -{\pa u\over\pa n} \; \in \Hdemig.\cr
}\right.
\ee
Interpretation of the Lagrangian multiplier:
$\lambda$ is, up to the sign, the force $\vgrad u.\vn$ needed on~$\Gamma$ for $u$ to stay equal to $\dird$ on~$\Gamma$.

With
\be
\left\{\eqalign{
& a(u,v) %=(\vgrad u,\vgrad v)_\Ldo + (u,v)_\Ldo\quad(
=(u,v)_\Huo,\cr
& b(v,\lambda) %= \la \gamma_0(v),\lambda\ra_{\Hdemig,\Hmdemig} \eqnote
= \la v,\lambda\ra_{\Hdemig,\Hmdemig},\cr
}\right.
\ee
\eref{eqluld2} has the appearance of~\eref{eqgenw0} with $V=\Huo$ and $Q=\Hmdemig$.
And $a\dd$ is bilinear (trivial) continuous and coercive (it is the $V=\Huo$-inner product),
and $b\dd$ is bilinear (trivial) continuous since
$|b(v,\lambda)| \le ||v||_\Hdemig ||\lambda||_\Hmdemig$ and $\gamma_0$ is continuous
(so $||v||_\Hdemig \le ||\gamma_0||\,||v||_\Huo$).

We have %$V'= \Huo'$ and 
$Q' = \Hdemig$,
%so $b(v,\lambda) = \la Bv,\lambda \ra_{\Hdemig,\Hmdemig} = \la B^t\lambda,v\ra_{\Huo',\Huo}= \la v,\lambda\ra_{\Hdemig,\Hmdemig}$, 
so
$B :
\left\{\eqalign{
\Huo &\rar \Hdemig \cr
v & \rar Bv = \gamma_0(v)
}\right\}
$
is linear continuous (since $b\dd$ is bilinear continuous) and surjective (definition of~$\Hdemig$),
thus $\Im(B^t)$ is closed in $ V'=\Huo'$,
with $B^t :
\left\{\eqalign{
\Hmdemig &\rar \Huo' \cr
\lambda & \rar B^t\lambda
}\right\}
$
defined by $\la B^t\lambda,v\ra_{\Huo',\Huo}= \la v,\lambda\ra_{\Hdemig,\Hmdemig}$.
Thus
\be
\label{eqlambda}
\exists k>0,\; \forall \lambda\in\Hmdemig,\; ||B^t\lambda||_{\Huo'} \ge k\,||\lambda||_{\Hmdemig/\Ker B^t}.
\ee
(That is,
$\ds
\exists k>0,\; \inf_{\lambda\in\Hmdemig}
\sup_{v\in \Huo} |b({v\over ||v||_\Huo||},{\lambda \over \lambda||_\Hmdemig})| \ge k
%    \sup_{v\in \Huo} |b({v\over ||v||_\Huo},{\lambda\over||\lambda||_\Hmdemig})| \ge k,
$.)

\debrem
\label{remli}
The computation of $\lambda$ may give disappointing results since the control for $\lambda$
is done with the $||.||_\Hmdemig$-norm, cf.~\eref{eqlambda} (not even a $||.||_\Ldg$ control).
So numerical problems are expected.
\finrem

\debrem
This mixed problem leads to ``transmission problems'' (or hybrid problems)
with ``mortar finite elements'', see Bernardi, Maday, Patera.
\finrem

The associated Lagrangian is (saddle point problem)
\be
\label{eqlul}
L(u,\lambda)
= \demi %||u||_\Huo^2 
(||\vgrad u||_\Ldo^2 + ||u||_\Ldo^2)
+ \la u - \dird,\lambda)_{\Hdemig,\Hmdemig} - (f,v)_\Ldo,
\ee

%%%%%%%%%%%%%%%%%%%%%%%%%%%%%%%%%%%%%%%%%%%%%%%%%%%%%%%%%%%%%%%

\subsection{Discrete problem}

Let $V_h\subset\Huo$ and $\Lambda_h\subset\Hmdemig$ be finite dimension spaces.
The discrete problem relative to~\eref{eqluld2} is:
Find $(u_h,\lambda_h)\in V_h\times\Lambda_h$ \st
\be
\label{eqluld2h}
\left\{\eqalignrll{
 & (u_h , v_h)_\Huo  + (v_h,\lambda_h)_\Ldg
    &= (f,v_h)_\Ldo,
           \qquad\forall v_h\in V_h,\cr
 & (u_h,\mu_h)_\Ldg &= (\dird,\mu_h)_\Ldg,
           \qquad\forall\mu_h\in\Lambda_h.\cr
}\right.
\ee

%%%%%%%%%%%%%%%%%%%%%%%%%%%%%%%%%%%%%%%%%%%%%%%%%%%%%%%%%%%%%%%

\subsection{Finite elements $P_k{-}C^0$: unstable}

%%%%%%%%%%%%%%%%%%%%%%%%%%%%%%%%%%%%%%%%%%%%%%%%%%%%%%%%%%%%%%%

\subsubsection{The discrete inf-sup condition}

Consider the Lagrange finite elements
$\Vh=P_k{-}C^0$ in~$\Omega$ and $\Lambda_h = \gamma_0(V_h)$ $P_k{-}C^0$ on~$\Gamma$.
The discrete (trace) operator
$
B_h = \gamma_0{}_{|\Gamma} :
\left\{\eqalign{
(V_h,||.||_\Huo) & \rar (\Lambda_h,||.||_\Hdemig) \cr
v_h & \rar \gamma_0(v_h) = v_{h|\Gamma}
}\right\}$
is continuous and surjective (trivial here with $P_k$-$C^0$ finite elements for both $\Vh$ and~$\Lambda_h$),
thus \eref{eqlambda} holds with some $k_h$ instead of~$k$,
and $Q_h$ instead of~$Q$.
But the control on~$\lambda_h$ is very weak (a $\Hmdemig$ control),
and $k_h$ a priori depends on~$h$.

%%%%%%%%%%%%%%%%%%%%%%%%%%%%%%%%%%%%%%%%%%%%%%%%%%%%%%%%%%%%%%%

\subsubsection{Barbosa et Hughes}

(See Barbosa and Hughes~\cite{barbosa-hughes}, Pitk\"aranta~\cite{pitkaranta},
Stenberg~\cite{stenberg1}).

A finite element mesh~$\calT_h$ is defined in~$\Omega$,
and the trace of this mesh on~$\Gamma$ will be used as a mesh on~$\Gamma$.

Barbosa and Hughes stabilize the Lagrangian multiplier~$\lambda$ with its value
$\lambda = -{\pa u\over\pa n}$, \cf~\eref{eqluld2i}.
Thus the problem now reads:
Find $(u_h,\lambda_h)\in V_h\times\Lambda_h$ \st,
for all $(v_h,\mu_h)\in V_h\times\Lambda_h$,
\be
\label{eqluld2hs}
\left\{\eqalign{
 & (u_h , v_h)_\Huo + \int_\Gamma v_h\lambda_h\,d\Gamma
   - \alpha h \int_\Gamma (\lambda_h{+}{\pa u_h\over\pa n}){\pa v_h\over\pa n}\,d\Gamma
    = (f,v_h)_\Ldo,\cr
 & \int_\Gamma u_h\mu_h\,d\Gamma
  -\alpha h\, \int_\Gamma (\lambda_h{+}{\pa u_h\over\pa n})\mu_h\,d\Gamma
   = \int_\Gamma \dird\mu_h\,d\Gamma,\cr
}\right.
\ee
with $\alpha$ a constant to be chosen,
corresponding to the saddle point of the modified Lagrangian
\be
\label{eqlulh}
  L_h(u,\lambda)=L(u,\lambda)
    - \alpha\, h\,||\lambda{+}{\pa u\over\pa n}||_\Ldg^2,
\ee
cf.~\eref{eqlul}. See Stenberg~\cite{stenberg1}.

We then get the ``penalized'' problem, written here as the matrix problem
\be
\label{eqpbmat2p}
  \pmatrix{A & B^t\cr B&-\alpha C}.\pmatrix{\vu\cr p}=\pmatrix{\vf\cr\vg}.
\ee

\debthm
(Barbosa and Hughes~\cite{barbosa-hughes}, Pitk\"aranta~\cite{pitkaranta}.)
The mesh~$\calT_h$ is supposed to be quasi-uniform, that is,
the following inverse inequality is true:
\be
\label{eqestinvd2}
   \exists C_i>0,\quad \forall v_h\in V_h,\quad
  h^\demi\; ||{\pa v_h\over\pa n}||_\Ldg\le C_i\;||\vgrad v_h||_\Ldo.
\ee
And $\alpha$ is supposed small enough (not to destroy the coercivity for~$u$), namely:
\be
   0<\alpha< {1\over C_i}.
\ee
Then the stabilized problem~\eref{eqluld2hs} is well posed, and
for $P_k$-$C^0$ finite elements, as soon as the exact solution $u$ is in $H^{k+1}(\Omega)$,
and we get the usual a priori estimate
$$
   ||u-u_h||_\Huo \le C h^k||u||_{H^{k+1}(\Omega)},
$$
$C$ being a constant independent of~$h$.
\finthm

\debdem
See Barbosa--Hughes~\cite{barbosa-hughes} and Stenberg~\cite{stenberg1}.
\findem

\debrem
The inequality~\eref{eqestinvd2} also reads
$$
 h\int_\Gamma (\vgrad v.\vn)^2\,d\Gamma \le C_i^2 \;\int_\Omega ||\vgrad v||^2_\RRn\;d\Omega,
$$
where $h$ on the left hand side is expected:
$(h\int_\Gamma)$ has a volume dimension, same dimension as $(\int_\Omega)$,
for a quasi-uniform mesh.
\finrem

%%%%%%%%%%%%%%%%%%%%%%%%%%%%%%%%%%%%%%%%%%%%%%%%%%%%%%%%%%%%%%%

\subsubsection{Multiplier elimination: Nitsche method}

\def\Lambdah{{\Lambda_h}}

Stenberg~\cite{stenberg1} has shown that Barbosa and Hughes~\cite{barbosa-hughes} method
is equivalent to Nitsche~\cite{nitsche} method when
$\Lambda_h=P_0(\Gamma_h)$ and $V_h=P_1(\Omega_h)$
(when the mesh on~$\Gamma$ is the trace of the mesh in~$\Omega$):
Find $u_h\in V_h$ \st, for all $v_h\in V_h$,
$$
(u_h,v_h)_\Huo -\la{\pa u_h\over\pa n},v_h\ra_\Gamma -\la{\pa v_h\over\pa n},u_h-\dird\ra_\Gamma
+ \gamma\sum_{E\in{\cal E}_h} {1\over h_E}\la u_h-\dird,v_h\ra_E = (f,v_h)_\Ldo ,
$$
for some $\gamma>0$,
\ie, find $u_h\in V_h$ \st, for all $v_h\in V_h$,
\be
\label{eqnitsche}
\left\{\eqalign{
 & (u_h,v_h)_\Huo -\la{\pa u_h\over\pa n},v_h\ra_\Gamma -\la{\pa v_h\over\pa n},u_h\ra_\Gamma
  + \gamma\sum_{E\in{\cal E}_h} {1\over h_E}\la u_h,v_h\ra_E\cr
 & \qquad \qquad \qquad \qquad = (f,v_h)_\Ldo -\la{\pa v_h\over\pa n},\dird\ra_\Gamma
  + \gamma\sum_{E\in{\cal E}_h} {1\over h_E}\la \dird,v_h\ra_E.\cr
}\right.
\ee
We then get $u_\Gamma=\dird$.
This method is simpler to compute since no Lagrangian multiplier intervenes.

\debprop
If~\eref{eqestinvd2}, if $\gamma>C_i$, if $V_h=P_k$-$C^0$ and $u\in H^{k+1}(\Omega)$,
then~\eref{eqnitsche} gives the usual result:
$$
   ||u-u_h||_\Huo\le C h^k||u||_{H^{k+1}(\Omega)}.
$$
\finprop

\debdem
See Stenberg~\cite{stenberg1}.
\findem

\noindent
{\bf Comparison of the method of Nitsche with the method of Barbosa and Hughes :}
~\eref{eqluld2hs}$_2$ gives
$$
 \lambda_h=-\Pi_\Lambdah({\pa u_h\over\pa n}) +{1\over\alpha h}( u_h{-}\dird).
$$
Thus~\eref{eqluld2hs}$_1$ becomes
$$
\eqalign{
& (u_h,v_h)_\Huo 
-\int_\Gamma \Pi_\Lambdah({\pa u_h\over\pa n})v_h\,d\Gamma
 - \int_\Gamma \Pi_\Lambdah({\pa v_h\over\pa n})u_h\,d\Gamma
 + {1\over\alpha h} \int_\Gamma  u_hv_h\,d\Gamma \cr
&\qquad\qquad\qquad\qquad\qquad\qquad\qquad\qquad -\alpha h (\int_\Gamma {\pa u_h\over\pa n}{\pa v_h\over\pa n}
         -\Pi_\Lambdah{\pa u_h\over\pa n}\Pi_\Lambdah{\pa v_h\over\pa n}\,d\Gamma)\cr
&\qquad = (f,v_h)_\Ldo -\int_\Gamma g\,\Pi_\Lambdah{\pa v_h\over\pa n}\,d\Gamma
  +  {1\over \alpha h}\int_\Gamma g v_h\,d\Gamma,
}
$$
With $\Lambda_h=P_0$ and $V_h=P_1$ we then get
$\Pi_\Lambdah({\pa u_h\over\pa n})={\pa u_h\over\pa n}$,
and then~\eref{eqnitsche}.

\newpage

\part{Theory}

Most of the results can be found in Brézis~\cite{brezis}.

\section{The open mapping theorem}

%%%%%%%%%%%%%%%%%%%%%%%%%%%%%%%%%%%%%%%%%%%%%%%%%%%%%%%%%%%%%%%%%%%%%%%%%%%%%%%%%%%

\subsection{Notations}

If $E$ and $F$ are linear spaces, a map $T:E\rar F$ is linear iff
$T(x_1+\lambda x_2)= T(x_1) + \lambda T(x_2)$ for all $x_1,x_2\in E$ and all $\lambda \in \RR$;
And then $T(x)$ is denoted $T.x$ or~$Tx$.

Let $(E,||.||_E)$ be a normed space.
Let $B_E(x,\rho)=\{ x'\in E ; ||x'-x||_E < \rho\}$ the ball of radius $\rho>0$ centered at $x\in E$.

If $(E,||.||_E)$ and $(F,||.||_F)$ are two normed spaces,
if $T :
\left\{\eqalign{
E & \rar F \cr
x & \rar T(x)
}\right\}
$
is a linear map, then $T$ is said to be continuous (or bounded) iff
\be
\label{eqTc}
\exists c>0,\quad \forall x\in E,\quad ||T.x||_F \le c\,||x||_E.
\qandthen ||T|| := \sup_{x\in B_E(0,1)} ||T.x||_F.
\ee

In the sequel, the space $E$ and $F$ will be Banach spaces (complete for the norm in use).

Let $\calL(E;F)$ be the set of linear continuous mapping from $E$ to~$F$.
Then
\be
||.||_{\calL(E;F)} :
\left\{\eqalign{
\calL(E;F) & \rar \RR \cr
T & \rar ||T||_{\calL(E;F)} := \sup_{x\in B_E(0,1)} ||T.x||_F \eqnote ||T||
}\right.
\ee
define a norm in~$\calL(E;F)$ (easy check),
and $(\calL(E;F),||.||)$ is a Banach space
(check: If $(T_n)_\NNs$ is a Cauchy sequence, that is $||T_n-T_m|| \matrarrow_{n,m\rar\infty}0$,
then, for any $x\in E$, $||(T_n - T_m)(x)||_F \matrarrow_{n,m\rar\infty}0$,
thus $(T_n(x))_\NNs$ is a Cauchy sequence in~$F$ complete, thus converge to a $y_x\in F$;
Then define $T : x \in E \rar T(x)=y_x$: It is easy to check that $T$ is linear and continuous
with $||T-T_n|| \matrarrow_{n\rar \infty}0$.)

Let $E' := \calL(E;\RR)$, called the dual of~$E$ (the set of linear continuous real valued functions,
$\RR$ being provided with its usual norm).
For $\ell \in E'$ and $x\in E$ denote:
\be
\ell(x) = \ell.x = \la\ell,x\ra_{E',E} \;\; \in \RR.
\ee
So, \cf~\eref{eqTc},
%Then the norm of $\ell\in\Es$ is given by
\be
\label{eqdefnell}
||\ell||_{E'}
= \sup_{x\in B_E(0,1)} |\ell.x|
= \sup_{x\in B_E(0,1)} |\la\ell,x\ra_{E',E}|
\ee
defines a norm in~$E'$ \st\ $(E',||.||_{E'})$ is a Banach space.

If $T\in \calL(E;F)$ (linear and continuous) then its adjoint is the
linear map $T':F'\rar E'$ characterized by:
\be
\label{eqdefadj}
T' : 
\left\{\eqalign{
F' & \rar E' \cr
\ell & \rar T'(\ell) \eqnote T'.\ell, \qwhere
\la T'.\ell,x \ra_{E',E} := \la \ell,T.x \ra_{F',F}, \quad \forall x\in E.
}\right.
\ee

\debprop
$T'$ is continuous with
\be
\label{eqnTp}
||T'|| = ||T||.
\ee
Thus $T' \in \calL(F',E')$.
\finprop

\debdem
$||T'.\ell||_{E'}
= \sup_{||x||_E\le 1} |\la T'.\ell , x \ra_{E',E}|
= \sup_{||x||_E\le 1} |\la \ell , T.x \ra_{F',F}|
\le \sup_{||x||_E\le 1} ||\ell||_{F'}||T.x||_F
= \sup_{||x||_E\le 1} ||T||\,||x||_E||\ell||_{F'}
= ||T||\,||\ell||_{F'}
$, thus $||T'|| \le ||T||$;
And similarly $||T.x||_F \le ||T'||\,||\ell||_{F'}$, thus $||T|| \le ||T'||$.
\findem

$\Ess = (\Es)' = \calL(E';\RR)$ is a Banach space (since $\RR$ is complete).
Let
\be
\label{eqreflex}
J : 
\left\{\eqalign{
E & \rar E'' = \calL(E';\RR) \cr
x & \rar J(x),\qwhere J(x)(\ell) := \ell.x,\; \forall x\in E.
}\right.
\ee
$J$ is linear (trivial),
is continuous, with $||J||
= \sup_{x\in B_E(0,1)} |J(x)|
= \sup_{x\in B_E(0,1)} ( \sup_{\ell\in B_\Es(0,1)}|J(x)(\ell)|)
= \sup_{x\in B_E(0,1)} ( \sup_{\ell\in B_\Es(0,1)}|\ell.x|)
= \sup_{x\in B_E(0,1)} ( ||x||_E) = 1
$,
and injective (one-to-one) since $J(x)=0$ implies $\ell.x=0$ for all $\ell\in E'$ that implies $x=0$.
%thanks to the Hahn--Banach theorem, see Brézis~\cite{brezis}.
Thus $J$ is a ``canonical injection''.

Thus $J(E)=\Im(E)$, the range or image of~$E$ by~$J$, can be identified to a subspace of~$E''$.

\debdef
\label{defrefl}
A Banach space $E$ is reflexive iff $J$ is bijective (= one-to-one and onto),
and then is identified with~$E$, denoted $E'' \simeq E$, and $J(x)$ is denoted~$x$.
\findef

(Remark: A Hilbert space is always reflexive,
and a reflexive Banach space ``almost'' behaves like a Hilbert space for computation purposes
(with the use of the bracket $\la.,.\ra_{E',E}$ similar to the use of a inner product).
There are however some substantial differences: \eg\ in a reflexive Banach space there exist closed subspaces without any complement, whereas in a Hilbert space any closed subspace has a complement (even an orthogonal one);
And this causes some theoretical difficulties treated in the sequel.)

\def\tT{{\widetilde T}}

%%%%%%%%%%%%%%%%%%%%%%%%%%%%%%%%%%%%%%%%%%%%%%%%%%%%%%%%%%%%%%%%%%%%%%%%%%%%%%%%%%%

\subsection{The open mapping theorem}
\label{secapo}

\begin{theorem}[Open mapping theorem]
\label{thmomt0}
Let $E$ and $F$ be Banach spaces.
If $T\in \calL(E;F)$ (linear and continuous) is surjective (= onto,
\ie\ $\Im(T)=F$), then
\be
\label{eq0ai}
\exists \gamma >0 \qst T(B_E(0,1)) \supset B_F(0,\gamma).
\ee
That is, if $T$ is linear continuous and surjective,
then any open set in~$E$ is transformed by~$T$ into an open set in~$F$.
So $T(B_E(0,1))$ is not ``flat'' (it contains an open set).

And the converse is true: if~\eref{eq0ai} then $T$ is surjective.
\finthm

\debdem
See Brézis~\cite{brezis}. Steps :
1- $T$ being onto, we have $\bigcup_{n\in\NNs} \overline{T(B_E(0,n))}=F$,
and Baire's Theorem gives the existence of a closed space $\overline{T(B_E(0,n))}$
containing an open set;
2- The linearity of~$T$ then implies that $\overline{T(B_E(0,1))}$ contains an open set $B_F(0,2\gamma)$
for some $\gamma>0$.
3- And $T$ being continuous and $E$ being complete we get
$T(B_E(0,1)) \supset B_F(0,\gamma)$.

Converse: $T(B_E(0,1))\supset B_F(0,\gamma)$, and $T$ is linear, so $T(E)=F$.
\findem

\debcor
\label{thmomt}
%Let $E$ and $F$ be banach spaces.
If $T\in \calL(E;F)$ is bijective, \ie\ injective (= one-to-one) and surjective (= onto),
then the linear map $T^{-1}:F \rar E$ is continuous, that is,
\be
\label{eq1ai}
\exists \gamma>0,\;
\forall y\in F,\quad ||T^{-1}.y||_E \le {1\over \gamma}||y||_F.
\ee
%and $||T^{-1}||_{\calL(F;E)} \le {1\over \gamma}$,
Thus
\be
\label{eq2ai}
\exists \gamma>0 ,\; 
\forall x\in E,\quad ||T.x||_{F} \ge \gamma ||x||_E.
\ee
\fincor

\debdem
Then $T$ bijective gives $T^{-1}(B_F(0,\gamma)) \subset B_E(0,1)$.
And $T^{-1}$ is linear since $T$ is,
thus $T^{-1}(B_F(0,1)) \subset B_E(0,{1\over \gamma})$.
So $y\in B_F(0,1)$ gives $||T^{-1}.y||_E \le {1\over \gamma}||y||_F$, \ie\ \eref{eq1ai}.
%Vrai pour tout $y\in B(0,1)$, donc $||T^{-1}||_{\calL(F;E)} \le {1\over \gamma}$.
Then $y=T.x$ gives~\eref{eq2ai} (bijectivity).
\findem

\debrem
If $T$ is bijective between Banach spaces,
then the problem: For $b\in F$ find $x\in E$ \st\ $T.x=b$ is well-posed,
that is, has a unique solution $x=T^{-1}.b$ \st\
$\exists c>0$ (independent of~$b$), $||x||_E \le c ||b||_F$
(the inverse $T^{-1}$ is continuous).
Indeed, the bijectivity of~$T$ gives a unique solution $x=T^{-1}.b$,
and~\eref{eq1ai} gives $||x||_E = ||T^{-1}.b||_E \le {1\over \gamma} ||b||_F$.
\finrem

\debrem
A linear continuous bijective mapping between two infinite dimensional Banach spaces behaves like
in finite dimension,
\eg\ like $T:\RR^2\rar\RR^2$ given by its matrix $\pmatrix{-2&0\cr0&3}$.
Here $||T||=3$,
$T^{-1}=\pmatrix{-\demi&0\cr0&{1\over3}}$,
and $||T^{-1}||=\demi = {1\over \gamma}$.
\finrem

\debrem
\label{exas1}
\def\elld{{\ell^2}}
The bijectivity between Banach spaces (complete spaces) is required:

Let $\ell^2=\{(x_n)_\NNs\in\RR^\NNs : \sum_{n\in\NNs} x_n^2<\infty\}$
(the space of finite energy sequences),
let $E=F=\ell^2$, and let $T : \ell^2\rar\ell^2$
be given by $T((x_n)_\NNs)=({x_n\over n})_\NNs$ for any $(x_n)\in\ell^2$,
that is, with $(e_n)_\NNs$ the canonical basis in~$\elld$, $T.e_n = {1\over n} e_n$
(the associated generalized matrix is the infinite diagonal matrix $\diag(1,\demi,...,{1\over n},...)$.) 	
Here $T$ is injective since $\Ker(T)=\{0\}$ (trivial),
but not surjective since $({1\over n})_\NNs \in \elld$ has no counterimage in~$\ell^2$
(it would be the constant sequence $(1)_\NNs \notin\elld$).
And its range $\RT$ is dense in~$\ell^2$: Indeed, if $(y_n)_\NNs \in \elld$ then let $x_n = ny_n$,
so that $(x_n)_\NNs \in \RR^\NNs$, and for $N\in\NNs$,
define the truncated sequence $(x_n^N)_\NNs$ by $x_n^N = x_n$ if $n\le N$ and
$x_n^N=0$ otherwise; then $(x_n^N)\in\elld$ (trivial) and
$\forall \eps>0$, $\exists N\in\NNs$, 
$||y_n - Tx_n^N||^2_\elld =\sum_{n=N+1}^\infty y_n^2 <\eps$).
Here $\RT$ is not closed in~$\elld$ (since dense and closed would imply $\RT=\elld$),
and $\RT$ is ``flat'' in~$\elld$, that is, $T(B_\elld(0,1))$ does not contain any open ball:
In this example it can be seen with the canonical basis $(e_n)_\NNs$ that verifies $T^{-1}.e_n = n e_n$,
so that $T^{-1}(B_\elld(0,1))$ is not bounded
(if one prefers, $T^{-1}(\demi\gamma e_n) = n \demi\gamma e_n \notin T(B_\elld(0,1))$
as soon as $n > {2\over \gamma}$ although $\demi\gamma e_n \in B_\elld(0,\gamma)$).
\finrem

%%%%%%%%%%%%%%%%%%%%%%%%%%%%%%%%%%%%%%%%%%%%%%%%%%%%%%%%%%%%%%%%%%%%%%%%%%%%%%%%%%%

%\subsection{The open mapping theorem applied when $\RT$ is closed}

\debcor
%Let $E$ and $F$ be banach spaces.
If $T\in \calL(E;F)$ (linear and continuous) is injective (=one-to-one), and if $\RT$ is closed in~$F$, then
% $T\in \calL(E;\RT)$ is bijective and
\be
\label{eq2ai2}
\exists \gamma>0 ,\;
\forall x\in E,\quad ||T.x||_{F} \ge \gamma ||x||_E.
\ee
\fincor

\debdem
For $y\in \RT$ denote
$||y||_\RT := ||y||_F$ for all $y\in \RT$.
So $||.||_\RT$ is a norm in~$\RT$ (trivial).
Then $\RT$ closed in~$F$ implies $(\RT,||.||_\RT) = (\RT,||.||_F)$ is a Banach space denoted~$\RT$.
%Let $||.||_\RT : \RT \rar \RR$ be defined by $||y||_\RT = ||y||_F$. Then $||.||_\RT$ is a norm in~$\RT$, since $||.||_F$ is.
Let
\be
T_R :
\left\{\eqalign{
E & \rar\RT \cr
x & \rar T_R(x) = T(x). \cr
}\right.
\ee
Then $T_R$ is linear continuous bijective between Banach spaces.
Thus \eref{eq2ai} gives
$\exists \gamma>0$, $\forall x\in E$, $||T_R.x||_\RT \ge \gamma ||x||_E$, \ie~\eref{eq2ai2}.
\findem

%%%%%%%%%%%%%%%%%%%%%%%%%%%%%%%%%%%%%%%%%%%%%%%%%%%%%%%%%%%%%%%%%%%%%%%%%%%%%%%%%%%

\subsection{Quotient space $E/\Ker(T)$, and open mapping theorem}

Let $E$ and $F$ be banach spaces.
Let $T\in \calL(E;F)$ (linear and continuous).
Then $K = \Ker(T) = T^{-1}(\{0\})$ (the kernel of~$T$) is linear subspace that is closed (since $T$ is continuous).

Consider the relation in~$E$ defined by: $x \sim y$ iff $x-y\in K$.
This is an equivalence relation (easy check).
Let $\EK = \{ Z\subset E : \exists x\in E,\; Z = x+K\}
=\{x+K : x\in E\}$ be the set of the equivalence classes (quotient space).
An element of~$\EK$ is denoted $\dx=x+K$.
In particular $\dot 0 = K$.

%$K$ being a linear space, $K+K = K$ and $\lambda K=K$ for $\lambda\in K$, and
The (usual) operators $+$ in~$\EK$ and $.$ on~$\EK$ are defined by, if $\dx = x+K$, $\dy = y+K$ and $\lambda\in\RR$,
\be
\dx+\dy=x+y+K, \qand \lambda.\dx = \lambda x+K,
\ee
definition independent of the $x' \in \dx$ and $y'\in \dy$ (easy check).
Then $(\EK,+,.)$ is a linear space (easy check) with $\dot 0$ the zero in~$\EK$.

\deblem
The canonical map $\pi :
\left\{\eqalign{
E & \rar \EK \cr
x & \rar \pi(x) := x+K =\dx \cr
}\right\}
$ is linear and surjective. 
\finlem

\debdem
Linearity: $\pi(x+\lambda y) = (x+\lambda y)+K = x+K + \lambda y + \lambda K
= \pi(x) + \lambda \pi(y)$ since $K$ is a linear space (so $K = K + \lambda K$).

Surjectivity: If $\dx \in \EK$ then $\exists x\in E$ \st\ $\dx = x+K =\pi(x)$ (definition of $\EK$).
\findem

For $\dx \in \EK$, define $||.||_\EK : \EK \rar \RR$ by
\be
\label{eqndx}
|| \dx ||_\EK = ||\pi(x)||_\EK := \inf_{x_0\in K}||x+x_0||_E \eqnote || x ||_\EK.
\ee

\deblem
$|| . ||_\EK$ is a norm in~$\EK$, and $(\EK,||.||_\EK)$ is a Banach space.

And $\pi$ is continuous with $||\pi||\le 1$.
\finlem

\debdem
$||\dx||_\EK = 0$
$\Leftrightarrow$ $\inf_{x_0\in K}||x+x_0||_E = 0$
$\Leftrightarrow$ $||x||_E \le 0$ since $0\in K$
$\Leftrightarrow$ $x=0$
$\Leftrightarrow$ $\pi(x)=0$ (since $\pi$ is linear)
$\Leftrightarrow$ $\dx = K = \dot0$.

$||\lambda \dx||_\EK = \inf_{x_0\in K}||\lambda x+x_0||_E
= \inf_{x_0\in K}||\lambda x+\lambda x_0||_E
= \inf_{x_0\in K}|\lambda|\,|| x+ x_0||_E
|\lambda|\,|| \dx||_\EK
$.

$||\dx + \dy||_\EK  = \inf_{x_0,y_0\in K}||x + y +x_0+y_0||_E
\le \inf_{x_0,y_0\in K}||x + x_0||_E +||y+y_0||_E
\le ||\dx||_\EK + ||\dx + \dy||_\EK $.

Thus $|| . ||_\EK$ is a norm in~$\EK$.

Let $(\dx_n)_\NNs$ be a Cauchy sequence in $\EK$,
that is, $||\pi(x_n - x_m)||_\EK = ||\dx_n - \dx_m||_\EK \mrar_{n,m\rar\infty}0$.
Let a subsequence, still denoted $(dx_n)$ \st\ $||\pi(x_{k+1} - x_k)||_\EK < {1\over 2^k}$ for all $k\in\NNs$.
Thus $\exists (y_k)_\NNs\in K$ \st\ $||x_{k+1} - x_k - y_k||_E < {2\over 2^k}$ for all $k\in\NNs$, see~\eref{eqndx}.
Then let $(z_k)_\NNs$ be defined by $z_1=0$ ad $y_k = z_{k+1} - z_k$.
Thus $||x_{k+1}-z_{k+1} - (x_k - y_k)||_E < {2\over 2^k}$ for all $k\in\NNs$,
thus $||x_{n+1}-z_{n+1} - (x_m - y_m)||_E \mrar_{n,m\rar\infty}0$,
thus $((x_n-z_n)_\NNs$ is a Cauchy sequence in~$E$,
thus converges to a limit $w\in E$. Thus $\pi(x_n-z_n) = \pi(x_n)-0$ converges to $\pi(w)\in \EK$,
and $\EK$ is closed.

$||\pi(x)||_\EK = \min_{x_0\in K}||x+x_0||_E$, and $0\in K$ (linear subspace),
thus $||\pi(x)||_\EK \le ||x||_E$.
\findem

Let 
\be
\tT :
\left\{\eqalign{
\EK &\rar F \cr
\dx & \rar \tT(\dx) := T(x) \qwhen x\in\dx,
}\right.
\ee
definition independent of $x \in \dx$ since $T(x+x_0) = T(x)$ for all $x_0 \in K$ ($= \Ker(T)$).
In other words, $\tT$ is characterized by $\tT\circ \pi = T$.

\deblem
$\tT$ is linear, injective and continuous with $||\tT|| = ||T||$.
\finlem

\debdem
With $\dx = x+K$ and $\dy = y+K$ we get $\dy + \lambda \dx = x+\lambda y +K$ since $K$ is a linear space,
thus $\tT(\dx+\lambda \dy) = T(x+\lambda y) = T(x)+\lambda T(y) = \tT(\dx)+\lambda \tT(\dy)$,
and $\tT$ is linear.

$\tT.\dx = 0$ $\Rightarrow$ $T(x+x_0)=0$ for all $x_0\in K$ $\Rightarrow$ $x+x_0\in K$ $\Rightarrow$ $x\in K$
$\Rightarrow$ $\dx=\dot0$, thus $\tT$ is injective.

Let $\dx\in \EK$. We have $||\tT(\dx)||_F = ||T.(x+x_0)||_F \le ||T||\,||x+x_0||_E$ for all $x_0\in K$,
thus $||\tT(\dx)||_F \le ||T||\,\min_{x_0\in K} ||x+x_0||_E = ||T||\,||\dx||_\EK$.
Thus $\tT$ is continuous, with $||\tT||\le ||T||$.

Let $x\in E$. We have $||T.x||_F = ||\tT.\dx||_F \le ||\tT||\, ||\dx||_\EK$, thus 
$||T.x||_F \le ||\tT||\,||x+x_0||_E$  for all $x_0\in K$, with $T.x=T(x+x_0)$ for all $x_0\in K$,
thus $||T(x+x_0)||_F \le ||\tT||\,||x+x_0||_E$ for all $x_0\in K$,
$||T.x||_F \le ||\tT||\,||x||_E$. Thus $||T|| \le ||\tT||$.
\findem

\debcor
Let $E$ and $F$ be banach spaces.
If $T\in \calL(E;F)$ (linear and continuous), and if $\RT$ is closed in~$F$,
then %$(\RT,||.||_F)$ is a Banach space (and $T:E \rar \RT$ is surjective by definition of~$\RT$), and
\be
\label{eq2ai2z}
\exists \gamma>0 ,\;
\forall x\in E,\quad ||T.x||_{F} \ge \gamma ||x||_{E/\Ker(T)}
\quad (= \gamma \inf_{x_0 \in \Ker(T)} ||x+x_0||_E).
\ee
%(The closed range theorem is not this corollary, but is theorem~\ref{thmcr}.)
\fincor

\debdem
$K:= \Ker(T) = T^{-1}(\{0\})$ is closed since $T$ is continuous.

$\RT$ is closed in~$F$, therefore $(\RT,||.||_F)$ is a Banach space denoted~$\RT$.
Then $\tT_R : \dx\in \EK \rar \tT_R(\dx) = T(x) \in \RT$ is linear continuous and bijective
between Banach spaces.
Thus \eref{eq2ai}  gives
$\exists \gamma>0$, $\forall \dx\in \EK$, $||\tT_R.\dx||_{F} \ge \gamma ||\dx||_\EK$, \ie~\eref{eq2ai2z}.
\findem

%%%%%%%%%%%%%%%%%%%%%%%%%%%%%%%%%%%%%%%%%%%%%%%%%%%%%%%%%%%%%%%%%%%%%%%%%%%%%%%%%%%

\subsection{The inf-sup condition}

\eref{eq2ai2z} is rewritten 
\be
\label{eq2aib}
\exists \gamma>0 ,\; 
\inf_{x\in E} {||T.x||_F \over ||x||_{E/\Ker(T)} } \ge \gamma .
\ee
(Light writing of $\exists \gamma>0 ,\; 
\inf_{x\in E-\{0\}} {||T.x||_F \over ||x||_{E/\Ker(T)} } \ge \gamma$.)

Consider $B\in \calL(E;F')$ (linear and continuous).
Then~\eref{eq2aib} gives 
\be
\label{eq2aic}
\exists \gamma>0 ,\; 
\inf_{x\in E} {||B.x||_{F'} \over ||x||_{E/\Ker(T)} } \ge \gamma,
\ee
Let $b\dd:E\times F \rar \RR$ be the bilinear form defined by, for all $(x,y)\in E\times F$,
\be
b(x,y) = \la B.x,y \ra_{F',F}.
\ee
Since $||B.x||_{F'} = \sup_{y\in F} { |\la B.x,y \ra_{F',F}| \over ||y||_F}$,
\eref{eq2aic} gives
\be
\label{eq2aid}
\exists \gamma>0 ,\; 
\inf_{x\in E} (\sup_{y\in F} {b(x,y) \over ||x||_{E/\Ker(T)} ||y||_F }) \ge \gamma,
\ee
named the inf-sup condition satisfied by $b\dd$

%%%%%%%%%%%%%%%%%%%%%%%%%%%%%%%%%%%%%%%%%%%%%%%%%%%%%%%%%%%%%%%%%%%%%%%%%%%%%%%%%%%

\section{Some spaces and their duals}
\label{secspaces}

%%%%%%%%%%%%%%%%%%%%%%%%%%%%%%%%%%%%%%%%%%%%%%%%%%%%%%%%%%%%%%%%%%%%%%%%%%%%%%%%%%%

\subsection{Divergence, Gradient, Rotationnal}

Let $(\ve_i)$ be a Euclidean basis in~$\RRn$,
let $\dd_\RRn$ be the associated inner product, and let $\vv.\vw := (\vv,\vw)_\RRn$.
Let $\Omega$ be an open set in $\RRn$.
Let $(\ve_i)$ be a basis in~$\RRn$ and ${\pa f \over \pa x^i} := df.\ve_i$.

The divergence operator is formally given by
\be
\label{eqdefdvg}
\dvg : 
\left\{\eqalign{
\calF(\Omega;\RRn) &\rar \RR \cr
\vv = \sumin v^i\ve_i & \rar \dvg\vv = \sumin {\pa v^i \over \pa x^i}.
}\right.
\ee
(The real value $\dvg\vv $ does not depend on the choice of the basis).

The gradient operator is formally given by
\be
\label{eqdefvgrad}
\vgrad : 
\left\{\eqalign{
\calF(\Omega;\RR) &\rar \RRn \cr
f & \rar \vgrad f = \sumin {\pa f\over \pa x^i}\ve_i.
}\right.
\ee

The rotationnal operator is formally given by
\be
\vcurl : 
\left\{\eqalign{
\calF(\Omega;\RR^3) &\rar \RR^3 \cr
\vv = \sumin v^i\ve_i & \rar \vcurl \vv
= ({\pa v_3 \over \pa x_2} - {\pa v_2 \over \pa x_3}) \ve_1
+({\pa v_1 \over \pa x_3} - {\pa v_3 \over \pa x_1}) \ve_2
+({\pa v_2 \over \pa x_1} - {\pa v_1 \over \pa x_2}) \ve_3.
}\right.
\ee
(In~$\RR^2$, $\curl : \vv\in\RR^2 \rar \curl\vv = {\pa v_2 \over \pa x_1} - {\pa v_1 \over \pa x_2} \in \RR$.)

%%%%%%%%%%%%%%%%%%%%%%%%%%%%%%%%%%%%%%%%%%%%%%%%%%%%%%%%%%%%%%%%%%%%%%%%%%%%%%%%%%%

\subsection{Some Hilbert spaces}
\label{secH0}

Let $\Omega$ be an open set in~$\RRn$, $n=1,2,3$, and $\Gamma =\pa\Omega$ be its boundary.
%Hilbert initiaux :
$$
\eqalignrllq{
&\Ldo = \{ f : \Omega\rar\RR : \int_\Omega f^2\,d\Omega <\infty\},
&
(f,g)_\Ld = \int_\Omega fg\,d\Omega.
\cr
&\Huo = \{f\in \Ldo : \vgrad f \in \Ldon\},
& (f,g)_\Hu = (f,g)_\Ld+(\vgrad f,\vgrad g)_\Ld.
\cr
&\Hdo = \{f\in \Huo : d^2f \in \Ldo^{n^2}\},
& (f,g)_\Hd = (f,g)_\Hu + (d^2f , d^2g)_\Ld.
\cr
&\Hdvgo = \{\vv\in \Ldon : \dvg \vv \in \Ldo\},
& (\vu,\vv)_\Hdvg = (\vu,\vv)_\Ld+(\dvg \vu,\dvg \vv)_\Ld.
\cr
&\Hcurlo = \{\vv\in \Ldo^3 : \vcurl \vv \in \Ldo^3\},
& (\vu,\vv)_\Hcurl = (\vu,\vv)_\Ld+(\vcurl \vu,\vcurl \vv)_\Ld.
\cr
}
$$
\comment{
\Eg\ if $(\ve_i)_{i=1,...,n}$ is a Euclidean basis, denoted $(\ve_i)$,
and if $\vv = \sumin v^i\ve_i$ (generic), then
$(f,g)\Ld = $
\be
(d^2 f,d^2 g)_\Ld = \int_\Omega \sum_{ij=1}^n {\pa ^2 f \over \pa x^i\pa x^j} {\pa ^2 g \over \pa x^i\pa x^j} \,d\Omega,
\qand
(\vu,\vv)_\Ld = \int_\Omega \sumin u^iv^i\,d\Omega.
\ee
And the norm is given by $||\vv||_E = \sqrt{(\vv,\vv)_E}$ (generic).
}

%%%%%%%%%%%%%%%%%%%%%%%%%%%%%%%%%%%%%%%%%%%%%%%%%%%%%%%%%%%%%%%%%%%%%%%%%%%%%%%%%%%

Integrations by parts	:
If $f\in\Huo$ and $\vv\in\Hdvgo$ then
\be
\int_\Omega \vgrad f . \vv \,d\Omega = -\int_\Omega f\, \dvg\vv \,d\Omega
+ \int_\Gamma f\,\vv.\vn\,d\Gamma.
\ee
If $\vv\in\Huo^n$ and $\vw\in\Hcurlo$ then
% avec $\la \vw\wedge\vn , \vv \ra_{\Hmdemig,\Hdemig} \eqnote \int_\Gamma (\vw\wedge\vn).\vv\,d\Gamma$ :
\be
\int_\Omega \vcurl\vv . \vw \,d\Omega
= +\int_\Omega \vv.\vcurl\vw \,d\Omega + \int_\Gamma \vv.(\vw\wedge\vn)\,d\Gamma.
\ee

%%%%%%%%%%%%%%%%%%%%%%%%%%%%%%%%%%%%%%%%%%%%%%%%%%%%%%%%%%%%%%%%%%%%%%%%%%%%%%%%%%%

\subsection{Some sup-spaces}
\label{secH}

Closures $\DO= C^\infty_c(\Omega)$ (space of $C^\infty$ functions with compact support):
$$
\eqalignrllq{
&\Huzo=\{f\in \Huo : f_{|\Gamma}=0\} = \overline{\DO}^\Hu, % = \Ker(\gamma_0),
& (f,g)_\Huz = (\vgrad f,\vgrad g)_\Ld.
\cr
&\Hdzo=\{f\in \Hdo : f_{|\Gamma}=0 \hbox{ et } \vgrad f.\vn_{|\Gamma} = 0 \}  = \overline{\DO}^\Hd,
& (f,g)_\Hdz = (d^2f , d^2g)_\Ld.
\cr
&\Hdvgzo = \{\vv\in \Hdvgo : (\vv.\vn)_{|\Gamma}=0\} = \overline{\DO^n}^\Hdvg,
& (\vu,\vv)_\Hdvgz = (\dvg \vu,\dvg \vv)_\Ld.
\cr
&\Hcurlzo = \{\vv\in \Hcurlo : (\vv\wedge\vn)_{|\Gamma}=0\} = \overline{\DO^3}^\Hcurl,
& (\vu,\vv)_\Hcurlz = (\vcurl \vu,\vcurl \vv)_\Ld.
\cr
}
$$
When $\Omega$ is bounded, the given semi-inner products are equivalent to the inner products of the embedding spaces.

%%%%%%%%%%%%%%%%%%%%%%%%%%%%%%%%%%%%%%%%%%%%%%%%%%%%%%%%%%%%%%%%%%%%%%%%%%%%%%%%%%%

\subsection{Trace operator $\gamma_0$ and the Hilbert space $\Hdemig$}

The trace mapping is
\be
\gamma_0 : 
\left\{\eqalign{
\Huo & \rar \Ldg, \cr
f & \rar \gamma_0(f) = f_{|\Gamma}.
}\right. 
\ee
%It is linear and continuous, see~\cite{brezis}.
(Same notation ifor $\gamma_0 : \vv \in \Huon  \rar \gamma_0(\vv) = \vv_{|\Gamma} \in \Ldg^n$.)
If $\Omega$ is regular, then
%On~a, supposant $\Omega$ régulier (%orienté et $\Gamma$ lipschitzienne $C^{1,1}$) :
\be
\Huzo = \Ker(\gamma_0).
\ee
With $\Im(T)$ the range of a mapping~$T$, let
\be
\label{eqngammaz}
\Im(\gamma_0) = \Hdemig, \qand
||d||_\Hdemig = \inf_{u\in\Huo : u_{|\Gamma} = d}||u||_\Hu.
\ee

\debprop
\label{propcgz}
Let $d\in\Hdemig$ and let $u_d \in \Huo$ be the solution of the Dirichlet problem: Find $u\in\Huo$ \st\
\be
\label{eqhdemigon}
\left\{\eqalign{
&{-}\Delta u + u = 0\quad\hbox{in }\Omega, \cr
& u_{|\Gamma}=d\quad\hbox{on }\Gamma ,\cr
}\right.
%\qe (u_d,v_0)_\Hu = 0,\quad \forall v_0 \in\Huzo.
\ee
then
\be
\label{eqhdemign}
||d||_\Hdemig = ||u_d||_\Hu,
\ee
and $||.||_\Hdemig$ is the norm of the inner product
\be
\label{eqcarhd}
%(\inf_{w\in\Huo : w_{|\Gamma} = d}||w||_\Hu=) \quad 
%||d||_\Hdemig = ||u_d||_\Huo, \qe
(c,d)_\Hdemig := (u_c,u_d)_\Huo.
\ee
Moreover $(\Hdemig,\dd_\Hdemig)$ is a Hilbert space.

Moreover $\gamma_0 : v\in (\Huo,||.||_\Huo) \rar v_{|\gamma} \in (\Hdemig,||.||_\Hdemig)$ is (linear) continuous.
\finprop

\debdem
Let $z_d\in\Huo$ be an counter image of $d\in\Hdemig$ (exists by definition of~$\Hdemig$, \cf~\eref{eqngammaz}).
So $\gamma_0(z_d) = d =z_{d|\Gamma}$.
Let $u=u_0+z_d \in \Huzo +z_d %\eqnote d+\Huzo
$. % (espace affine passant par~$z_d$, d'espace vectoriel associé~$\Huzo$).
With~\eref{eqhdemigon} we get
\be
\left\{\eqalign{
&{-}\Delta u_0 + u_0 = \Delta z_d - z_d \quad\hbox{dans }\Hmuo, \cr
& u_{0|\Gamma}=0\quad\hbox{dans }\Hdemig .\cr
}\right\}
\ee
Thus $u_0 \in \Huzo$ satisfies
\be
\label{eqhdemigonv}
(u_0,v_0)_\Huo = -(z_d,v_0)_\Huo ,\qquad\forall v_0\in\Huzo,
\ee
and the Lax--Milgram theorem gives the existence of a solution $u_0\in\Huzo$.
Then we check that $u_d=u_0+z_d$ is independent of the chosen~$z_d$:
If $z_d'$ satisfies $\gamma(z_d')=d$, if the associate solution is~$u_0'$,
if $u_d'=u_0'+z_d'$, then  $u_d-u_d'\in\Huzo$ and $(u_d-u_d',v_0)_\Huo = 0$ for any $v_0\in\Huzo$,
so $u_d-u_d'=0$ (Lax--Milgram theorem).
Moreover \eref{eqhdemigonv} tells that $u_0+z_d = u_d \perp_\Hu \Huzo$.
Thus we get, for any $v_0\in \Huzo$,
$$
||u_d+v_0||_\Huo^2 = ||u_d||_\Huo^2 + ||v_0||_\Huo^2 + 2(u_d,v_0)_\Huo
=||u_d||_\Huo^2 + ||v_0||_\Huo^2 + 0 .
$$
So $\inf_{w\in\Huo\atop w_{\Gamma}=d}||w||_\Huo^2 = \inf_{v_0\in\Huzo}||u_d+v_0||_\Huo^2 = ||u_d||_\Huo^2$,
denoted $||d||_\Hdemig$.
Then we define $(c,d)_\Hdemig$ as in~\eref{eqcarhd},
and $\dd_\Hdemi$ is trivially bilinear symmetric positive (is an inner product).

Then we check that $(\Hdemig,||.||_\Hdemi)$ is complete:
If $(d_n)_\NNs \in \Hdemig$ is a Cauchy sequel in~$\Hdemig$
and if $u_{d_n}$ is the solution of~\eref{eqhdemigon},
then $(u_{d_n})_\NNs$ is a Cauchy sequel in $(\Huo,||.||_\Hu)$, \cf~\eref{eqhdemign},
so converges in~$\Huo$ (since $\Huo$ is complete) toward some $u\in\Huo$.
Then let $d:=\gamma_0(u) \in\Hdemig$.
Since $(u_{d_n}, v_0)_\Huo = 0$ for any $v_0\in\Huzo$, \cf~\eref{eqhdemigon},
we get $(u_d, v_0)_\Huo = 0$ for any $v_0\in\Huzo$  (continuity of an inner product relatively to itself).
Thus $u_d$ is the solution of~\eref{eqhdemigon},
and $||d-d_n||_\Hdemi = ||u-u_{d_n}||_\Hu \matrarrow_{n\rar\infty}0$.

And $\gamma_0 : (\Huo,||.||_\Hu) \rar (\Hdemig,||.||_\Hdemig$ (linear) satisfies,
for $u\in\Huo$, with $d:=\gamma_0(u)$ and $u_d$ solution of~\eref{eqhdemigon},
$||\gamma_0(u)||_\Hdemig = ||u_d||_\Huo \le ||u_d+u_0||_\Huo$ for any $u_0\in\Huzo$,
thus with $u_0=u-u_d$ we get $||\gamma_0(u)||_\Hdemig = ||u||_\Huo$,
so $\gamma_0$ is bounded ($||\gamma_0||\le 1$).
\findem

%%%%%%%%%%%%%%%%%%%%%%%%%%%%%%%%%%%%%%%%%%%%%%%%%%%%%%%%%%%%%%%%%%%%%%%%%%%%%%%%%%%

\subsection{Some other trace operators}

\be
\gamma_1 :
\left\{\eqalign{
\Hdo & \rar \Ldg, \cr
f & \rar \gamma_1(f) = (\vgrad f)_{|\Gamma}.\vn \eqnote {\pa f \over \pa\vn}_{|\Gamma}.
}\right.
\ee
\be
\gamma_n :
\left\{\eqalign{
\Hdvgo & \rar \Hmdemig, \cr
\vv & \rar \gamma_n(\vv) = \gamma_0(\vv).\vn \eqnote (\vv.\vn)_{|\Gamma}
}\right.
\ee
(and the divergence operator enables the control of $\vv.\vn$ on~$\Gamma$),
\be
\vgamma_t :
\left\{\eqalign{
\Hcurlo & \rar \Hmdemig^3, \cr
\vv & \rar \vgamma_t(\vv) = \gamma_0(\vv) \wedge \vn \eqnote (\vv\wedge\vn)_{|\Gamma}
}\right.
\ee
(and the rotational operator enables the control $\vv\wedge\vn$ on~$\Gamma$.

%%%%%%%%%%%%%%%%%%%%%%%%%%%%%%%%%%%%%%%%%%%%%%%%%%%%%%%%%%%%%%%%%%%%%%%%%%%%%%%%%%%

\subsection{Some Dual spaces}

Banach spaces:
\be
\label{eqdetb}
\eqalignrllq{
& (\Ldo)' \simeq \Ldo\quad\hbox{(usual identification)}
& ||f||_{\Ldo'} = ||f||_\Ldo .
\cr
&\Ldzo = \{ f \in\Ldo  : \int_\Omega f\,d\Omega = 0\} \simeq \Ldo/\RR,
& ||f||_\Ldz = \inf_{c\in\RR} ||f+c||_\Ld.
\cr
&\Hmuo = \Huzo',
& ||f||_\Hmu = \sup_{v\in\Huzo} {\la f , v \ra \over ||v||_\Huz}.
\cr
& \Hmdemig = (\Hdemig)',
& ||\mu||_\Hmdemig = \sup_{\lambda\in\Hdemig} {|\la\mu,\lambda\ra| \over ||\lambda||_\Hdemig}.
\cr
}
\ee
We have identified $(\Ldo)'$ with~$\Ldo$ thanks to the Riesz representation theorem
in $(\Ldo,\dd_\Ld)$, that is,
\be
\label{eqRieszLd}
\forall \ell \in \Ldo',\; \exists ! f\in\Ldo,\; \forall g \in \Ldo,\;\la \ell, g\ra_{\Ld',\Ld} = (f,g)_\Ldo, \qand
||f||_\Ldo = ||\ell||_{\Ldo'}.
\ee
Thus $\Huzo \subset \Huo \subset \Ldo \simeq \Ldo' \subset (\Huo)' \subset \Hmuo$.
And $\Ldo$ is named the ``pivot space'' (a~central space in distribution theory of Schwartz~\cite{schwartz}).

% $\Ldo/\RR$ est l'ensemble des fonctions~$\Ldo$ définies à une constante près, identifié à~$\Ldzo$. Et la norme $||.||_\Ldz$ donnée est la norme quotient usuelle.

\debprop
Let $\lambda\in \Hmdemig$ and let $w_\lambda\in\Huo$ be the solution of the Nenmann problem: Find $w\in\Huo$ \st\
\be
\label{eqhmde}
\left\{\eqalign{
 &  -\Delta w + w =0\quad\hbox{dans }\Omega,\cr
 & {\pa w\over\pa n} = \lambda\quad\hbox{sur }\Gamma,\cr
}\right.
%\} \qs (w,v)_\Huo = \la\lambda,v\ra_{\Hmdemig,\Hdemig},\quad\forall v\in\Huo.
\ee
then
\be
\label{eqhmde0}
||\lambda||_\Hmdemig = ||w_\lambda||_\Huo .
%\quad  (=\sup_{d\in\Hdemig} |\la \lambda , {d \over ||d||_\Hdemig} \ra_{\Hmdemi,\Hdemi}| ),
\ee
\finprop

\debdem
\eref{eqhmde} reads $(w,v)_\Huo = \la\lambda,v\ra_{\Hmdemig,\Hdemig}$ for all $v\in\Huo$,
thus~\eref{eqhmde} has a unique solution $w_\lambda$ 
(Lax--Milgram theorem:
The bilinear form given by $a(u,v)=(u,v)_\Huo$ is  trivially $\Huo$-continuous and coercive,
and the linear form given by $\ell(v)= \la\lambda,v\ra_{\Hmdemig,\Hdemig}$
is continuous since $|\ell(v)| \le ||\lambda||_\Hmdemig||\gamma_0(v)||_\Hdemig
\le ||\lambda||_\Hmdemig||v||_\Huo$, \cf\ prop.~\ref{propcgz}). And :
\be 
\eqalign{
||\lambda||_\Hmdemig 
= &\sup_{d\in \Hdemig} { |\la \lambda , d \ra_{\Hmdemig,\Hdemig}|  \over ||d||_\Hdemig}
 \quad (\hbox{definition})\cr
= & \sup_{d\in \Hdemig} { |\la \lambda , u_d \ra_{\Hmdemig,\Hdemig}|  \over ||u_d||_\Huo}
 \quad (\hbox{\cf\ \eref{eqcarhd}})\cr
= &\sup_{d\in \Hdemig} { |(w_\lambda , u_d)_\Huo| \over ||u_d||_\Huo}
 \quad (\hbox{\cf\ \eref{eqhmde}})\cr
% = &\sup_{u\in \Huo} |(w_\lambda , {u\over ||u||_\Huo})_\Huo|  \quad (\hbox{surjectivité de $\gamma_0$})\cr
\le & ||w_\lambda||_\Huo,
 \quad (\hbox{Cauchy--Schwarz in~$\Huo$}).\cr
}
\ee
Et $\gamma_0(w_\lambda) \in \Hdemig$ gives
\be
\eqalign{
||\lambda||_\Hmdemig  
\ge & {|\la \lambda , \gamma_0(w_\lambda) \ra_{\Hmdemig,\Hdemig} \over ||\gamma_0(w_\lambda)||_\Hdemig}|
\quad \hbox{(by definition of the sup)} \cr
\ge & {|(w_\lambda , w_\lambda \ra_\Huo|\over ||\gamma_0(w_\lambda)||_\Hdemig}
\quad \hbox{(\cf~\eref{eqhmde})} \cr
\ge & %{ ||w_\lambda||_\Huo^2 \over  ||w_\lambda||_\Huo}= 
||w_\lambda||_\Huo
\quad \hbox{(since $||\gamma_0(w_\lambda)||_\Hdemig \le ||w_\lambda||_\Huo$ \cf~\eref{eqngammaz})}. \cr
}
\ee
Thus~\eref{eqhmde0}.
\findem

%%%%%%%%%%%%%%%%%%%%%%%%%%%%%%%%%%%%%%%%%%%%%%%%%%%%%%%%%%%%%%%%%%%%%%%%%%%%%%%%%%%

\subsection{Dual of $\Huo$ and $\Hdvgo$ (characterizations)}

\debthm
 {\bf Dual of~$\Huo$.}
%Si $F\in (\Huo)'$ alors il existe $(f,\vu)\in\Ldo \times \Ldon$ t.q. :
\be
\label{eqdhu}
\ell \in (\Huo)'  \quad\Rightarrow\quad
\exists (f,\vu)\in\Ldo \times \Ldon :
\la \ell,\psi\ra
= (f,\psi)_\Ld + (\vu,\vgrad\psi)_\Ld
\quad\forall \psi\in\Huo.
\ee
{\bf Dual of~$\Huzo$.}
\be
\label{eqdhuzo}
\ell\in \Hmuo  \quad\Rightarrow\quad
\exists (f,\vu)\in\Ldo \times \Ldon \qtq \ell = f - \dvg\vu.
\ee
And if $\Omega$ is bounded then we can choose $f=0$ (with $(\Huzo,\dd_\Huz)$).
\finthm

\debdem (from Brézis~\cite{brezis}.)
Characterization of~$\Huon$.
Define $Z := \Ldo \times \Ldon$ provided with the inner product
$((\phi,\vu),(\psi,\vv))_Z = (\phi,\psi)_\Ld + (\vu,\vv)_\Ld$
so that $Z$ is a Hilbert space.
Define $T :
\left\{\eqalign{
\Huo & \rar Z \cr
\psi & \rar T\psi = (\psi,\vgrad \psi)
}\right\}
$.
So $||\psi||_\Hu = ||T\psi||_Z = ||(\psi,\vgrad \psi)||_Z$,
and $T : (\Huzo,\dd_\Huz) \rar (\Im(T),||.||_Z)$ is an isometry.
Let $\ell\in \Huo'$. Define
$$
\Phi_{\Im(T)} :
\left\{\eqalign{
\Im(T) & \rar \RR \cr
(\psi,\vv{=}\vgrad\psi) & \rar
\la \Phi_{\Im(T)},(\psi,\vv)\ra_{Z',Z}
= \la \ell,T^{-1}(\psi,\vv)\ra_{\Hu',\Hu}
= \la \ell,\psi\ra_{\Hu',\Hu}.
}\right.
$$
$\Phi_{\Im(T)}$ is linear (trivial) and continuous since $\ell$ and $T^{-1}$ are.
With Hahn--Banach theorem, extend $\Phi_{\Im(T)}$ to~$Z$, so that we get a linear countinous form
$
\Phi_Z :
\left\{\eqalign{
Z & \rar \RR \cr
(\psi,\vv) & \rar \la \Phi_Z,(\psi,\vv)\ra
}\right\}
$.
%(avec $||\Phi_Z||_{Z'} = ||\Phi_{\Im(T)}||_{Z'}$).
Then the Riesz representation theorem gives:
$\exists(\phi,\vu)\in Z$ \st\ $\la \Phi_Z,(\psi,\vv)\ra = ((\phi,\vu),(\psi,\vv))_Z
= \int_\Omega \phi\,\psi\,d\Omega + \int_\Omega \vu.\vv\,d\Omega$ for all $(\psi,\vv)\in Z$.
Then take $(\psi,\vv{=}\vgrad\psi)\in \Im(T)$ to get~\eref{eqdhu}.

Similar proof for~\eref{eqdhuzo}.
\findem

\debthm
 {\bf Dual de~$\Hdvgo$.}
\be
\label{eqhdvgd}
F\in (\Hdvgo)'  \quad\Rightarrow\quad
\exists(\vf,\phi)\in\Ldon \times \Ldo \hbox{ \st\ }
\la F,\vv\ra
= (\vf,\vv)_\Ld + (\phi,\dvg\vv)_\Ld
,\quad\forall \vv\in\Hdvgo.
%\quad (\int_\Omega \vf.\vv\,d\Omega + \int_\Omega \phi.\dvg\vv\,d\Omega).
%\qo f\in\Ldo\qe \vg\in\Ldon.
\ee
{\bf Dual de~$\Hdvgzo$.}
In particular
\be
\label{eqhdvgzd}
F\in \Hdvgzo'  \quad\Rightarrow\quad
\exists (\vf,\phi)\in\Ldon \times \Ldo  \hbox{ \st\ } F = \vf - \vgrad\phi.
\ee
And if $\Omega$ is bounded we can choose $\vf=0$ (with $(\Hdvgzo,\dd_\Hdvgz)$).
\finthm

\debdem (Similar to the proof of~\eref{eqdhu}.)
Define $Z = \Ldon \times \Ldo$ provided with the inner product
$((\vu,p),(\vv,q))_Z =(\vu,\vv)_\Ld + (p,q)_\Ld$
so that $(Z,\dd_Z)$ is a Hilbert space.
Define $T :
\left\{\eqalign{
\Hdvgo & \rar Z \cr
\vv & \rar T\vv = (\vv,\dvg \vv)
}\right\}
$.
So $||\vv||_\Hdvg = ||T\vv||_Z = ||(\vv,\dvg \vv)||_Z$,
and $T : (\Hdvgo,||.||_\Hdvg) \rar (\Im(T),||.||_Z)$ is an isometry.
Let $F\in \Hdvgo'$.
The mapping $(\vv,q{=}\dvg\vv)\in \Im(T) \rar \la F,T^{-1}(\vv,q)\ra_{\Hdvg',\Hdvg} = \la F,\vv\ra_{\Hdvg',\Hdvg}$
is a linear form (trivial) that is continuous since $F$ and $T^{-1}$ are.
With Hahn--Banach theorem, extend it to~$Z$ to get a linear continuous form named
$\Phi : (\vv,q)\in Z \rar \la \Phi,(\vv,q)\ra$.
Then the Riesz representation theorem gives:
$\exists(\vu,p)\in Z$ \st\ $\la \Phi,(\vv,q)\ra = ((\vu,p),(\vv,q))_Z
= \int_\Omega \vu.\vv\,d\Omega + \int_\Omega pq\,d\Omega$
for all $(\vv,q)\in Z$.
An choose $(\vv,q{=}\dvg\vv)\in \Im(T)$ to get~\eref{eqhdvgd}.

Similar proof for~\eref{eqhdvgzd}.
\findem

%%%%%%%%%%%%%%%%%%%%%%%%%%%%%%%%%%%%%%%%%%%%%%%%%%%%%%%%%%%%%%%%%%%%%%%%%%%%%%%%%%%

\subsection{Kernel of the trace operators}

$\Omega$ is supposed to be a regular open set.

\be
\eqalignrllq{
&\Ker(\gamma_0) = \Huzo,
& \Im(\gamma_0) = \Hdemig\hbox{ dense in }\Ldg. \cr
& \Ker(\gamma_1) \cap \Ker(\gamma_0) = \Hdzo,
& \Im(\gamma_1) = \Hdemig. \cr
& \Ker(\gamma_n) = \Hdvgzo,
& \Im(\gamma_n) = \Hmdemig. \cr
& \Ker(\vgamma_t) = \Hcurlzo,
& \Im(\vgamma_t) = \Hmdemig^3. \cr
}
\ee

%%%%%%%%%%%%%%%%%%%%%%%%%%%%%%%%%%%%%%%%%%%%%%%%%%%%%%%%%%%%%%%%%%%%%%%%%%%%%%%%%%%

\subsection{Poincaré--Friedrichs}

(See \eg\ manuscript ``Eléments finis'', or Raviart--Thomas~\cite{raviart-thomas},
or Ciarlet~\cite{ciarlet}...

If $\Omega$ is bounded (at least in one direction), then we have Poincaré's inequality in~$\Huzo$:
There exists $c_\Omega>0$ s.t
\be
\label{eqpoinc}
\forall v\in\Huzo,\quad ||v||_\Ld \le c_\Omega ||\vgrad v||_\Ld,
\ee
and the norms $||v||_\Huo$ and $||\vgrad v||_\Ldo$ are equivalent in~$\Huzo$
(this space is closed in~$\Huo$ it is the closure of~$\DO$ in~$\Huo$).

And if $\Omega$ is bounded then there exists $c_\Omega>0$ s.t:
\be
\label{eqpoinc2}
\forall v\in\Huzo\bigcap\Hdo,\qquad ||v||_\Hdo \le c_\Omega||\Delta v||_\Ldo.
\ee
and the norms $||v||_\Hdo$ and $||\Delta v||_\Ldo$ are equivalent in~$\Huzo\bigcap\Hdo$
(this space is not closed in~$\Huo$).

%%%%%%%%%%%%%%%%%%%%%%%%%%%%%%%%%%%%%%%%%%%%%%%%%%%%%%%%%%%%%%%%%%%%%%%%%%%%%%%%%%%

\subsection{$\Ldon$ Decomposition (Helmholtz)}

\def\XO{{X(\Omega)}}
\def\XNO{{X_N(\Omega)}}
\def\XTO{{X_T(\Omega)}}

Let $\Omega\subset\RRn$ an open regular bounded set.
 %With $\dvg:\Hdvgo \rar \Ldo$ we have
Let $\dvg:\Hdvgo \rar \Ldo$, so $\Ker(\dvg) = \{\vv\in\Hdvgo : \dvg\vv=0\}$. And let
\be
\Ker(\dvg)_0 = \Ker(\dvg) \cap \Hdvgzo =  \{\vv\in\Ker(\dvg) : (\vv.\vn)_{|\Gamma}=0\},
\ee
the subspace of incompressible functions with $\Gamma$ impervious.

\debthm
\be
\label{eqHe0}
\left\{\eqalign{
& \Ldon = \vgrad(\Huzo) \oplus^{\perp_\Ld} \Ker(\dvg),\cr
& \Ldon = \vgrad(\Huo) \oplus 	^{\perp_\Ld} \Ker(\dvg)_0,
}\right.
\ee
i.e, for any $\vf\in\Ldon$ there exists $(\phi,\vw)\in \vgrad(\Huzo) \times \Ker(\dvg)$ for~\eref{eqHe0}$_1$,
and there exists $(\phi,\vw)\in \vgrad(\Huo) \times \Ker(\dvg)_0$ for~\eref{eqHe0}$_2$, \st\
\be
\label{eqHe1}
\vf = \vgrad \phi + \vw, \qwith (\vgrad\phi , \vw)_\Ld= 0.
\ee
%and $||\vgrad\phi||_\Ld + ||\vw||_\Ld \le 2 ||\vf||_\Ld$.
\finthm

\debdem
Let $\vf\in\Ldon$.

For~\eref{eqHe0}$_1$, consider the solution of the homogenous Dirichlet problem:
Find $\phi\in\Huzo$ \st\ $\Delta\phi = \dvg\vf$ (distribution),
meaning, find $\phi\in\Huzo$ \st\ 
$(\vgrad\phi,\vgrad\psi)_\Ld = (\vf,\vgrad\psi)_\Ld$
for all $\psi \in\Huzo$.
The Lax--Milgram theorem gives a unique solution $\phi \in \Huzo$.

Let $\vw = \vf - \vgrad\phi \in \Ldon$.
So $(\vw,\vgrad\psi)_\Ld %= (\vf-\vgrad\phi,\vgrad\psi)_\Ld 
= 0$ for all $\psi \in\Huzo$,
by definition of~$\phi$,
thus $\vw \perp_\Ld \vgrad(\Huzo)$ and
$\dvg\vw = 0 \in \Hmuo$, and $0\in\Ldo$, thus $\vw\in\Hdvgo$ and $\vw\in\Ker(\dvg)$;
Thus $\vf = \vw+\vgrad\phi \in \Ker(\dvg) \oplus^{\perp_\Ld} \vgrad(\Huzo)$,
thus~\eref{eqHe0}$_1$.

For~\eref{eqHe0}$_2$, consider the solution of the homogenous Neumann problem:
Find $\phi\in\Huo$ \st\ 
$\int_\Omega \vgrad\phi . \vgrad\psi \,d\Omega = \int_\Omega \vf.\vgrad\psi\,d\Omega$
for all $\psi \in\Huo$.
The Lax--Milgram theorem gives a unique solution $\phi \in \Huo/\RR$ (\ie\ up to a constant),
moreover with $\phi\in\Hdo$ (regularity result thanks to $\vf\in\Ldo$),
so that $-\la \Delta \phi,\psi\ra_{(\Huo)',\Huo} + \int_\Gamma \vgrad\phi(x).\vn(x)\psi(x)\,d\Gamma
= -\la \dvg\vf, \psi\ra_{(\Huo)',\Huo}$ for all $\psi\in\Huo$.
In particular $\psi\in\Huzo$ gives $\Delta \phi = \dvg\vf \in (\Huo)'$,
and we are left with $\int_\Gamma \vgrad\phi(x).\vn(x)\psi(x)\,d\Gamma$ for all $\psi\in\Huo$,
thus for all $\psi_{|\Gamma}\in\Hdemig$,
thus $\vgrad\phi.\vn_{|\Gamma} = 0$.

Let $\vw = \vf - \vgrad\phi \in \Ldon$.
Thus $(\vw,\vgrad \psi)_\Ld = 0$ for all $\psi\in\Huo$, by definition of~$\phi$,
thus $\vw \perp \vgrad(\Huo)$.
And $\dvg\vw = \dvg\vf - \Delta\phi = 0$, thus $\dvg\vw\in\Ldo$ and $\vw\in\Ker(\dvg)$.
With $\int_\Omega \vw . \vgrad\psi \,d\Omega = 0$ for all $\psi\in\Huo$,
thus $\int_\Omega \vw.\vn\,\psi\,d\Gamma=0$ for all $\psi\in\Huo$,
and $\vw.\vn=0 \in \Hmdemig$ (since $\Hdemig$ is dense in~$\Ldg$),
thus $\vw\in \Ker(\dvg)_0$,
thus $\vf = \vgrad \phi + \vw \in \vgrad(\Huo) \oplus 	^{\perp_\Ld} \Ker(\dvg)_0$,
thus~\eref{eqHe0}$_2$. % et~\eref{eqHe1}.
%And the orthogonal decomposition gives the last inequality.
\findem

%Voir par exemple Schweizer~\cite{schweizer}.

%%%%%%%%%%%%%%%%%%%%%%%%%%%%%%%%%%%%%%%%%%%%%%%%%%%%%%%%%%%%%%%%%%%%%%%%%%%%%%%%%%%
%%%%%%%%%%%%%%%%%%%%%%%%%%%%%%%%%%%%%%%%%%%%%%%%%%%%%%%%%%%%%%%%%%%%%%%%%%%%%%%%%%%

\section{A surjectivity of the gradient operator}

See \eg\ Girault--Raviart~\cite{girault-raviart}.
We deal here with infinite dimensional spaces.	
The surjectivity of $\vgrad$ is need for a Stokes like problem, see~\eref{eqgenw1}.

%%%%%%%%%%%%%%%%%%%%%%%%%%%%%%%%%%%%%%%%%%%%%%%%%%%%%%%%%%%%%%%%%%%%%%%%%%%%%%%%%%%

\subsection{The theorem}
\def\calD{{\cal D}}

Let $\Omega$ be an open regular set in~$\RRn$.
Let $(\ve_i)$ be a given Cartesian~$\RRn$,
and $\vn(x)=\sumin n_i(x)\ve_i$ be the outer normal unit to~$\Gamma$ at~$x$.

$\Hmuo = (\Huzo)' = \calL(\Huzo;\RR)$ is the set of continuous linear forms
defined on~$\Huzo$, \cf~\eref{eqdetb}.

With~\eref{eqRieszLd}, $\Ldo'$ is identified to~$\Ldo$,
and $\Hmuo \supset \Ldo' = \Ldo \supset \Huzo$.

And if $g\in\Ldo$, then ${\pa g\over \pa x_i} \in \Hmuo$, and, for all $\phi\in\Huzo$,
\be
\la {\pa g\over \pa x_i},\phi\ra_{\Hmu,\Huz} := -\int_\Omega g(x){\pa\phi \over \pa x_i}(x)\,dx,
\ee
see Schwartz~\cite{schwartz}.
In particular,
if $p\in \Ldon$ then $\vgrad p \in \Hmuon$ and, for all $\vv\in\Huzon$,
\be
\label{eqgradpv}
\la \vgrad p,\vv\ra_{\Hmu,\Huz} = -\int_\Omega p(x)\dvg\vv(x)\,dx = -(p,\dvg\vv)_\Ld.
\ee

\debthm
\label{thmrgo}
The range of the gradient operator
$\vgrad :
\left\{\eqalign{
\Ldo & \rar \Hmuo \cr
p & \rar \vgrad p
}\right\}$
is closed, and its kernel $\Ker(\vgrad)$ is the set of constant functions.
\finthm

\debproof
The proof of this quite difficult theorem is given in the next~\S.
\finproof

And the open mapping theorem, \cf~\eref{eq2ai}, then gives the needed result for the Stokes like problem,
\cf~\eref{eqgenw1}:

\debcor
\be
\label{eqthmrgo}
\exists \beta >0,\; \forall p\in \Ldo,\quad ||\vgrad p||_\Hmu \ge \beta ||p||_{\Ldo/\RR},
\ee
that is $\exists\beta>0$,
$\ds \inf_{p\in\Ldo}\sup_{v\in\Huzo}{ |(\dvg\vv,p)_\Ld | \over ||\vv||_{\Huzo/\Ker(\dvg)} ||p||_\Ldzo} \ge\beta$
(inf-sup condition).
\fincor

\comment{
And the closed range theorem gives:

\debcor
The range of
$\dvg : 
\left\{\eqalign{
\Huzon & \rar \Ldo \cr
\vv & \rar \dvg\vv \cr
}\right\}
$
is closed, thus
\be
\exists \beta >0,\; \forall \vv\in \Huzon,\quad ||\dvg \vv||_\Ldo \ge \beta ||\vv||_{\Huz/\Ker(\dvg)}.
\ee
\fincor
}

%%%%%%%%%%%%%%%%%%%%%%%%%%%%%%%%%%%%%%%%%%%%%%%%%%%%%%%%%%%%%%%%%%%%%%%%%%%%%%%%%%%

\subsection{Steps for the proof}

%%%%%%%%%%%%%%%%%%%%%%%%%%%%%%%%%%%%%%%%%%%%%%%%%%%%%%%%%%%%%%%%%%%%%%%%%%%%%%%%%%%

\subsubsection{Equivalent norms in~$\Hmuo$}

$\Omega$ being bounded, the Poincaré inequality gives:
\be
\label{eqipoinc}
\exists c_\Omega\in\RR,\;\forall q\in\Huzo,\; ||q||_\Ld\le c_\Omega||q||_\Huz.
\ee

\comment{
In $(\Ldo,\dd_\Ld)$, the Riesz representation theorem gives:
\be
\label{eqR}
\forall \ell\in\Ldo',\;\exists \vell\in\Ldo,\;\forall v\in\Ldo,\quad \ell.v = (\vell,v)_\Ldo
\qand ||\ell|| = ||\vell||_\Ldo.
\ee
And we thus identify $\Ldo$ and~$\Ldo'$ (the usual $\dd_\Ldo$ inner product being implicit),
and $\vell$ is denoted~$\ell$.
}

Let $q\in\Ldo$, and let $\ell_q\in\Hmuo$ be defined on~$\Huzo$ by
\be
\label{eqspr100}
\forall \psi\in\Huzo,\quad \la \ell_q, \psi\ra_{\Hmu,\Huz} := (q,\psi)_\Ldo.
\ee
Thus $\ell_q$ is trivially linear, and, with~\eref{eqipoinc},
\be
\label{eqspr10}
\forall \psi\in\Huzo,\quad
|\la \ell_q,\psi \ra_{\Hmu,\Huz}| = |(q,\psi)_\Ldo| \le  ||q||_\Ld||\psi||_\Ld \le c_\Omega ||q||_\Ld||\psi||_\Huz.
\ee
Thus $\ell_q$ is continuous, thus $\ell_q\in\Hmuo$,
$\Ldo$ is considered to be a subspace in~$\Hmuo$.

\debprop
If $q\in\Ldo$, then
\be
\label{eqspr1}
||\ell_q||_\Hmu \le c_\Omega ||q||_\Ld,\quad ||\vgrad q||_\Hmu \le ||q||_\Ld.
\ee
Thus the injection $
\left\{\eqalign{
\Ldo &\rar \Hmuo \cr
q &\rar \ell_q \cr
}\right\}
$
and the gradient operator $
\left\{\eqalign{
\Ldo &\rar \Hmuon \cr
q &\rar \vgrad q \cr
}\right\}
$
are continuous. In particular, with~\eref{eqspr10},
\be
\hbox{if } q\in\Ldo \hbox{ then } \ell_q \eqnote q.
\ee
(The space $\Ldo$ is the pivot space.)
\finprop

\debdem
\def\vphi{{\vec\phi}}
\eref{eqspr1}$_1$ is given by~\eref{eqspr10}.
Let $q\in\Ldo$, with~\eref{eqgradpv} we get $\vgrad q\in\Hmuon$.
We have, for all $\vphi\in\Huzon$, \cf~\eref{eqgradpv},
$$
|\la \vgrad q,\vphi\ra_{\Hmu,\Huz}| = |-(q,\dvg\vphi)_\Ld|
\le ||q||_\Ld||\dvg\vphi||_\Ld
\le ||q||_\Ld||\grad \vphi||_{(\Ld)^{n^2}} \le ||q||_\Ld||\vphi||_{(\Huz)^n}.
$$
Thus \eref{eqspr1}$_2$.
\findem

Let :
\be
\label{eqspr13}
||.||_+:
\left\{\eqalign{
\Ldo & \rar \RR \cr
v & \rar ||v||_+ = ||v||_\Hmu + ||\vgrad v||_\Hmu.
}\right.
\ee
\comment{
This is trivially a norm in~$\Ldo$, and it satisfies, \cf~\eref{eqspr1},
\be
\label{eqspr12}
||q||_\Hmu + ||\vgrad q||_{(\Hmu)^n} \le (1{+}c_\Omega)||q||_\Ldo.
\ee
}

\debcor
In~$\Ldo$ the norms $||.||_\Ld$ and $||.||_+$ are equivalent norms:
\be
\exists c_1,c_2>0,\; \forall v\in\Ldo,\quad c_1||v||_+ \le ||v||_\Ld \le c_2||v||_+.
\ee
\fincor

\debproof
\eref{eqspr13} trivially defines a norm in~$\Ldo$, and~\eref{eqspr1} gives $c_1 = {1\over 1+c_\Omega}$.

Let $Z=(\Ldo,||.||_+)$.
Thanks to~${1\over c_1}$, $Z$ is a Banach space.
Then consider the canonical injection $I_+ : v\in (\Ldo,||.||_\Ldo) \rar I_+(v)=v\in (\Ldo,||.||_Z)$:
it is the algebraic identity and thus is bijective.
And $I_+$ is continuous (thanks to ${1\over c_1}$).
Thus $I_+{}^{-1} : v\in (\Ldo,||.||_Z) \rar I_+(v)=v\in (\Ldo,||.||_\Ldo)$ is continuous
(open mapping theorem~\ref{thmomt0}).
Then let $c_2 = ||I_+{}^{-1}||$.
\finproof

%%%%%%%%%%%%%%%%%%%%%%%%%%%%%%%%%%%%%%%%%%%%%%%%%%%%%%%%%%%%%%%%%%%%%%%%%%%%%%%%%%%

\subsubsection{Rellich theorem  $\Ldo \rar\Hmuo$}

Reminder: An operator $\kappa\in\calL(E;F)$ is compact iff $\overline{\kappa(B_E(0,1))}$ is compact in~$F$.

\deblem
\label{lemkappa}
Let $E$ and $F$ be Banach spaces. If $\kappa\in\calL(E;F)$ is compact,
then it dual $\kappa':F'\rar E'$ is compact.
\finlem

\debdem
Let $(\ell_n)\in B_{F'}(0,1)$.
We have to prove that the sequence $(T'.\ell_n)\in T'(B_{F'})\subset E'$
has a converging subsequence.
Let $K=\overline{T(B_E(0,1))}$. $K$ is a compact in~$F$ since $T$ is compact.
Then consider the restriction $\phi_n=\ell_n{}_{|K}:K\rightarrow\RR$.
So $(\phi_n)_\NNs$ is a sequence in $C^0(K;\RR)$,
and $(\phi_n)_\NNs \subset B_{F'}(0,1)$ is a bounded set in~$F'$.
Moreover $(\phi_n)_\NNs$ is equicontinuous
since $\ell_n$ is linear continuous $||\ell_n(y)||\le||\ell_n||\,||y||_F\le||y||_F$).
Thus the set $(\phi_n)_\NNs$ is relatively compact in $C^0(K;\RR)$ (Ascoli theorem, see Brézis~\cite{brezis}).
Thus we can extract
a convergent subsequence $(\phi_{n_k})_{k\in\NNs}$ in~$C^0(K;\RR)$. % , donc de Cauchy.
Thus, $T(B_E)$ being relatively compact and thus bounded, we have
$$
\sup_{x\in B_E}|\la \ell_{n_k}-\ell_{n_m},T.x\ra|\matrarrow_{k,m\rar\infty}0.
$$
Thus $||T'.\ell_{n_k}-T'.\ell_{n_m}||_{E'}\rar 0$. 
Thus $E'$ being a Banach space, since $E$ is,
$(T'.\ell_{n_k})_{k\in\NNs}$ converges in~$E'$. Thus the set $\overline{(T'.\ell_{n_k})_{k\in\NNs}}$ is compact,
thus $T'$ is compact.
\findem

\debthm (Rellich)
\label{thmReld}
The canonical injection $T:v\in\Ldo\rightarrow v\in\Hmuo$ is compact.
\finthm

\debdem
$I_{10} : v\in\Huzo\rar v\in\Ldo$ is compact, Rellich theorem see Brézis~\cite{brezis},
thus $I_{10}' : v\in \Ldo \rar v\in \Hmuo$ is compact, \cf~Lemma~\ref{lemkappa}.
\findem

%%%%%%%%%%%%%%%%%%%%%%%%%%%%%%%%%%%%%%%%%%%%%%%%%%%%%%%%%%%%%%%%%%%%%%%%%%%%%%%%%%%

\subsubsection{Petree--Tartar compactness theorem}

Let $E$ and $F$ be two Banach spaces, and $T\in\calL(E;F)$ (linear and continuous).
The purpose is to prove that the range of~$T$ is eventually closed.
But to use theorem~\ref{thmomt0} and~\eref{eq0ai} to prove it, can be difficult.
It can be easier to find a compact operator $\kappa:E\rar G$, where $G$ is a Banach space, \st\
\be
\label{eq0aiPT}
\exists \gamma>0 ,\; 
\forall x\in E,\quad ||T.x||_{F} + ||\kappa.x||_G \ge \gamma ||x||_E.
\ee

\debthm
\label{thmPT}
Let $E$, $F$ and $G$ be three Banach spaces,
let $T\in\calL(E;F)$ be injective (one-to-one), and $\kappa\in\calL(E;G)$ be compact.
If~\eref{eq0aiPT} holds then~\eref{eq0ai} holds, and thus $\RT$ is closed.

(If $T$ is not injective, consider $E/\Ker(T)$.)
\finthm

\debproof
\def\nk{{n_k}}
Suppose~\eref{eq0ai} is false.
Thus there exists a sequence $(x_n)_{n\in\NNs}$ in~$E$ \st\
$||x_n||_E=1$ and $||T.x_n||\matrarrow_{n\rightarrow\infty}0$, \cf~\eref{eq2ai}.
And $\kappa$ being compact and $(x_n)_{n\in\NNs}$ being bounded, the sequence $(\kappa.x_n)_{n\in\NNs}$ 
has a convergent subsequence $(\kappa.x_\nk)_{k\in\NNs}$ that converges in the Banach space~$G$.
With $T$~continuous, $\kappa$~compact and the hypothesis~\eref{eq0aiPT}, we get
$$
\gamma ||x_{n_i}-x_{n_j}||_E\le ||T.x_{n_i}-T.x_{n_j}||_F+||\kappa.x_{n_i}-\kappa.x_{n_j}||_G
\matrarrow_{i,j\rar\infty} 0+0=0.
$$
Thus $(x_{n_k})_{k\in\NNs}$ is a Cauchy sequence in the Banach space~$E$, so converges to a limit $x \in E$.
Since $||T.x_\nk||\matrarrow_{k\rightarrow\infty}0$ and $T$ is continuous,
we get  $||T.x|| = 0$, thus $x=0$ since $T$ is injective.
But $||x_\nk||_E=1$ implies $||x||_E=1$. Absurd, thus~\eref{eq0ai} is true.
\finproof

%%%%%%%%%%%%%%%%%%%%%%%%%%%%%%%%%%%%%%%%%%%%%%%%%%%%%%%%%%%%%%%%%%%%%%%%%%%%%%%%%%%

\subsubsection{The range of $\vgrad : \Ldo \rar\Hmuon$ is closed}

We can now prove theorem~\ref{thmrgo}.
Let $T=\vgrad : \Ldo \rar \Hmuon$,
and $\kappa$ the canonical injection $\Ldo \rar \Hmuo$.
Since $T$ is linear continuous, \cf~\eref{eqspr1}, and $\kappa$ is compact, \cf\ Rellich theorem~\ref{thmReld},
the Petree--Tartar theorem~\ref{thmPT} implies that the range of~$T$ is closed.
\findem

%%%%%%%%%%%%%%%%%%%%%%%%%%%%%%%%%%%%%%%%%%%%%%%%%%%%%%%%%%%%%%%%%%%%%%%%%%%%%%%%%%%
%%%%%%%%%%%%%%%%%%%%%%%%%%%%%%%%%%%%%%%%%%%%%%%%%%%%%%%%%%%%%%%%%%%%%%%%%%%%%%%%%%%

\section{The closed range theorem}

The results and full proofs can be found \eg\ in Brézis~\cite{brezis}.

\comment{

\subsubsection{The open mapping theorem applied to $T'$}

$T'$ is defined in~\eref{eqdefadj}.

\deblem
\label{coromt}
If $T \in \calL(E;F)$ (linear continuous) is bijective, 
then $T' \in\calL(F';E')$ (linear continuous) is bijective, therefore
\be
\label{eq2aii}
\exists \gamma'>0 ,\; 
\forall \ell \in F',\quad ||T'.\ell||_{E'} \ge \gamma' ||\ell||_{F'}.
\ee
\comment{
If $T \in \calL(E;F)$ (linear continuous) is injective \st\ $\RTp$ is closed then
\be
\label{eq2aii2}
\exists \gamma'>0 ,\; 
\forall \ell \in F',\quad ||T'.\ell||_{E'} \ge \gamma' ||\ell||_{F'}.
\ee
If $T \in \calL(E;F)$ (linear continuous) is \st\ $\RTp$ is closed then
\be
\label{eq2aii3}
\exists \gamma'>0 ,\;
\forall \ell \in F',\quad ||T'.\ell||_{E'} \ge \gamma' ||\ell||_{F'/\Ker(T')}.
\ee
}
\finlem

\debdem
$T'$ is linear and continuous, \cf~\eref{eqnTp}.

Injectivity: If $\ell \in F'$ and $T'.\ell = 0$, then
$\la T'.\ell,x\ra_{E',E} = 0 = \la \ell,T.x\ra_{F',F}$ for all $x\in E$,
thus $0 = \la \ell,y\ra_{F',F}$ for all $y\in F$ (bijectivity), thus $\ell=0$.

Surjectivity: Let $k \in E'$;
We look for $\ell \in F'$ \st\ $k=T'.\ell$, that is \st\ $k=\ell.T$
since $\la k,x\ra = \la T'.\ell,x\ra = \la \ell,T.x\ra = \ell(T.x) = (\ell \circ T)(x)$ for all~$\ell$
and $\ell\circ T$ is denoted $\ell.T$ (with linear maps).
So, $T$ being bijective, $\ell := k.T^{-1}$ fits.
Check:
$\la T'.\ell , x \ra_{E',E}
= \la \ell, T.x\ra_{F',F}
= \la k.T^{-1}, T.x\ra_{F',F}
= (k.T^{-1}).(T.x) =k.x$
for all $x\in E$, thus $T'.\ell = k$.

And $E'$ and $F'$ are Banach spaces.
Therefore~\eref{eq2ai} applied to~$T'$ gives~\eref{eq2aii}.
\comment{
And~\eref{eq2ai2} applied to~$T'$ gives~\eref{eq2aii2}.

An~\eref{eq2ai2z} applied to~$T'$ gives~\eref{eq2aii3}.
}
\findem

}

%%%%%%%%%%%%%%%%%%%%%%%%%%%%%%%%%%%%%%%%%%%%%%%%%%%%%%%%%%%%%%%%%%%%%%%%%%%%%%%%%%%

\subsubsection{The closed range theorem}

Let $T\in L(E;F)$, so $T' \in \calL(F';E')$, \cf~\eref{eqdefadj}. We have
\be
\label{eqcr1}
\eqalign{
%& \Im(T') = \{\ell \in E' \;\hbox{\st}\; \exists m\in F',\;\ell = T'.m\} \subset E', \cr
&\Ker(T') = \{m\in F' \;\hbox{\st}\; T'.m=0\} = \{m\in F' \;\hbox{\st}\; m.T=0\} \subset F',
}
\ee
since $m\in \Ker(T')$ $\Leftrightarrow$ $T'.m=0$
$\Leftrightarrow$ $\la T'.m,x\ra_{E',E} = 0 = \la m, T.x \ra_{F',F} = m(T.x) = (m\circ T)(x)$ for all $x\in E$
$\Leftrightarrow$ $m.T=0$, with $m.T$ the notation of $m\circ T$ when the maps are linear maps.

If $M\subset E$ is a linear subspace in~$E$ then the dual orthogonal of~$M$ is
\be
\label{eqMo}
M^o := \{\ell\in E' : \la\ell,x\ra_{E',E}=0,\;\forall x\in M\} \quad (\subset E')
\ee
(the subspace of $E'$ of linear forms vanishing on~$M$).
$M^o$ is a linear subspace in~$E'$ (trivial).
And $M^o$ is closed in~$E'$:
Indeed if $(\ell_n)_\NNs$ is a Cauchy sequence in~$M^o$, so
$||\ell_n- \ell_m||_{E'} \matrarrow_{n,m\rar\infty} 0$,
then, for $x\in E$ the sequence $(\ell_n(x))$ is a Cauchy sequence in~$\RR$, thus convergence
toward a real named~$\ell(x)$; This defines a function $\ell:E\rar \RR$.
And $\ell(x_1+\lambda x_2)
= \lim_{n\rar\infty} \ell_n(x_1+\lambda x_2)
= \lim_{n\rar\infty} \ell_n(x_1)+\lambda \lim_{n\rar\infty}(x_2)
=\ell(x_1) + \lambda \ell(x_2)$, thus $\ell$ is linear,
and, for $x\in E$, $\ell$ is continuous at~$x$ since
$|\ell.x' - \ell.x| \le |(\ell-\ell_n).x' + (\ell-\ell_n).x| + |\ell_n.x' - \ell_n.x|
\le(||\ell-\ell_n|| + ||\ell_n||) ||x'-x||_E$ with $||\ell_n||\le ||\ell_N||+1$ for $N$ large enough
and $n\ge N$.

If $N\subset F'$ is a linear subspace in~$F'$ then let
\be
\label{eqNo}
N^\perp := \{y \in F : \la m,y\ra_{F',F}=0,\;\forall m \in N\} \quad (\subset F).
\ee
Then $N^\perp$ is a linear subspace in~$F$ (trivial) that is closed in~$F$ (similar proof than for $M^o$).
To be compared with, \cf~\eref{eqMo},
\be
\label{eqNo2}
N^o = \{y'' \in F'' : \la y'',m\ra_{F'',F'}=0,\;\forall m \in N\} \quad (\subset F'').
\ee

\debrem
If $F$ is reflexive, that is $F''\simeq F$ (identification) then $N^o \simeq N^\perp$ (identification).
Indeed, with $J$ the canonical isomorphism given in~\eref{eqreflex}, 
if $y \in N^\perp$ then let $y''=J(y)\in F$,
so for all $m \in N$ we have $0 = m.y = y''.m$, and thus $y''\in N^o$;
And  if $y'' \in N^o$ then let $y\in F$ \st\ $J(y)=y''$ (thanks to the reflexivity),
then for all $m \in N$ we have $0=y''.m = m.y$, and thus $y\in N^\perp$.
\finrem

\begin{theorem}[Closed range theorem]
\label{thmcr}
Let $E$ and $F$ be Banach spaces and $T\in \calL(E;F)$ (linear and continuous).
Then %$T'\in \calL(F',E')$ (that is, the dual operator is linear and continuous), and
the following properties are equivalent:

(i) $\RT$ is closed in~$F$,

(ii) $\RTp$ is closed in~$E'$,

(iii) $\RT = \Ker(T')^\perp$,

(iv) $\RTp = \Ker(T)^o$.
\finthm

We then deduce, with~\eref{eq2ai2z}:

\debcor
If $\RT$ is closed in~$F$ then $\Im(T')$ is closed in~$E$', thus
\be
\label{eq2ai2z2}
\exists \gamma'>0,\quad \forall \ell\in F',\quad ||T'.\ell||_{E'} \ge  \gamma'\,||\ell||_{F'/\Ker(T')} .
\ee
\fincor

\debproof
The full proof of theorem~\ref{thmcr} (even for unbounded operators with dense domain of definition) can be found \eg\ in Brézis~\cite{brezis}
or Yosida~\cite{yosida} (for locally convex spaces that are metrizable and complete).
We give here the proof in the simplified case of $T$ a linear continuous mapping between two Banach spaces
(sufficient for our needs).
We need some lemmas:
%\finproof

\deblem
If $E$ is a Banach space and $M$ is a linear subspace in~$E$, then
\be
%M \subset 
\overline M = (M^o)^\perp .
\ee
\finlem

\debdem
If $x \in M$ then
$\la \ell,x\ra_{E',E} = 0$ for all $\ell\in M^o$,
thus $x\in (M^o)^\perp$, \cf~\eref{eqNo}.
And $(M^o)^\perp$ being closed we get $\overline M \subset (M^o)^\perp$.

Conversely: Suppose $x_0\in (M^o)^\perp$ and $x_0\notin \overline M$; 
Then $\{x_0\}$ being compact and $\overline M$ being closed and convex (it is a linear subspace),
there exists a hyperplane that strictly separates $x_0$ and~$\overline M$ (geometric form or the Hahn--Banach theorem),
that is there exists $\ell \in E'$ and $\alpha\in \RR$ \st\
$\la \ell,x \ra_{E',E} < \alpha < \la \ell,x_0 \ra_{E',E} $ for all $x\in M$.
And $M$ being a linear space, taking $-x\in M$, it follows that $\la \ell,x \ra_{E',E} = 0$ for all $x\in M$,
thus $\ell \in M^o$. And $\la \ell,x \ra_{E',E} = 0$ for all $x\in M$ implies $\la \ell,x_0 \ra_{E',E} >0$
with $x_0\notin \overline M$, thus $\ell\notin M^o$. Absurd, thus $(M^o)^\perp \subset \overline M$.
\findem

\deblem
If $G$ and $L$ are two closed subspaces in a Banach space, then
\be
\label{eqG1}
G \cap L = (G^o + L^o)^\perp, \qand
G^o \cap L^o = (G+L)^o.
\ee
\finlem

\debdem
\eref{eqG1}$_1$:
If $x\in G \cap L$, and if $m = g + \ell \in G^o + L^o$,
then $m.x = g.x + \ell.x = 0+0$, thus $x\in (G^o + L^o)^\perp$.

Conversely, we have $(G^o + L^o)^\perp \subset (G^o)^\perp$
(particular case of: If $Y \subset Z$ then $Z^\perp \subset Y^\perp$),
and $(G^o)^\perp = G$ since $G$ is closed,,
thus if $x\in (G^o + L^o)^\perp$ then $x\in G$; And similarly $x\in L$,;
Thus $x\in G \cap L$.

Similar proof for~\eref{eqG1}$_2$.
\findem

Let $E\times F$ be equipped with the (usual) norm $||(x,y)||_{E\times F} = \max(||x||_E,||y||_F)$,
so $E\times F$ is a Banach space.

\deblem
If $T$ is continuous, then its graph
\be
G(T) = \{(x,y)\in E\times F\hbox{ \st\ } \exists x\in E,\; y=T.x\}=\{(x,T.x)\in E\times F\}
\ee
is closed in $E\times F$.
\finlem

\debdem
If $((x_n,T.x_n)))_\NNs$ is a Cauchy sequence in~$G(T)$, then,
$E$ being a Banach space, $(x_n)_\NNs$ converges toward a $x \in E$, thus,
$T$ being continuous, $T.x_n$ convergence toward $T.x \in F$, so $(x,T.x)\in G(T)$.
\findem

We have 
\be
\label{eqGTp}
G(T')=\{(m,T'.m)\in F'\times E'\}.
\ee

\deblem
If $T$ is continuous, then
$G(T')$ is closed in $F'\times E'$, and
\be
\label{eqGTo0}
(m,\ell) \in G(T') \; \Longleftrightarrow\; (-\ell,m)\in G(T)^o.
\ee
\finlem

\debdem
Let $((m_n,T'.m_n))_\NNs$ be a sequence in~$G(T)$ \st\
$(m_n,T'.m_n) \matrarrow_{n\rar\infty} (m,z)\in F'\times E'$.
Thus $m_n \matrarrow m\in F'$ and $T'.m_n \matrarrow_{n\rar \infty} k\in E'$,
that is $\la T'.m_n,x \ra_{E',E} \matrarrow_{n\rar \infty} \la k,x \ra_{E',E}\in \RR$ for all $x\in E$.
And we have to check that $k=T'.m$.
For all $x\in E$, we have $\la T'.m_n,x \ra_{E',E} = \la m_n,T.x \ra_{F',F}$,
thus $\la m_n,T.x \ra_{F',F} \matrarrow_{n\rar \infty} \la k,x \ra_{E',E}$,
\ie\ $\la T'm_n,x \ra_{F',F} \matrarrow_{n\rar \infty} \la k,x \ra_{E',E}$,
thus $T'm_n, \matrarrow_{n\rar \infty} k \in F'$. So $G(T')$ is closed.

$(m,\ell) \in G(T')$
$\Leftrightarrow$ 
$\ell=T'.m$
$\Leftrightarrow$ 
$\la \ell,x \ra_{E',E} = \la T'.m,x \ra_{E',E} = \la m,T.x\ra_{E',E}$ for all $x\in E$
$\Leftrightarrow$ 
$\la \ell,x \ra_{E',E} - \la m,T.x\ra_{E',E} = 0$ for all $x\in E$
$\Leftrightarrow$ 
$\la (\ell,-m),(x,T.x) \ra_{E'\times F',E\times F} = 0$ for all $x\in E$
$\Leftrightarrow$
$(\ell,-m) \in G(T)^o$. % (and $G(T)^o$ is a linear space).
\findem

Define
\be
L := E\times \{0\}.
\ee
$L$ is closed in $E\times F$ since $E$ and $\{0\}$ are, and
\be
\label{eqGo}
L^o = \{0\} \times F'.
\ee
Indeed $L^o
= \{(\ell,m) \in E'\times F' : \la(\ell,m),(x,0)\ra_{E'\times F',E\times F} = 0,\, \forall x \in E\}
= \{(\ell,m) \in E'\times F' : \la\ell,x\ra_{E',E} +0
= 0,\, \forall x \in E\} = \{0\} \times F'$.

\deblem
\be
\label{eqGK1}
\Ker(T) \times \{0\} = G(T) \cap L
\ee
\be
\label{eqGK2}
E \times \RT = G(T)+L
\ee
\be
\label{eqGK3}
\{0\} \times \Ker(T') = G(T)^o \cap L^o
\ee
\be
\label{eqGK4}
\Im(T') \times F' = G(T)^o + L^o
\ee
\finlem

\debdem
$(x,y)\in \Ker(T) \times \{0\}$ iff $Tx=0$ and $y=0$;
And $(x,y)\in G(T)\cap L$ iff $y=Tx$ and $y=0$, thus~\eref{eqGK1}.

$(x_1,y_1) \in E\times \RT$ iff $(x_1,y_1) = (x_1,T.x_1')$ for some $x_1'\in E$;
%$x\in E$ and $\exists x'\in E$ \st\ $y=T.x'$,
And $(x_2,y_2) \in G(T)+L$ iff $\exists x_2'\in E$ and $\exists x_2''\in E$ \st\
$(x_2,y_2)= (x_2',T.x_2') + (x_2'',0) = (x_2'+x_2'',T.x_2') = (x_3,T.(x_3-x_2'))$, thus~\eref{eqGK2}.

$(\ell,m) \in \{0\} \times \Ker(T')$ iff $\ell=0$ and $m \in \Ker T'$;
And $(\ell,m) \in G(T)^o \cap L^o$
iff $(-m,\ell)\in G(T')$, \cf~\eref{eqGTo0}, and $(\ell,m) \in L^o$,
\ie\ iff $\ell = -T'.m$ and $\ell=0$,
\ie\ iff $\ell=0$ and $m \in \Ker T'$, thus~\eref{eqGK3}.

$(\ell,m) \in \Im(T') \times F'$ iff $\exists k\in F'$ \st\ $\ell = T'.k$ and $m\in F'$;
And $(\ell,m)\in G(T)^o + L^o$ iff $\exists (\ell_1,m_1)\in G(T)^o$ and $\exists (\ell_2,m_2) \in L^o$
\st\ $\ell = \ell_1+\ell_2$ and $m=m_1+m_2$,
\ie, with~\eref{eqGTo0} and~\eref{eqGo},
iff $\exists m_1\in F'$ (and then $-\ell_1=T'.m_1$) and $m_2 \in F'$
\st\ $\ell = -T.m_1 +0$ and $m=m_1+m_2$,
\ie\ iff $\ell \in \Im(T')$ and $m\in F'$, thus~\eref{eqGK4}.
\findem

\debcor
\label{corthmcr}
\be
\label{eqct1}
\Ker(T) = \Im(T')^\perp,
\ee
\be
\label{eqct2}
\Ker(T') = \RT^o ,
\ee
\be
\label{eqct3}
(\Ker(T))^o = \overline{\Im(T')},
\ee
\be
\label{eqct4}
\Ker (T')^\perp = \overline{\RT}.
\ee
\fincor

\debdem
\eref{eqGK4} gives $R(T')^\perp \times \{0\} = (G(T)^o + L^o)^\perp = G(T) \cap L$, \cf~\eref{eqG1},
thus $= \Ker(T) \times \{0\}$, \cf~\eref{eqGK1}, thus~\eref{eqct1}. Thus~\eref{eqct3}.

\eref{eqGK2} gives $\{0\} \times \RT^o = (G(T) + L)^o = G^o \cap L^o$, \cf~\eref{eqG1},
thus $= \{0\} \times \Ker (T')$, \cf~\eref{eqGK3}, thus~\eref{eqct2}. Thus~\eref{eqct4}.
\findem

\noindent
{\bf Proof of theorem~\ref{thmcr}:} apply corollary~\ref{corthmcr}.

\comment{

%%%%%%%%%%%%%%%%%%%%%%%%%%%%%%%%%%%%%%%%%%%%%%%%%%%%%%%%%%%%%%%%%%%%%%%%%%%%%%%%%%%

\subsubsection{Complement}

In
\be
\RTp = \calL(\RT;\RR)
\ee
the dual norm $||.||_\RTp$ is defined by, for $\ell\in \RTp$,
\be
||\ell||_\RTp = \sup_{y\in\RT \atop ||y||_\RT\le 1} |\la \ell,y\ra_{\RTp,\RT}|
\quad (= \sup_{y\in\RT \atop ||y||_F\le 1} |\la \ell,y\ra_{F',F}|)
.
\ee
Thus, if $\RT$ is closed in~$F$ then $(\RTp,||.||_\RTp)$ is a Banach space since $(\RT,||.||_\RT)$ is.

If $M\subset F$ is a linear subspace in~$F$ then let
\be
\label{eqMo2}
M^o := \{\ell\in F' : \la\ell,y\ra_{F',F}=0,\;\forall y\in M\} \quad (\subset F'),
\ee
and $M^o$ is called the dual orthogonal of~$M$
(the subspace of $F'$ of linear forms vanishing on~$M$).
Then $M^o$ is a linear subspace in~$F'$ (trivial).
And $M^o$ is closed in~$F'$:
Indeed if $(\ell_n)_\NNs$ is a Cauchy sequence in~$M^o$, so
$||\ell_n- \ell_m||_{F'} \matrarrow_{n,m\rar\infty} 0$,
then, for $y\in F$ the sequence $(\ell_n(y))$ is a Cauchy sequence in~$\RR$, thus convergence
toward a real named~$\ell(y)$; This defines a function $\ell:F\rar \RR$.
And $\ell(y_1+\lambda y_2) = \lim_{n\rar\infty} \ell_n(y_1+\lambda y_2)$ and the $\ell_n$ are linear forms,
and $\ell$ is continuous since
$|\ell.y_2 - \ell.y_1| \le |(\ell-\ell_n).y_2 + (\ell-\ell_n).y_1| + |\ell_n.y_2 - \ell_n.y_1|
\le(||\ell-\ell_n|| + ||\ell_n||) ||y_2-y_1||_F$ and $\ell_n$ is continuous.

Let 
\be
Z : 
\left\{\eqalign{
F' & \rar \RT' = \calL(\RT;\RR) \cr
\ell & \rar Z(\ell) = \ell_{|\RT},
}\right.
\ee
that is $Z(\ell)(y) = \ell.y$ for all $y\in \RT$,
so $Z(\ell)(T.x) = \ell.(T.x)$ for all $x\in E$.

\deblem
If $\RT$ is closed, then
$Z$ is linear, continuous,
$\Ker(Z) = \Ker(T') = \RT^o$, %($=\{\ell\in F' : \la \ell,y\ra_{F',F} = 0,\; \forall y\in \RT\}$).
and $\Im(Z) = \RTp$, \ie\ $Z$ is surjective. 
\finlem

\debdem
Linearity: $Z(\ell+\lambda m) = (\ell+\lambda m)_{|\RT} = \ell_{|\RT}+\lambda m_{|\RT} = Z(\ell)+\lambda Z(m)$.

Continuity: $||Z(\ell)||_\RTp = ||\ell_{|\RT}||_\RTp
=\sup_{y\in \RT,\atop ||y||_F\le 1}|\ell(y)| \le ||\ell||_{F'}$, thus $||Z||\le 1$.

Next, $\ell\in \Ker Z$ iff $Z(\ell)=0$,
\ie\ iff, for all $x\in E$, $\ell.(T.x)=0 = \la \ell,T.x\ra_{F',F} = \la T'.\ell,x\ra_{E',E}$,
\ie\ iff $T'.\ell=0$,
\ie\ iff $\ell \in \Ker(T')$.

And $\ell \in \Ker(T')$ iff $T'.\ell=0$ ($\in E'$),
\ie\ iff, for all $x\in E$, $0 = \la T'.\ell,x \ra_{E',E} = \la \ell,T.x\ra_{F',F}$,
\ie\ iff $\ell\in \RT^o$.

Surjectivity:
We apply Hahn--Banach: % theorem (extension of a linear form defined on a linear subspace):
$\RT'$ is a linear subspace in~$F'$, then for $m:\RT' \rar \RR$ linear and continuous,
it exists $\ell\in F'$ that extends~$m$, and in particular $\ell.y = m.y$ for all $y\in \RT$,
\ie\ $Z(\ell) = m$. Thus $Z$ is surjective.
\findem

\comment{
\deblem
If $\ell\in E'$, $m\in F'$ and $T\in \calL(E;F)$ then $T'\in\calL(F';E')$ satisfies
\be
\label{eqltm}
\ell = T'.m \; \Longleftrightarrow\; \ell = m.T.
\ee
(We used linearity to denote $T' \circ m = T.m$ and $m \circ T = m.T$.) Thus
\be
\ell \in \Im(T') \; \Longleftrightarrow\; \ell\in E'\;\hbox{ and }\;\exists m\in F',\; \ell = m.T.
\ee
\finlem

\debdem
$k=T'.\ell$ $\Leftrightarrow$
$k.x = \la k,x\ra_{E',E} = \la T'.\ell,x\ra_{E',E} = \la \ell,T.x\ra_{F',F}
=\ell.(T.x)$ for all $x\in E$
$\Leftrightarrow$ $k=\ell.T$.
\findem
}

%%%%%%%%%%%%%%%%%%%%%%%%%%%%%%%%%%%%%%%%%%%%%%%%%%%%%%%%%%%%%%%%%%%%%%%%%%%%%%%%%%%

\subsubsection{Corollary: The closed range case}

\debcor
If $T$ is injective %and if the range $\RT$ is closed in~$F$, 
then $T'$ is surjective, that is $\Im(T') = E'$.
\fincor

\debdem
Let $\ell\in E'$. We look for $m\in F'$ \st\ $T'.m = \ell$,
that is $\la T'.m,x \ra = \la m,T.x \ra = \la \ell,x\ra$ for all $x\in E$.

\findem

\debcor
and $T'$ is bijective from $\RT/\Ker(T')$ to~$E'$.
\fincor

\debdem
\def\tell{{\tilde\ell}}
$G=(\RT,||.||_F)$ is a Banach space,
thus $\tT : E\rar G$ defined by $\tT(x)=T(x)$ is linear continuous bijective, \cf\ corollary~\ref{coromt},
thus $\Im(\tT')=E'$.

And $\Im(\tT') = \{m\in F',\; \exists\}$

And $\Im(T')=\Im(\tT')$: Indeed,
1- If $m \in \Im(T')$ then $\exists \ell\in E'$ \st\ $T'.\ell =m$,
that is, $\la \ell,T.x\ra_{E',E} = \la m,x\ra_{E',E}$ for all $x\in E$

And for $\ell\in F'$ we have $\la T'.\ell,x\ra_{E',E} = \la \ell,T.x \ra_{F',F}$.
Then let $\tell \in G'$ be defined by $\la \tell , y \ra_{F',F} = \la \ell , y \ra_{F',F}$ for all $y\in \RT$.

And $T'$ is continuous, thus $\Ker(T')$ is closed in~$F'$, and
$T': F'/\Ker(T') \rar \RTp$ is bijective.

Then $k\in \Ker(T')$ iff $T'.k=0$, \ie\ $\la T'.k,x\ra_{E',E} = 0 = \la k,T.x\ra$ for all $x\in E$.

\findem

\def\tTp{{\widetilde{T'}}}

\deblem
$\RTp$ is a closed subspace in~$E'$.
\finlem

\debdem

\findem

Next define
$\tTp : 
\left\{\eqalign{
F'/\Ker(T') & \rar \RTp \cr
\dot m &\rar \tTp(\dot m) := T'(m)
}\right\}
$ when $m\in \dot m$, thus characterized by $\tTp \circ \pi = T'$.

\debthm
If $\RT$ is closed in $F$, then $\tTp$ is linear, continuous, injective, surjective. %, and $||\tT.m||_{\RTp} = ||m||_{F'/\Ker(T')}$.
Thus,
\be
\label{eq2ai2zp}
\exists \gamma>0 ,\;
\forall m\in F',\quad ||T'.m||_{\RTp} \ge \gamma ||m||_{F'/\Ker(T')}.
\ee
\finthm

\debdem
$\tTp$ is linear since $T'$ is.
Next
$\tTp(\dot m)=0 = T'(m)$ implies $m\in \Ker(T')=\dot 0$, thus $\dot m=\dot0$ and $\tTp$ is injective,
and $\tTp$ is continuous and surjective since $Z$~is. And $\RTp$ is closed since

Thus, \cf~\eref{eq2ai2z},
\be
%\label{eq2ai2zp}
\exists \gamma>0 ,\;
\forall \dot m\in F'/\Ker(T'),\quad ||\tTp.\dot m||_{\RTp} \ge \gamma ||\dot m||_{F'/\Ker(T')}.
\ee

\findem

If $N\subset E'$ is a linear subspace in~$E'$ then let
\be
\label{eqNo}
N^\perp := \{x \in E : \la\ell,x\ra_{E',E}=0,\;\forall \ell \in N\} \quad (\subset E)
\ee
Then $N^\perp$ is a linear subspace in~$E$ (trivial) that is closed in~$E$ (similar proof than for $M^o$).
To be compared with % $N^o \in E''$:
\be
\label{eqNo2}
N^o = \{x'' \in E'' : \la x'',\ell\ra_{E'',E'}=0,\;\forall \ell \in N\} \quad (\subset E'')
\ee

%%%%%%%%%%%%%%%%%%%%%%%%%%%%%%%%%%%%%%%%%%%%%%%%%%%%%%%%%%%%%%%%%%%%%%%%%%%%%%%%%%%

\subsubsection{Application to mixed problems}

Let $E$ and $F$ be Banach spaces and let $B \in \calL(E;F')$ (linear and continuous).
We then have $B' \in \calL(F'';E')$. In particular, if $F''$ is reflexive,
that is, if $F'' \simeq F$ \cf\ definition~\ref{defrefl},  then $B' \in \calL(F;E')$.

\debcor
\label{coromt}
If $B \in \calL(E;F')$ (linear continuous) is bijective, then
\be
\label{eq2aiB1}
\exists \gamma>0 ,\; 
\forall x\in E,\quad ||B.x||_{F'} \ge \gamma ||x||_E.
\ee
Moreover, if $F$ is reflexive, then $B'\in \calL(F;E')$ is bijective and
\be
\label{eq2aiB2}
\exists \gamma'>0 ,\; 
\forall y\in F,\quad ||B'.y||_{E'} \ge \gamma' ||y||_F.
\ee
\fincor

\debdem
\eref{eq2aiB1} is given by~\eref{eq0ai}.
And~\eref{eq2aii} and $B'\in\calL(F'';E')$ give $||B'.y||_{E'} \ge \gamma' ||y||_{F''}$ for all $y\in F''$.
Thus is $F$ is reflexive we get~\eref{eq2aiB2}.
\findem

\debcor
If $\Im(B)$ is closed in~$F'$,
then
\be
\label{eq2aiB2}
\exists \gamma>0 ,\quad
\forall x\in E,\quad ||B.x||_{F'} \ge \gamma ||x||_{E/\Ker(B)}.
\ee
Moreover $F$ is reflexive, then $\Im(B')$ is closed in~$E'$ and
\be
\label{eq2aiB22}
\exists \gamma'>0 ,\quad
\forall y \in F,\quad ||T'.y||_{E'} \ge \gamma' ||y||_{F/\Ker(B')} .
\ee
\fincor

\debdem
$F'$ being a Banach space (with its usual norm \cf~\eref{eqdefnell}),
\eref{eq2ai2} gives~\eref{eq2aiB2}.

Then define $\tT : D \rar G$ by $\tT(\dot x) = T(y)$
for any $y = x+x_0 \in \dot x=x+\Ker(T)$.
Thus $\tT$ bijective (injective and surjective by definition of the kernel and of the range).
And $\tT$ is linear (trivial) and continuous (with $||\tT|| = ||T||$).
Thus \eref{eq2ai} applied to $\tT$ gives~\eref{eq2ai2}.

Then $\tT ' : G' = \calL(G;\RR) \rar D' = \calL(D;\RR)$ is linear continuous and bijective since $\tT$~is,
\cf\ corollary~\ref{coromt}, with $||\tT'.\ell||_{D'} \le {1\over \gamma'} ||\ell||_F$ for any $\ell \in G'$,
so

$||\tT '.\ell||_{D'} \ge \gamma' ||\ell||_{G'}$.
\findem

}

%%%%%%%%%%%%%%%%%%%%%%%%%%%%%%%%%%%%%%%%%%%%%%%%%%%%%%%%%%%%%%%%%%%%%%%%%%%%%%%%%%%

\section{A well-posed mixed problem}

\def\Vs{{V'}}
\def\Qs{{Q'}}

%%%%%%%%%%%%%%%%%%%%%%%%%%%%%%%%%%%%%%%%%%%%%%%%%%%%%%%%%%%%%%%%%%%%%%%%%%%%%%%%%%%

\subsection{Notations}

Let $V$ and $Q$ be two Banach spaces.
Let $b\dd : V\times Q \rar \RR$ be a bilinear form.
$b\dd$ is said to be continuous (or bounded) iff
\be
\label{eqbcont0}
\exists c>0,\quad \forall (v,q)\in B_V(0,1) \times B_Q(0,1),\quad |b(v,q)| \le c.
\ee
Then let
\be
\label{eqbcont1}
||b|| := \sup_{v\in B_V(0,1)\atop q\in B_Q(0,1)} |b(v,q)|.
\ee
And let $\calL(V,Q;\RR)$ be the space of bilinear and continuous forms
with its (usual) norm given by~\eref{eqbcont1}.

If $b\dd \in \calL(V,Q;\RR)$ (bilinear and continuous), then define
\be
B:
\left\{\eqalign{
V & \rar Q' \cr
v & \rar Bv % \qst \forall q\in Q,\;  = b(v,q).
}\right\}
\qand
B^t:
\left\{\eqalign{
Q & \rar V' \cr
q & \rar B^tq
}\right\}
\ee
by
\be
b(v,q) = \la Bv,q \ra_{Q',Q} =  \la B^tq , v \ra_{V',V}.
\ee
Thus $B$ and $B^t$ are linear (trivial) and continuous with 
\be
||B|| = ||B^t|| = ||b||.
\ee
Indeed
$||Bv||_\Vs
= \sup_{q\in B_Q(0,1)} |\la Bv,q \ra_{Q',Q}| 
= \sup_{q\in B_Q(0,1)} |b(v,q)| 
\le  \sup_{q\in B_Q(0,1)} ||b||\,||v||_V||q||_Q = ||b||\,||v||_V$
gives $||B|| \le ||b||$ (continuity), and
$|b(x,y)| = |\la Bv,q \ra_{Q',Q}| \le ||Bx||_\Qs ||y||_Q \le ||B||\,||x||_V ||y||_Q$ gives $||b|| \le ||B||_\Vs$.
So $B \in \calL(V;Q')$. Idem pour $B^t$.

Suppose $Q$ reflexive, \cf\ definition~\ref{defrefl}, then the dual $B'\in\calL(Q'';V')$ of $B\in\calL(V;Q')$,
defined by $\la B'v,\ell\ra_{V',V} = \la v,B\ell\ra_{Q'',Q}$ for all $v\in Q''$ and $\ell\in V'$ \cf~\eref{eqdefadj}, is identified to~$B^t$:
\be
\label{eqQr}
\calL(Q'';V') \ni B' \simeq B^t \in \calL(Q;V').
\ee

Suppose $V$ reflexive, \cf\ definition~\ref{defrefl}, 
then the dual $(B^t)'\in\calL(V'';Q')$ of $B^t\in\calL(Q;V')$,
defined by $\la (B^t)'v,q\ra_{Q',Q} = \la v,B^t\ell\ra_{V'',V'}$ for all $v\in Q''$ and $\ell\in V'$
\cf~\eref{eqdefadj}, is identified to~$B^t$:
\be
\label{eqVr}
\calL(V'';Q') \ni (B^t)' \simeq B \in \calL(V;Q').
\ee

%%%%%%%%%%%%%%%%%%%%%%%%%%%%%%%%%%%%%%%%%%%%%%%%%%%%%%%%%%%%%%%%%%%%%%%%%%%%%%%%%%%

\subsection{The mixed problem}

Let $a\dd : V\times V \rar \RR$ and $b\dd : V\times Q \rar \RR$ be bilinear forms.
Let $f\in V'$ and $g\in Q'$ (linear forms).
A mixed problem is a problem of the type: Find $(u,p) \in V\times Q$ \st\
\be
\label{eqw0}
\left\{\eqalignrll{
& a(u,v) + b(v,p) &= \la f,v \ra_{V',V}, \quad \forall v \in V , \cr
& b(u,q) &= \la g,q\ra_{Q',Q}, \quad \forall q \in Q,
}\right.
\ee
\cf~\eref{eqgenw0}.

%%%%%%%%%%%%%%%%%%%%%%%%%%%%%%%%%%%%%%%%%%%%%%%%%%%%%%%%%%%%%%%%%%%%%%%%%%%%%%%%%%%

\subsection{The inf-sup conditions}

For the existence (and control) of~$p$, we suppose that the range $\Im(B^t)$ of $B^t: Q\rar V'$ is closed,
that is, \cf~\eref{eq2ai},
\be
\label{eqism1}
\exists \beta>0 ,\; \forall q\in Q,\; ||B^tq||_{V'} \ge \beta  ||q||_{Q/\Ker(B^t)},
\ee
\ie,
\be
\label{eqism2}
\exists \beta>0 ,\;\forall q\in Q,\;
\sup_{v \in V}{ b(v,q) \over ||v||_{V/\Ker(B)} } \ge \beta  ||q||_{Q/\Ker(B^t)},
\ee
also written as the inf-sup condition
$\inf_{q\in Q}\sup_{v \in V}{ b(v,q) \over ||v||_{V/\Ker(B)} ||q||_{Q/\Ker(B^t)} } \ge \beta$.

For the existence (and control) of~$u$, we suppose the range $\Im(B)$ of $B:V\rar Q'$ is closed,
that is, we suppose, \cf~\eref{eq2ai},
\be
\label{eqism3}
\exists \beta>0 ,\; \forall v\in V,\; ||Bv||_{Q'} \ge \beta  ||v||_{V/\Ker(B)},
\ee
\ie,
\be
\label{eqism4}
\exists \beta>0 ,\;\forall v\in V,\;
\sup_{q \in Q}{ b(v,q) \over ||q||_{Q/\Ker(B^t)} } \ge \beta  ||v||_{V/\Ker(B)},
\ee
also written as the inf-sup condition
$\inf_{v\in V}\sup_{q \in Q}{ b(v,q) \over ||v||_{V/\Ker(B)} ||q||_{Q/\Ker(B^t)} } \ge \beta$.

\medskip
Remark: With~\eref{eqQr} or~\eref{eqVr}, the reflexivity of $Q$ or~$V$ gives that
\eref{eqism3} implies~\eref{eqism1} or \eref{eqism1} implies~\eref{eqism3}.

\comment{
\begin{corollary}[Inf-sup condition]
\label{coris0}
If $E$ and $F$ are Banach spaces,
if $B \in \calL(E;F')$ (linear and continuous), and if $\Im(B)$ is closed, then
\be
\label{eqcig0}
\exists \gamma>0 ,\; 
\forall x\in E ,\; 
\sup_{\ell\in F'}{|\la \ell,T.x\ra_{F',F}| \over ||\ell||_{F'} }  \ge \gamma ||x||_E,
\ee
that is,
\be
\label{eqcig}
\exists \gamma>0 ,\; 
\inf_{x\in E}(\sup_{\ell\in F'}
{|\la \ell,T.x\ra_{F',F}| \over ||x||_E ||\ell||_{F'}} ) \ge \gamma.
\ee
And if $F$ is reflexive, that is $F''\simeq F$ (identification with~$J$, \cf~\eref{eqreflex}),
% in particular if $E$ is a Hilbert space, 
then
\be
\label{eqcig1bc}
\exists \gamma'>0 ,\; \forall \ell\in F',\;
\sup_{x\in E} {|\la T'.\ell,x\ra_{E',E}| \over ||x||_E} \ge \gamma'  ||\ell||_{F'},
\ee
that is,
\be
\label{eqcig1}
\exists \gamma'>0 ,\; \inf_{\ell\in F'}(\sup_{x\in E}
{|\la \ell,T.x\ra_{F',F}| \over ||x||_E ||\ell||_{F'}} ) \ge \gamma'.
\ee
\fincor

}

%%%%%%%%%%%%%%%%%%%%%%%%%%%%%%%%%%%%%%%%%%%%%%%%%%%%%%%%%%%%%%%%%%%%%%%%%%%%%%%%%%%

\subsection{The theorem for mixed problem}

\debthm.
\label{thmis}
{\bf Hypotheses: }
(i) $(V,\dd_V)$ is a Hilbert space, $(Q,||.||_Q)$ is a reflexive Banach space, $f\in V'$, and $g\in Q'$.

(ii) The bilinear form $a\dd$ is continuous on $V$, \cf~\eref{eqbcont0}, and coercive on~$\Ker(B)$,
that is,
\be
\label{eqcic}
\exists \alpha>0, \; \forall v\in \Ker(B),\quad a(v,v) \ge \alpha ||v||_V^2.
\ee

(iii) The bilinear form $b\dd$ is continuous on $V\times Q$, \cf~\eref{eqbcont0}, 
and $B$ is surjective (= onto), so we have~\eref{eqism3} and then~\eref{eqism1} since $Q$ is reflexive.

\comment{
so,
\be
\label{eqcis}
\exists \beta>0, \; 
\left\{\eqalign{
& \forall v\in V,\; \exists q\in Q,\quad  | b(v,q) | \ge \beta ||v||_{V/\Ker(B)}||q||_Q, \cr
& \forall q\in Q',\; \exists v\in V,\quad  | b(v,q) | \ge \beta ||q||_{Q/\Ker(B^t)}||v||_V, \cr
}\right.
\ee
}

\medskip
\noindent
{\bf Conclusion: }
Problem~\eref{eqw0} has a unique solution $(u,p) \in V \times Q/\Ker B^t$ that
depends continuously on~$f$ and~$g$, and more precisely,
with $C_a = (1 + {||a|| \over \alpha})$,
\be
\label{eqcis2}
\left\{\eqalign{
& ||u||_V \le {1\over \alpha}||f||_{V'} + {C_a\over \beta}||g||_{Q'}, \cr
& ||p||_{Q/\Ker B^t} 
\le {C_a\over \beta} \Bigl(||f||_{V'} + {||a||\over \beta}||g||_{Q'}\Bigr).
}\right.
\ee
\finthm

\debdem
Let $u_g\in V$ \st\ $B.u_g=g$, exists since $B$ is surjective,
and $||u_g||_{V/\Ker(B)} \le {1\over \beta}||g||_{Q'}$, \cf~\eref{eqism3}.

Let $u_0\in \Ker(B)$ be the solution of the problem:
Find $u_0\in \Ker(B)$ \st\
\be
\label{equ0}
a(u_0,v_0) = \la f,v_0\ra_{V',V} - a(u_g,v_0), \quad \forall v_0\in \Ker(B).
\ee
The Lax--Milgram theorem tells that~\eref{equ0} is well-posed:
Indeed, $(\Ker B,\dd_V)$ is a Hilbert space,
$a\dd$ is bilinear continuous coercive, and
$F : v_0 \in \Ker(B) \rar F(v_0) := \la f,v_0\ra_{V',V} - a(u_g,v_0)$
is linear (trivial) and continuous on~$\Ker(B)$,
with $||F||_{V'} \le ||f||_{V'} + ||a||\,||u_g||_V$ (easy check).
So $u_0$ exists, is unique, and
$||u_0||_V \le {1\over\alpha} ||F||_{V'}$, that is,
$||u_0||_V
\le {1\over\alpha} (||f||_{V'} + ||a||\,||u_g||_V)
\le {1\over\alpha} (||f||_{V'} + ||a||\,{1\over \beta}||g||_{Q'})
$.

%And $a(u_0,u)  = \la f,u_0\ra_{V',V} - a(u_g,u_0)$,
%gives $\alpha ||u_0||_V^2 \le ||F||_{V'} ||u_0||_V$,
%so $||u_0||_V \le {1\over\alpha} ||F||_{V'}$.

Then let $u:=u_0+u_g$.
So $a(u,v) = \la f,v\ra_{V',V}$, \cf~\eref{equ0}, and $u_0\in \Ker(B)$ and $Bu_g=g$ give
$b(u,q) = b(u_0,q) + b(u_g,q) = 0 + \la g,q\ra_{Q',Q}$,
therefore $u$ is as solution of~\eref{eqw0}.

Moreover $u$ is independent of~$u_g$:
If si $u_g'$ also satisfies $Bu_g'=g$, if $u_0'\in \Ker(B)$ is the associated solution,
if $u'=u_0' + u_g'$, then $u-u' = u_0 - u_0' + u_g - u_g' \in \Ker(B)$ (since $B(u_g - u_g')=g-g=0$)
and $a(u-u',v_0)=0$ for all $v_0\in \Ker(B)$, thus $u-u'=0$ (coercitivy of~$a\dd$ on~$\Ker(B)$), and $u=u'$.
Thus $u=u_0+u_g\in V$ exists and is unique.

And $||u||_V \le ||u_0||_V + ||u_g||_V
\le {1\over\alpha} (||f||_{V'} + ||a||\,{1\over \beta}||g||_{Q'}) + {1\over \beta}||g||_{Q'}$,
that is~\eref{eqcis2}$_1$.

Then we look for $p$ solution of $b(v,p) = a(u,v) - \la f,v\ra_{V',V}$ for all $v\in V$.
Let $L(v) := a(u,v) - \la f,v\ra_{V',V}$.
So if $p$ exists then $L(v) = b(v,p)$, thus $L$~vanishes on~$\Ker(B)$, \ie, $L\in (\Ker(B))^\circ$.
% Let $(\Ker(B))^\circ = \{ m\in V' : \la m,v_0\ra_{V',V} = 0,\; \forall v_0\in \Ker(B)\}$ (the orthogonal of $\Ker(B)$ with respect to the duality).
And $(\Ker(B))^\circ = \overline{\Im(B^t)}$, so $(\Ker(B))^\circ = \Im(B^t)$ (closed ranged theorem~\ref{thmcr}).
Thus there exists $p\in Q$ \st\ $L=B^tp$.
And $||B^tp||_{V'} \ge \beta ||p||_{Q/\Ker B^t}$, \cf~\eref{eqism1}.
Then~\eref{eqcis2}$_2$.
\findem

%%%%%%%%%%%%%%%%%%%%%%%%%%%%%%%%%%%%%%%%%%%%%%%%%%%%%%%%%%%%%%%%%%%%%%%%%%%%%%%%%%%

\subsection{The saddle point problem}

Let $L : V\times Q \rar \RR$ be defined by
\be
\label{eqlag}
\calL(v,q) = \demi a(v,\vv) +  b(v , q) - (f,v)_\Ld - (g,q)_\Ld.
\ee
If $a\dd$ is symmetric then $L$ is the Lagrangean bilinear form associated to the mixed problem~\eref{eqw0}.
And the associated optimization problem is: Find $(u,p)\in V \times Q$ (saddle point) \st\
\be
\label{eqlaggen}
\calL(u,p) = \inf_{v\in V} (\sup_{q\in Q} \calL(v,q)).
\ee
If $(u,p)$ is a solution of~\eref{eqlaggen}, then, $a\dd$ being symmetric,
\be
\left\{\eqalign{
& \forall v \in V,\quad {\pa \calL \over \pa v}(u,p).v 
= \lim_{h\rar0} {\calL(u{+}h v , p) - \calL(u,,) \over h}
= a(u,v) + b(u,q) - \la f,v \ra, \cr
& \forall q \in Q,\quad {\pa \calL \over \pa q}(u,p).q 
= \lim_{h\rar0} {\calL(u , p{+}h q) - \calL(u,p) \over h}
= b(u,q) - \la g,q\ra. \cr
}\right.
\ee
So $(u,p)$ is solution of~\eref{eqw0}.

%%%%%%%%%%%%%%%%%%%%%%%%%%%%%%%%%%%%%%%%%%%%%%%%%%%%%%%%%%%%%%%%%%%%%%%%%%%%%%%%%%%
%%%%%%%%%%%%%%%%%%%%%%%%%%%%%%%%%%%%%%%%%%%%%%%%%%%%%%%%%%%%%%%%%%%%%%%%%%%%%%%%%%%

\section{The surjectivites of the divergence operator}
\label{secdiv}

Let $\Omega$ be an open bounded set in $\RRn$.

Let $b\dd$ be defined by
$b :
\left\{\eqalign{
V\times Q & \rar \RR \cr
(\vv,q) & \rar b(v,q) = \int_\Omega \dvg\vv(x)q(x)\,d\Omega
}\right\}
$
where $V$ and~$Q$ are appropriate Banach spaces, see below ($b\dd$ is bilinear).

Let $B : V \rar Q'$ be the associated operator defined by $\la Bv,q\ra_{Q',Q} = b(v,q)$,
and $B$ will be denoted~$\dvg$ (notation of distribution of L. Schwartz).

Then the operator $B^t : Q \rar V'$ is defined by $\la B^tq,v\ra_{V',V}= \la Bv,q\ra_{Q',Q} = b(v,q)$.

The integration by parts, if legitimate, gives
\be
%\label{eqipp}
b(\vv,q) = \la B\vv,q\ra_{Q',Q} = \la B^tq,\vv\ra_{V',V}= - \int_\Omega dq(x).\vv(x)\, d\Omega
+ \int_\Gamma q(x)\,\vv(x).\vn(x)\,d\Gamma.
\ee

%%%%%%%%%%%%%%%%%%%%%%%%%%%%%%%%%%%%%%%%%%%%%%%%%%%%%%%%%%%%%%%%%%%%%%%%%%%%%%%%%%%

\subsection{The divergence operator $\dvg : \Hdvgo \rar \Ldo$ is surjective}
\label{secdvg1}

Here $b(\vv,q)=(\dvg\vv,q)_\Ld$.

\debthm
\label{thmdsd}
The linear mapping
$\dvg : 
\left\{\eqalign{
\Hdvgo &\rar \Ldo \cr
\vv & \rar \dvg\vv, % = \Tr(d\vv),
}\right\}
$
is continuous and surjective.
And the open mapping theorem gives, \cf~\eref{eq2ai2z},
\be
\label{eqdiv11}
\exists \beta > 0 ,\;
\forall \vv\in \Hdvgo,\quad ||\dvg(\vv)||_\Ldo \ge \beta ||\vv||_{\Hdvg/\Ker(\dvg)},
\ee
or
\be
\label{eqdiv12}
\exists \beta > 0 ,\;
\forall \vv\in \Hdvgo,\; \exists p\in\Ldo,
\quad  (\dvg(\vv),p)_\Ldo \ge \beta ||\vv||_{\Hdvg/\Ker(\dvg)} ||p||_\Ldo,
\ee
also written as the inf-sup inequality 
$\ds \exists \beta>0, \;
\inf_{v\in\Hdvgo}\sup_{p\in\Ldo} {|(\dvg\vv,p)_\Ld| \over ||\vv||_{\Hdvg/\Ker(\dvg)} ||p||_\Ld}  \ge \beta$.
\finthm

\debdem
Since $||\dvg\vv||_\Ld \le ||\vv||_\Hdvg$, $\dvg$ is continuous.
Let $f\in\Ldo$. Let $p\in\Huzo$ be the solution of $\Delta p=f$ (Lax--Milgram theorem).
So $\dvg(\vgrad p) = f \in \Ldo$, thus $\vgrad p \in \Hdvgo$; Then let $\vv = \vgrad p \in\Hdvgo$.
Thus $\dvg\vv=f$, and $\dvg$ is surjective.

And~\eref{eqism3} gives~\eref{eqdiv11}, thus~\eref{eqdiv12}.
\findem

%%%%%%%%%%%%%%%%%%%%%%%%%%%%%%%%%%%%%%%%%%%%%%%%%%%%%%%%%%%%%%%%%%%%%%%%%%%%%%%%%%%

\subsection{The divergence operator $\dvg : \Hdvgzo \rar \Ldzo$ is surjective}
\label{secdvg2}

Here $b(\vv,q)=(\dvg\vv,q)_\Ld$.

\debthm
The linear mapping
$\dvg : 
\left\{\eqalign{
\Hdvgzo &\rar \Ldzo \cr
\vv & \rar \dvg\vv, % = \Tr(d\vv),
}\right\}
$
is continuous and surjective.
And the open mapping theorem gives, \cf~\eref{eq2ai2z},
\be
\label{eqcisdiv30}
\exists \beta > 0 ,\;
\forall \vv\in \Hdvgzo,\quad ||\dvg(\vv)||_\Ldo \ge \beta ||\vv||_{\Hdvgz/\Ker(\dvg)},
\ee
also written as the inf-sup inequality 
$\ds 
\exists \beta>0, \quad
\inf_{p\in\Ldzo}\sup_{\vv\in\Hdvgzo} {|(\dvg\vv,p)_\Ld| \over ||\vv||_{\Hdvgz/\Ker(\dvg)} ||p||_\Ldz}  \ge \beta.
$
\finthm

\debdem
Since $||\dvg\vv||_\Ld \le ||\vv||_\Hdvgzo$, $\dvg$ is continuous.
Let $f\in\Ldzo$. Let $p\in\Huo/\RR$ be the solution of
$(\vgrad p,\vgrad q)_\Ld=(f,q)_\Ld$ for all $q\in\Huo/\RR$, \cf\ the Lax--Milgram Theorem in~$\Huo/\RR$
(the hypothesis $f\in\Ldzo$, that is $(f,1_\Omega)_\Ld=0$ ($=(\vgrad p,\vgrad 1_\Omega)_\Ld$),
is mandatory and is called the compatibility condition).
And $(\vgrad p,\vgrad q)_\Ld=(f,q)_\Ld$ for all $q\in\Huo/\RR$ gives $(\vgrad p.\vn)_{|\Gamma}=0$.
Thus with $\vv=\vgrad p$, we have $\vv\in\Hdvgzo$ and
$\dvg(\vgrad p)= f \in\Ldo$, so $\dvg$ is surjective from $\Hdvgzo$ to~$\Ldo$.
So we get~\eref{eqcisdiv30}, \cf~\eref{eqism3}.
\findem

%%%%%%%%%%%%%%%%%%%%%%%%%%%%%%%%%%%%%%%%%%%%%%%%%%%%%%%%%%%%%%%%%%%%%%%%%%%%%%%%%%%

\subsection{The divergence operator $\dvg : \Ldon \rar \Hmuo$ is surjective}
\label{secdvg3}

Here $b(\vv,q)= \la\dvg\vv , q\ra_{\Hmu,\Huz} = -\la \vv,dq\ra_{\Ldo,\Ldo} := -\int_\Omega dq(x).\vv(x)\,d\Omega$ (distributions of L. Schwartz) for all $q\in\Huzo$.

\debthm
\label{thmdvghmu}
The linear mapping
$\dvg :
\left\{\eqalign{
\Ldon &\rar \Hmuo\cr
\vv & \rar \dvg\vv,
}\right\}
$
is continuous and surjective.
And the open mapping theorem gives, \cf~\eref{eq2ai2z},
\be
\label{eqcisdiv3b}
\exists \beta > 0 ,\;
\forall \vv\in \Ldo,\quad ||\dvg(\vv)||_\Hmu \ge \beta ||\vv||_{\Ldo/\Ker(\dvg)},
\ee
also written as the inf-sup inequality 
$\ds 
\exists \beta>0, \quad
\inf_{p\in\Huzo}\sup_{\vv\in\Ldo} {|b(\vv,q)| \over ||\vv||_{\Ldo/\Ker(\dvg)} ||p||_\Huz}  \ge \beta.
$
\finthm

\debdem
$||\dvg\vu||_{\Hmu}
= \sup_{\phi\in\Huzo} {|\la \dvg\vu,\phi| \ra \over ||\phi||_\Huz}
= \sup_{\phi\in\Huzo} {|(\vu,\vgrad\phi)_\Ld \over ||\phi||_\Huz}
\le ||\vu||_\Ld
$ (Cauchy--Schwarz in~$\Ldo$), therefore $\dvg$ is continuous.
Let $\ell\in\Hmuo$. Thus there exists $f\in \Ldo$ and $\vu\in\Ldon$ \st\ $\ell = f + \dvg\vu$, \cf~\eref{eqdhuzo}.
Let $\vw\in\Hdvgo$ \st\ $\dvg\vw = f$, \cf\ thm.~\ref{thmdsd}.
So $\ell = \dvg(\vw+\vu)$ with $\vu+\vw\in\Ldon$, and $\dvg$ is continuous.
So we get~\eref{eqcisdiv3b}, \cf~\eref{eqism3}.
\findem

%%%%%%%%%%%%%%%%%%%%%%%%%%%%%%%%%%%%%%%%%%%%%%%%%%%%%%%%%%%%%%%%%%%%%%%%%%%%%%%%%%%

\subsection{The divergence operator $\dvg : \Huzon \rar \Ldzo$ is surjective}
\label{secdvg4}

Here $b(\vv,q)=(\dvg\vv,q)_\Ld$.

\debthm
The linear mapping
\be
\label{eqdivhuz}
\dvg : 
\left\{\eqalign{
\Huzon &\rar \Ldzo \cr
\vv  & \rar \dvg\vv , % = \Tr(d\vv),
}\right\}\quad\hbox{is continuous and surjective}.
\ee
And the open mapping theorem gives, \cf~\eref{eq2ai2z},
\be
\label{eqcisdiv}
\exists \beta > 0 ,\;
\forall \vv\in \Huzon,\quad ||\dvg(\vv)||_\Ldz \ge \beta ||\vv||_{\Huzon/\Ker(\dvg)},
\ee
also written as the inf-sup inequality 
$\ds
%\label{eqcisdiv2}
\exists \beta>0,\quad \inf_{p\in\Ldzo}\sup_{\vv\in\Huzon}
{|(\dvg\vv,q)_\Ld| \over ||\vv||_{\Huz/\Ker(\dvg)} ||p||_\Ldz)_\Ld} \ge \beta.
$
\finthm

\debdem
$\dvg = \vgrad' : \Huzo \rar \Ldzo$ is the dual operator of the gradient operator $\vgrad : \Ldo \rar \Hmuon$.
Since the range of the $\vgrad$ is closed, \cf\ theorem~\ref{thmrgo},
the range of the $\dvg$ operator is closed, \cf\ the closed range theorem~\ref{thmcr}.
(Remark:
For any $\vv\in\Huzon$ we have
$\int_\Omega \dvg\vv\,d\Omega = \int_\Gamma \vv.\vn\,d\Gamma = 0$, so $\Im(\dvg) \subset \Ldzo$.)
\findem

With $\dvg:\Huzon \rar \Ldo$ we have $\Ker(\dvg) = \{\vv\in\Huzon : \dvg\vv=0\}$, and
\be
\Ker(\dvg)^{\perp_\Huz}
:= \{\vv \in \Huzon : (\vv,\vw)_\Huz = 0,\; \forall \vw\in\Huzon,\; \dvg\vw=0\}.
\ee
Let $\Delta^{-1} :
\left\{\eqalign{
\Hmuo & \rar \Huzo \cr
f & \rar u = \Delta^{-1}f
}\right\}
$, that is, $u \in \Huzo$ solves the Dirichlet problem $\Delta u = f$.

\debcor
\be
\Ker(\dvg)^{\perp_\Huz}
= \{ \vv = \Delta^{-1}(\vgrad q),\; q\in\Ldo\} \quad (= \Delta^{-1}(\vgrad(\Ldo)) ),
\ee
that is $\vv\in \Ker(\dvg)^{\perp_\Huz}$  iff $\Delta\vv$ derives from a potential $q\in\Ldo$.

And $\Huzon = \Ker(\dvg) \oplus^{\perp_\Huz} \Ker(\dvg)^{\perp_\Huz}$ give a decomposition of~$\Huzon$.
\fincor

\debdem
Let $A := \{\vv\in \Huzon : \vv = \Delta^{-1}(\vgrad q),\; q\in\Ldo\}$.
So $\vv\in A$ iff $\vv\in \Huzon$ and $\exists q\in \Ldo$, $\Delta\vv = \vgrad q$,
i.e. $(\vv,\vw)_\Huz = (\grad \vv,\grad\vw)_\Ld = (q,\dvg\vw)_\Ld$ for all $\vw\in\Huzon$.

$\bullet$ $A \subset \Ker(\dvg)^{\perp_\Huz}$:
Let $\vv\in A$. Thus $\exists q\in\Ldo$ \st\
$(\vv,\vw)_\Huz = (q,\dvg\vw)_\Ld$ for all $\vw\in\Huzon$.
Thus $(\vv,\vw)_\Huz =0$ for all $\vw\in\Ker(\dvg)$,
thus $\vv \in \Ker(\dvg)^{\perp_\Huz}$.

$\bullet$ $\Ker(\dvg)^{\perp_\Huz} \subset A$:
Let $\vv\in \Ker(\dvg)^{\perp_\Huz}$.
%, so $\vv\in\Huzon$ satisfies $(\grad\vv,\vgrad\vw_0)_\Ld = 0 % = \la \Delta \vv,\vw\ra_{\Hmu,\Huz} $ for all $\vw_0\in\Ker(\dvg)$.
We look for $q\in\Ldzo$ \st\ $\Delta \vv = \vgrad q$:
thus we look for $q\in\Ldzo$ \st\ $\Delta \vv = \vgrad q$,
that is $(q,\dvg\vz)_\Ld = -(\Delta\vv,\vz)_{\Hmu,\Huz}$ for all $\vz\in\Huzo$.

The operator $B=\dvg : \vz\in \Huzo/\Ker(\dvg) \rar \dvg\vz \in \Ldzo$ is linear continuous bijective,
with $||\dvg\vz||_\Ldzo \le ||B||\,||\vz||_{ \Huzo/\Ker(\dvg)}$.

Its inverse $B^{-1} : \psi\in\Ldzo \rar B^{-1}\psi\in \Huzo/\Ker(\dvg)$ is linear continuous bijective,
with $||B^{-1}\psi||_{\Huzo/\Ker(\dvg)} \le ||B^{-1}||\,||\psi||_\Ldz$.

%So we have to find $q\in \Ldzo$ s.t.$(q,\dvg\vz)_\Ld = -(\Delta\vv,\vz)_{\Hmu,\Huz}$ for all $\vz\in\Huzo$.

Let $a\dd : (q,\psi)\in \Ldzo\times \Ldzo \rar a(q,\psi) = (q,\psi)_\Ld$:
$a\dd$ is trivially bilinear continuous coercive in $(\Ldo,\dd_\Ld)$.

Let $\ell : \psi \in\Ldzo \rar \ell(\psi)= -(\Delta\vv,B^{-1}\psi)_{\Hmu,\Huz} = (\grad\vv,\grad\psi)_\Ld 
\in \RR$: $\ell$ is trivially linear,
and is continuous since $|\ell(\psi)|
\le ||\Delta\vv||_\Hmu ||B^{-1}\psi||_\Ld
\le ||\Delta\vv||_\Hmu||B^{-1}||\,||\psi||_\Ldz$.

Thus the Lax-Milgram theorem gives the existence of~$q$. And the first point shows
that if $\Delta\vv= \vgrad q$ then $\vv \perp_\Huz \Ker(\dvg)$.
\findem

%%%%%%%%%%%%%%%%%%%%%%%%%%%%%%%%%%%%%%%%%%%%%%%%%%%%%%%%%%%%%%%
%%%%%%%%%%%%%%%%%%%%%%%%%%%%%%%%%%%%%%%%%%%%%%%%%%%%%%%%%%%%%%%

\newpage
\setcounter{section}{0}
\def\thesection{\Alph{section}}

%%%%%%%%%%%%%%%%%%%%%%%%%%%%%%%%%%%%%%%%%%%%%%%%%%%%%%%%%%%%%%%
%%%%%%%%%%%%%%%%%%%%%%%%%%%%%%%%%%%%%%%%%%%%%%%%%%%%%%%%%%%%%%%

%\newpage
\section{Singular value decomposition (SVD)}

We want to estimate the $\beta$ inf-sup constant, \cf~\eref{eqism4}.
%(in fact the $\beta_h$ numerical inf-sup constant).
Consider a $m*n$ rectangular matrix~$B$.
We look for its singular value $\sigma_i$, that is, we look for
a $m*n$ ``diagonal'' matrix~$\sigma$, \ie\ \eg\ in the case $m<n$,
$$
  \Sigma=\diag_{m,n}(\sigma_1,...\sigma_p)=\pmatrix{
\sigma_1&0   &... &  &0 &0& \ldots & 0\cr
0 & \sigma_1 & 0  &  &0&\vdots & &\vdots\cr
\vdots & & \ddots & \ddots \cr
  & &       &  & 0      \cr
0 & &\ldots &  & \sigma_p & 0&\ldots&0\cr},
$$
and for two matrices $U$~$m*m$ and $V$~$n*n$ \st
$$
\Sigma=U^T.B.V,
$$
with $U^T$ the $U$ transposed matrix.

\debprop
Let $B$ be a $m*n$ real matrix.
If $\lambda_i$ is an eigenvalue of the $n*n$ matrix $B^T.B$,
then $\lambda_i$ is positive and is an eigenvalue of the $m*m$ matrix $B.B^T$.

If $\lambda_i$ is an eigenvalue of the $m*m$ matrix $B.B^T$,
then $\lambda_i$ is positive and is an eigenvalue of the $n*n$ matrix $B^T.B$.

Let $\sigma_i = \sqrt{\lambda_i}$.
Let $(\vv_i)_{1,...,n}$ be an orthonormal basis of eigenvectors of~$B^T.B$ associated to the eigenvalues~$\lambda_i$, and let $V$ be the (orthonormal) matrix whose $j$-th column is~$\vv_j$.
Let $(\vu_i)_{1,...,n}$ be an orthonormal basis of eigenvectors of~$B.B^T$ associated to the eigenvalues~$\lambda_i$, and let $U$ be the (orthonormal) matrix whose $j$-th column is~$\vu_j$.
And let $\Sigma= \diag_{m,n}(\sigma_1,...,\sigma_p)$ where $p=\min(m,n)$.
Then the singular value decomposition of~$B$ is
\be
\label{eqsutav}
  \Sigma=U^T.B.V,\qie B=U.\Sigma^T.V^T.
\ee
Thus, if ${\rm rank}(B)=r$ and $\sigma_1\ge...\ge\sigma_r>0$ (and $\sigma_i=0$ pour $i>r$), then
\be
\label{eqsutav2}
   B=\sum_{i=1}^r \sigma_i \,\vu_i.\vv_i^T.
\ee
Moreover,
$\pmatrix{\vu_j\cr\vv_j}\in \RR^{m+n}$, $j=1,...,p$, is an eigenvector of
$\pmatrix{0&B\cr B^T&0}$ associated to the eigenvalue~$\sigma_j$. 
%Et si $B$ est une matrice carré $n*n$ symétrique positive (réelle), alors les $\sigma_i$ sont ses valeurs propres, et $U=V$.
\finprop

\debdem
$B^T.B$ is symmetric real, thus diagonalisable.
Moreover $B^T.B$ is non negative since $\vx^T.(B^T.B).\vx=(B.\vx)^T.(B.\vx)=||B.\vx||^2\ge0$.
Let $\lambda_1,...,\lambda_n$ be its eigenvalues,
and $\lambda_1\ge...\ge\lambda_n(\ge0)$, even if you have to renumber them.
Let $\vv_i$ be associated eigenvectors constituting an orthonormal basis in~$\RRn$,
and let $V$ be the orthonormal matrix which columns are made of the~$\vv_i$'s.
Suppose $(B)\le p=\min(m,n)$, so that ${\rm rank}(B^T.B)\le p$ and $\lambda_{p+1}=...=\lambda_n=0$.
Then
$$
   \diag_{n,n}(\lambda_1,...,\lambda_p,0,...,0)=V^T.B^T.B.V\quad\hbox{$n*n$ matrix}.
$$
Same steps for $B.B^T$ with eigenvalues $\mu_1\ge..\ge\mu_m\ge0$
and the associated orthonormal matrix~$U$:
$$
   \diag_{m,m}(\mu_1,...,\mu_p,0,...,0)=U^T.B.B^T.U\quad\hbox{$m*m$ matrix}.
$$
With $B.B^T.\vu_i=\mu_i\vu_i$ we get $B^T.B.B^T.\vu_i=\mu_iB^T.\vu_i$,
thus  $B^T.\vu_i$ is an eigenvector for $B^T.B$ associated to the eigenvalue~$\mu_i$.
And $B^T.B.\vv_j=\lambda_j\vv_j$ tells that $\mu_i$ is one the the~$\lambda_j$.

Remark: If $\Sigma=\diag_{m,n}(\sigma_1,...,\sigma_p)=U^T.B.V$ then
$\Sigma.\Sigma^T=\diag_{m,m}(\sigma_1^2,...,\sigma_p^2)=
(U^T.B.V).(U^T.B.V)^T=U^T.(B.B^T).U$,
and the $\lambda_i=\sigma_i^2$ are indeed the eigenvalues of~$B.B^T$
associated to the eigenvectors of~$U$.
Idem for $B^T.B$.
And the matrices $U$ and $V$ are the matrices made of the column vectors $\vu_j$ and~$\vv_j$.

Existence of the decomposition:
Let $\lambda_i$, $i=1,...,n$, be the eigenvalues of $B^T.B$.
Suppose $\lambda_1\ge...\ge\lambda_r>0$, and $\lambda_{r+1}=...=\lambda_n=0$.
Let $\sigma_j=\sqrt{\lambda_j}$. 

Then let
\be
\label{eqvujbvj}
\vu_j={B.\vv_j\over\sigma_j}\in\RR^m,\quad 1\le j\le r.
\ee
The~$\vu_j$ are (orthonormal) eigenvectors of~$B.B^T$:
Indeed
$(B.B^T).\vu_j={B.(B^T.B).\vv_j\over\sigma_j}={\lambda_j B.\vv_j\over\sigma_j}=\lambda_j\vu_j$.
And $\vu_i^T.\vu_j={\vv_i^T.(B^T.B\vv_j)\over\sigma_i\sigma_j}
=\lambda_j{\vv_i^T.\vv_j\over\sigma_i\sigma_j}=\delta_{ij}$,
and the~$\vu_j$ are the normalized vectors~$B.\vv_j$.
We then complete $(\vu_j)_{j=1,..,r}$ to get an orthonormal basis in~$\RRm$
(\eg\ with Gram--Schmidt method).
Let $U$ be the $m*m$ matrix made of the columns vector~$\vu_j$.

Let $\Sigma=U^T.B.V$. So $[\Sigma_{ij}]=[\vu_i^T.B.\vv_j]
=[\sigma_j\vu_i^T.\vu_j]=[\sigma_j\delta_{ij}]$
if $j\le r$, and vanishes if $j>r$ (since $B.\vv_j=0$).
Thus~\eref{eqsutav}. And then~\eref{eqsutav2}.

And $B\vv_j=\sigma_j\vu_j$ for $j\le r$, \cf~\eref{eqvujbvj}, thus
$B^T.B.\vv_j=\sigma_j B^T.\vu_j=\lambda_j\vv_j$ for $j\le r$,
so $B^T.\vu_j=\sigma_j\vv_j$ for $j\le r$.
And if $j>r$ then $B^T.\vu_j=0$ since $\vu_j\in(\Im(B))^\perp=\Ker(B^T)$.
Thus $\pmatrix{0&B\cr B^T&0}.\pmatrix{\vu_j\cr\vv_j}=\sigma_j\pmatrix{\vu_j\cr\vv_j}$.

And if $B$ is a symmetric positive real matrix, then $B^T.B=B^2=B.B^T$,
so $B^T.B=B.B^T$; With $B\ge0$ thus its eigenvalues are non negative,
$\sigma_i=+\sqrt{\lambda_i}$, thus the $\sigma_i$ are the singular values.
\findem

\debrem
If $m>n$ and $j\ge m+1$ then the $\vu_j$ are useless, and $U$ is computed as a $m*n$ matrix
(and $U^T$ as a $n*m$ matrix). And $\Sigma$ is then a $n*n$ matrix.
This method is called the ``Thin SVD''.
\finrem

\debcor
$\rang(B)=r$, $\Ker(B)=\Vect\{\vv_{r+1},...\vv_n\}$ and $\Im(B)=\Vect\{\vu_1,...\vu_r\}$.
\fincor

\debdem
Apply~\eref{eqsutav2}.
\findem

\debcor
Let $k\le r{-}1$ and $B_k=\sum_{i=1}^k \sigma_i\vu_i.\vv_i^T$.
Then
$$
\ds\min_{Z : \rang Z=k}||B-Z||=\sigma_{k+1}=||B-B_k||,
$$
where
$||Z||=\sup_{\vx\ne\vec0}{||Z.\vx||_\RRm\over||\vx||_\RRn}$ is the usual norm.

This gives a numerical measure of the rank of~$B$:
If $\sigma_{k+1}$ is of the precision order of the computer,
then the numerical rank o~$B$ is~$k$.
\fincor

\debdem
We get $U^T.B_k.V=\diag_{m,n}(\sigma_1,...,\sigma_k,0,...)$ (easy check).

Thus $U^T.(B-B_k).V=\diag_{m,n}(0,...,0,\sigma_{k+1},...,\sigma_r,0,...)$,
thus $||B-B_k||=\sigma_{k+1}$, and in particular $\ds\min_{Z : \rang Z=k}||B-Z||\le \sigma_{k+1}$.

Let $Z$ be a $m*n$ matrix with rank~$k$. Thus $\dim\Ker Z=n-k$.
Let $E=\Ker Z\bigcap\Vect\{\vv_1,...\vv_{k+1}\}$. So $\dim(E)\ge 1$
(intersection of a dimension $n{-}k$ space with a dimension $k{+}1$ space in~$\RR^n$).

Let $\vx\in E$ \st\ $||\vx||_\RRn=1$; Then
$||(B{-}Z).\vx||_\RRm^2=||B.\vx||_\RRm^2
=||\sum_{i=1}^{k+1} \sigma_i(\vv_i^T.\vx)\vu_i||^2
=\sum_{i=1}^{k+1} \sigma_i^2(\vv_i^T.\vx)^2$,
the $\vu_i$ being orthonormal vectors. Thus
$||(B{-}Z).\vx||_\RRm^2\ge \sigma_{k+1}\sum_{i=1}^{k+1}(\vv_i^T.\vx)^2=\sigma_{k+1}||\vx||^2
=\sigma_{k+1}$. So
$\ds\min_{Z : \rang Z=k}||B-Z||\ge \sigma_{k+1}$.
\findem

%%%%%%%%%%%%%%%%%%%%%%%%%%%%%%%%%%%%%%%%%%%%%%%%%%%%%%%%%%%%%%%

\section{Application: The discrete inf-sup condition}

Example of $\dvg\vu=0$, corresponding to $B$ a rectangular $m*n$ matrix
(computation of $b(\vv_h,q_h)=0$).

Here $B^T$ stands for~$[B_h]$, \cf~\eref{eqgenh0m},
and we compute the singular values of~$B^T$,
i.e. the eigenvalues of $\pmatrix{0&B^T\cr B&0}$. We get:

\debprop
Let $\sigma_r>0$ be the smallest positive eigenvalue of~$B$.
We have (value of the inf-sup constant)
$$
  \inf_{q_h\in Q_h}\sup_{v_h\in V_h}{b(v_h,q_h)\over||v_h||_V||q_h||_Q}=\sigma_r.
$$
\finprop

\debdem
$B=\sum_{i=1}^r \sigma_i\vu_i.\vv_i^T$ gives
$B.\vx=\sum_{i=1}^r \sigma_i(\vv_i^T.\vx)\,\vu_i$, so
$\vy^T.B.\vx=\sum_{i=1}^r \sigma_i(\vv_i^T.\vx)\,(\vy^T.\vu_i)$.

Let $\vx=\sum_{j=1}^nx^j\vv_j$
and $\vy=\sum_{k=1}^ny^k\vu_k$. Thus
$\vy^T.B.\vx=\sum_{i=1}^r \sigma_ix^i\,y^i$.

Thus for $||x||=1$, and $\vy$ being fixed with $||\vy||=1$,
the sup is given by $x^i={\sigma_iy^i\over(\sum_i\sigma_i^2y_i^2)^\demi}$,
and gives $\vy^T.B.\vx={1\over(\sum_i\sigma_i^2y_i^2)^\demi}\sum_{i=1}^r \sigma_i^2y_i^2
=(\sum_i\sigma_i^2y_i^2)^\demi$.
Thus the inf for~$\vy$ is given with $\vy=\vu_r$, and gives $\vy^T.B.\vx=\sigma_r$.
\findem

%%%%%%%%%%%%%%%%%%%%%%%%%%%%%%%%%%%%%%%%%%%%%%%%%%%%%%%%%%%%%%%%%%%%%%%%%%%%%%%%%%%%%%%%

%\newpage

 \catcode`@=11
\def\thebibliography#1{\section*{References \@mkboth
 {References}{References}}\list
% {[\arabic{enumi}]}{\settowidth\labelwidth{[#1]}\leftmargin\labelwidth
 {\arabic{enumi}}{\settowidth\labelwidth{#1 }\leftmargin\labelwidth
 \advance\leftmargin\labelsep
 \usecounter{enumi}}
 \def\newblock{\hskip .11em plus .33em minus .07em}
 \sloppy\clubpenalty4000\widowpenalty4000
 \sfcode`\.=1000\relax}
\def\@cite#1#2{[{#1\if@tempswa , #2\fi}]}
\def\@biblabel#1{[#1]\hfill}
\def\@bibitem#1{\item\if@filesw \immediate\write\@auxout
       {\string\bibcite{#1}{\the\c@enumi}}\fi\ignorespaces}
 \catcode`@=12
%}

\end{document}